\documentclass[english,12pt]{article}  
 
\usepackage{cite}             
\usepackage{latexsym}         
\usepackage{amsmath}          
\usepackage{amssymb}          
\usepackage{amsxtra}          
\usepackage[latin1]{inputenc} 
\usepackage{eucal}            
\usepackage{longtable}        
\usepackage{exscale}          
\usepackage{bbm}              
\usepackage{theorem}          
\usepackage{paralist}         
\usepackage{stmaryrd}         
\usepackage{a4wide}
\usepackage{float}
\usepackage{url}
\usepackage{enumitem}
\usepackage{tensor}
\usepackage{tikz}             
\usepackage{tikz-cd}          
\usetikzlibrary{matrix, arrows, shapes, plotmarks, decorations.pathreplacing, calligraphy}
\usepackage{leftidx}          
\usepackage{rotating}
\usepackage{hyperref}

\numberwithin{equation}{section}

\newcommand{\I}          {\mathrm{i}}
\newcommand{\C}        {\mathbb{C}}
\newcommand{\N}        {\mathbb{N}}
\newcommand{\Z}        {\mathbb{Z}}

\newcommand{\E}          {\mathrm{e}}
\newcommand{\cc}[1]      {\overline{{#1}}}

\newcommand{\Vect}   {\mathsf{Vect}}
\newcommand{\vect}   {\mathsf{vect}}

\newcommand{\tr}         {\operatorname{\mathsf{tr}}}
\newcommand{\chara}         {\mathrm{char}}
\newcommand{\Ad}{\mathrm Ad}
\newcommand{\FS}{\mathrm{FS}}

\newcommand{\PQ}{\mathrm{PQ}}
\newcommand{\PDQ}{\mathrm{P(DQ)}}
\newcommand{\CentQ}{\mathrm{CentQ}}
\newcommand{\CQ} {\mathrm{CQ}}  

\newcommand{\opp}             {{\mathrm{op}}}

\newcommand{\vses}             {{\mathrm{v-ses}}} 

\def\boti  {\boxtimes} 

\newcommand{\circtensor} {\setbox0\hbox{\large$\circlearrowleft$} \rlap{\hbox to\wd0{\hss$\times$\hss}}\box0}
\newcommand{\circtensorsmall} {\setbox0\hbox{$\circlearrowleft$} \rlap{\hbox to\wd0{\hss$\times$\hss}}\box0}

\newcommand{\id}         {\operatorname{\mathsf{id}}}

\newcommand{\unit}         {\operatorname{\mathbbm{1}}}

\newcommand{\act}        {\operatorname{\triangleright}}

\newcommand{\ract}       {\operatorname{\triangleleft}}

\newcommand{\Cat}[1]         {\operatorname{\mathcal{#1}}}

\newcommand{\op}[1]           {\Cat{#1}^{\operatorname{\mathrm{op}}}}
\newcommand{\rev}           {{\operatorname{\mathrm{rev}}}} 

\newcommand{\Fun}         {\operatorname{\mathrm{Fun}}}

\newcommand{\twoFunls}         {\operatorname{\mathrm{2Fun}^{\mathrm{lax,str}}}}
\newcommand{\twoFunol}         {\operatorname{\mathrm{2Fun}^{\mathrm{oplax}}}}

\newcommand{\Hom}        {\operatorname{\mathsf{Hom}}}

\newcommand{\End}        {\operatorname{\mathsf{End}}}
\newcommand{\Gl}        {\operatorname{\mathsf{Gl}}}
\newcommand{\Gr}        {\operatorname{\mathsf{Gr}}}
\newcommand{\Sl}        {\operatorname{\mathsf{Sl}}}
\newcommand{\SL}        {\operatorname{\mathsf{SL}}}

\newcommand{\im}         {\operatorname{\mathsf{im}}}

\newcommand{\Lincat}         {\operatorname{\mathrm{LinCat}}}

\newcommand{\lincatses}         {\operatorname{\mathrm{linCat}^{\mathrm{ses}}}}
\newcommand{\Mses}         {\operatorname{\mathcal{M}}^{\mathrm{ses}}} 

\newcommand{\Bimod}[5]{\sideset{^{\scriptscriptstyle{#1}}_{\scriptscriptstyle{#2}}}{^{\scriptscriptstyle{#4}}_{\scriptscriptstyle{#5}}}{\operatorname{#3}}}

\newcommand{\AM}  {\Bimod{}{\Cat{A}}{\Cat{M}}{}{}}
\newcommand{\AMp}  {\Bimod{}{\Cat{A}}{\Cat{M}}{'}{}}
\newcommand{\AMpp}  {\Bimod{}{\Cat{A}}{\Cat{M}}{''}{}}

\newcommand{\AN}  {\Bimod{}{\Cat{A}}{\Cat{N}}{}{}}

\newcommand{\AAA}  {\Bimod{}{\Cat{A}}{\Cat{A}}{}{\Cat{A}}}

\newcommand{\Mod}      {\operatorname{\mathsf{Mod}}}
\newcommand{\Modoplax}      {\operatorname{\mathsf{Mod}^{\mathrm{oplax}}}}

\newcommand{\BiModoplax}      {\operatorname{\mathsf{Bimod}^{\mathrm{oplax}}}}
\renewcommand{\mod}      {\operatorname{\mathsf{mod}}}
\newcommand{\comod}      {\operatorname{\mathsf{comod}}}
\newcommand{\modln}      {\operatorname{\mathsf{mod}^{ln}}}
\newcommand{\Clnz}      {\operatorname{\mathcal{C}^{ln}_{z}}}

\newcommand{\Modln}      {\operatorname{\mathsf{Mod}^{ln}}}

\newcommand{\Quiv}      {\operatorname{\mathsf{Quiv}}}
\newcommand{\fusqui}      {\operatorname{\mathsf{FusQuiv}}}
\newcommand{\Dfusqui}      {\operatorname{\mathsf{DrinFusQuiv}}}

\newcommand{\Rep}      {\operatorname{\mathsf{Rep}}}
\newcommand{\Comod} {\operatorname{\mathsf{Comod}}}

\newcommand{\cent}     {\mathcal{Z}}
\newcommand{\Cent}     {\mathrm{Cent}} 
\newcommand{\PCent}     {\mathrm{PCent}} 

\newcommand{\ev}[1]   {\operatorname{\mathsf{ev}}_{#1}}
\newcommand{\coev}[1]   {\operatorname{\mathsf{coev}}_{#1}}

\newcommand{\qoti}    {\oast}    
\newcommand{\otik}    {\otimes_{\Bbbk}}    

\makeatletter
\newcommand{\Rmnum}[1]{\expandafter\@slowromancap\romannumeral #1@}
\makeatother

\theoremheaderfont{\normalfont\bfseries}
\theorembodyfont{\itshape}
\newtheorem{lemma}{Lemma}[section]
\newtheorem{proposition}[lemma]{Proposition}
\newtheorem{theorem}[lemma]{Theorem}
\newtheorem{corollary}[lemma]{Corollary}
\newtheorem{definition}[lemma]{Definition}

\newenvironment{theorem-n}[1][Theorem]{\begin{trivlist}
  \item[\hskip \labelsep {\bfseries #1}]\itshape }{\end{trivlist}}

\theorembodyfont{\rmfamily}

\newtheorem{example}[lemma]{Example}
\newtheorem{remark}[lemma]{Remark}

\newenvironment{theorem-1}[1][Theorem 1]{\begin{trivlist}
  \item[\hskip \labelsep {\bfseries #1}]}{\end{trivlist}}
\newenvironment{theorem-2}[1][Theorem 2]{\begin{trivlist}
  \item[\hskip \labelsep {\bfseries #1}]}{\end{trivlist}}
\newenvironment{theorem-3}[1][Theorem 3]{\begin{trivlist}
  \item[\hskip \labelsep {\bfseries #1}]}{\end{trivlist}}
\newenvironment{cor-n}[1][Corollary]{\begin{trivlist}
  \item[\hskip \labelsep {\bfseries #1}]}{\end{trivlist}}



\newcommand{\refitem}[1] {~\textit{\ref{#1})}}

\makeatletter
\newcommand\qedsymbol{\hbox{$\boxempty$}}

\newcommand\qed{\relax\ifmmode\boxempty\else
  {\unskip\nobreak\hfil\penalty50\hskip1em\null\nobreak\hfil\qedsymbol
    \parfillskip=\z@\finalhyphendemerits=0\endgraf}\fi}
\makeatother
\newenvironment{proof}[1][{}]{\par\noindent Proof{#1}. }{\qed}
\newenvironment{theoremlist}{\begin{enumerate}}{\end{enumerate}}
\newenvironment{remarklist}{\begin{enumerate}}{\end{enumerate}} 
\newenvironment{examplelist}{\begin{enumerate}}{\end{enumerate}}
\newenvironment{lemmalist}{\begin{enumerate}}{\end{enumerate}}
\newenvironment{propositionlist}{\begin{enumerate}}{\end{enumerate}}
\newenvironment{definitionlist}{\begin{enumerate}}{\end{enumerate}}

\usepackage{color}

\definecolor{DarkViolet} {rgb}{0.580392,0.000000,0.827450}

\tikzset{
  doublearrow/.style={draw, thin, double distance=3pt, ->, >=implies},
  ldoublearrow/.style={draw, thin, double distance=3pt, <-, >=implies},
  thirdline/.style={draw, thin, ->, >=implies}, 
  lthirdline/.style={draw, thin, <-, >=implies},
  equality/.style={draw, thin, double distance=3pt, -},
  shorten <>/.style={shorten >=#1,shorten <=#1}
}
\makeatletter
\def\latearrow#1#2#3#4{%
  \toks@\expandafter{\tikzcd@savedpaths\path[/tikz/commutative diagrams/every arrow,#1]}%
  \global\edef\tikzcd@savedpaths{%
    \the\toks@%
    (\tikzmatrixname-#2)
    to%
    node[/tikz/commutative diagrams/every label] {$#4$}
    (\tikzmatrixname-#3)
    ;}}
\makeatother

\begin{document}
   

\vskip 22mm{}

\begin{center}

  {\Large\bf Fusion Quivers}

  \vskip 12mm{}

  {\large \  \ Gregor Schaumann$^{*}\,$
  }
\vskip 6mm
{\today}
\vskip 12mm

 \it$^*$ 
 Mathematische Physik, \ Institut f\"ur Mathematik, \\
  Universit\"at W\"urzburg \\
 Emil-Fi\-scher-Stra\ss e 31, \ D\,--\,97\,074\, W\"urzburg

\end{center}

\vskip 3.2em

\noindent{\sc Abstract}\\[3pt]
We develop a categorical approach to quivers and their modules. Naturally this leads to a notion of an action of
a monoidal category on quivers. Using this, we construct for a large class of quivers rigid monoidal structures on
their categories of modules. This fusion product on the quiver modules induces a graded ring structure with duality and trace on the moduli spaces
of semisimple quiver modules. Our approach allows to consider a class of relations on such fusion quivers that are compatible with the rigid monoidal structure.
In particular we obtain a class of preprojective algebras with fusion product on their modules. 

\newpage


\tableofcontents
\newpage




\section{Introduction}

The term ``fusion quiver'' indicates that these objects lie in the intersection of two mathematical areas: quivers and fusion categories.

Quivers, i.e. finite directed graphs, are mathematical objects of great interest far beyond graph theory and combinatorics.
Since the pioneering works \cite{GabrielQ} and  \cite{KacQ},
the main interest lies in the modules $\mod(Q)$ over a given quiver $Q$, see \cite{NakaV}, \cite{LusQuiv}, \cite{reineke08}.
A main theme in this line of research is that certain  equivalence classes of quiver modules form interesting geometric objects.
An important tool thereby is the preprojective algebra of a quiver $Q$ which is obtained as a certain quotient of the double of $Q$.

The other relevant mathematical object, fusion categories
(see e.g. \cite{EGNObook}) are more complicated and harder to construct, but albeit of their categorical nature still concrete: A fusion category $\Cat{A}$ is a finite semisimple category with monoidal product, such that each object has a left and right dual (such a category is called rigid monoidal). By semisimplicity it suffices to know these structures on the simple objects, which implies that a fusion category contains finitely many data.
Fusion categories arise as categories of modules over semisimple Hopf algebras \cite{ENOfus}, certain vertex operator algebras (see \cite{HL1} and other works in their series), and, especially when equipped with a braiding,  have applications to constructions in topological and conformal field theories \cite{Resh,FFRS-fact}. A main tool in studying these categories and for the applications is the graphical calculus for (braided) monoidal categories.
The central questions of this article are: 

\vspace{1em}
\emph{Which quivers $Q$ allow for a rigid monoidal structure on $\mod(Q)$? How to construct such structures?}

\vspace{1em}
We believe that these questions are of interest for both communities: When studying quivers, one usually focuses on a subclass of all quivers, which are of interest due to their algebraic or geometric properties.
In \cite{CRHopf} for instance, Hopf quivers are studied, whose path algebra has the structure of a Hopf algebra.

This article defines fusion quivers as quivers with rigid monoidal structure on their categories of modules and studies their  algebraic and geometric properties. They generalize Hopf quivers  and the  moduli spaces  of semisimple modules acquire an associative multiplication.

For the theory of fusion categories the questions are relevant for the following reasons. 
A main problem in studying fusion categories is their rigidity in the sense that they have no deformations as semisimple monoidal categories. This makes it difficult to construct examples and to classify them.
Fusion quivers provide in the following sense a non-semisimple thickening of a fusion category:
For a fusion quiver $Q$ we demand that the modules in $\mod(Q)$ which all arrows acting as zero, form a rigid monoidal subcategory-- which is thus a fusion category
  and a retract of $\mod(Q)$,  see  Equation \ref{eq:retract-intro}. 
Thus, a  fusion quiver allows for bridges between semisimple and non-semisimple rigid monoidal categories and 
  might also have applications to field theories. Furthermore given a fusion quiver one can try to identify relations on the quiver that are compatible with the rigid monoidal structure to construct further interesting rigid monoidal categories that extend the fusion category in a certain sense. 

  We continue to summarize the main contributions of this work. 
  \paragraph{Quivers as endofunctors}
We start with the following well-known basic connection  between quivers and linear categories \cite{KapVoe}:
With $I$ the finite set of vertices of a quiver of rank $n$, a quiver is specified by a $n \times n$-matrix $Q_{ij}$ with non-negative integer  entries.
On the other hand, a linear functor $F: \Cat{M} \rightarrow \Cat{M}$ for a finite semisimple category $\Cat{M}$ 
is up to isomorphism fixed by the same data: The integer matrix $(\dim_{\Bbbk}(x_{j},F(x_{i})))_{ji}$ with $x_{i}$, $i \in I$ representatives for the simple objects.
Thus there is an obvious  bijection 
\begin{equation}
  \label{eq:main-correspondence}
  \begin{split}
    \{ \text{Quivers with vertices $I$} \}_{/ \sim}\quad &\stackrel{1:1}{\longleftrightarrow} \quad  \{ \text{Linear Endofunctors $F: \Cat{M}\rightarrow \Cat{M}$} \}_{/ \sim}\\
          Q_{ij} \quad &  \longleftrightarrow  \quad (\dim_{\Bbbk}(x_{j},F(x_{i})))_{ji}
  \end{split}
\end{equation}
In this sense we regard an endofunctor $F$ as a  coordinate free description of a quiver.
This point of view is implicit in few articles on the subject, for instance \cite{EOGraphs} and \cite{MOVQuiv}. 
On the one hand this motivates to treat an  endofunctor  $F: \Cat{M} \rightarrow \Cat{M}$ of a general linear category $\Cat{M}$ as a generalized quiver.  
On the other hand the generalization raises the questions to phrase modules over quivers and morphisms of quivers in these categorical terms.
Here we offer a bicategorical approach and systematically translate the two notions.
Indeed,  considering quivers as particular endofunctors, suggests a natural notion of morphism between quivers $Q_{1},Q_{2}$: A functor $F$ between the underlying categories together with a natural transformation $Q_{2}F \rightarrow F Q_{1}$. This leads naturally to the bicategory of quivers.
Also the quiver modules have a natural description in terms of an endofunctor $Q: \Cat{M} \rightarrow \Cat{M}$:
By considering as a module over $Q$ a pair $(m,f)$ of an object $m \in \Cat{M}$ with a morphism $f:Q(m)\rightarrow m$ we recover and generalize the category of modules over a quiver.
The class of locally nilpotent modules over $Q$, which feature for instance in the construction of the  semi-canonical basis by Lusztig, naturally fits in this framework. 

\paragraph{Actions of monoidal categories on quivers}

Given a monoidal category $\Cat{A}$, the relevant notion of a module over $\Cat{A}$ is that of a module category
(e.g. \cite[Sec. 7]{EGNObook}): It is a category $\Cat{M}$ on which the monoidal category $\Cat{A}$ acts functorially.
As a benefit  of the definition of a bicategory of quivers, we have a notion of quivers on $\Cat{M}$ that are compatible with the action of $\Cat{A}$ on $\Cat{M}$.
These so called module quivers correspond to module endofunctors $Q$ on $\Cat{M}$ and the action of $\Cat{A}$ induces an action on the semisimple moduli spaces of $Q$.

As a particular case one can consider the action of $\Cat{A}$ on itself and obtains as bimodule quivers the bimodule endofunctors of $\Cat{A}$, which correspond to the Drinfeld center $\cent(\Cat{A})$ of $\Cat{A}$.
This leads  an important class of quivers for this article: Given a multi-fusion category $\Cat{A}$ and an object $z \in \cent(\Cat{A})$, the \emph{Drinfeld quiver} of $z$ is the endofunctor
\begin{equation}
  Q_{z}=z \otimes -: \Cat{A} \rightarrow \Cat{A}, 
\end{equation}
given by the left tensor product of $z$ on $\Cat{A}$.
As consequence, $Q_{z}$ is naturally  an  $\Cat{A}$-bimodule quiver.
A module $(a,f)$ over $Q_{z}$ has the graphical expression
$
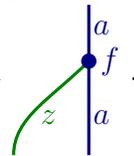
\begin{tikzpicture}[very thick,scale=1,color=blue!50!black, baseline]
\draw (0,-1) -- (0,1); 
\draw[color=green!50!black] (-1,-1) .. controls +(0,0.5) and +(-0.5,-0.5) .. (0,0.25);
\fill[color=blue!50!black] (0,0.25) circle (2.9pt) node[right] (meet2) {{$f$}};
\draw[color=blue!50!black] (-0.1,-0.5) node[right] (A1) {{$a$}};
  \draw[color=blue!50!black] (-0.1,0.7) node[right] (A1) {{$a$}};
    \draw[color=green!50!black] (-0.8,-0.5) node[right] (A1) {{$z$}};
  \end{tikzpicture}.
  $

    The quivers that arise in this manner appear in special cases under the name of McKay graphs or McKay correspondence in the literature, see e.g. \cite{LusQuiv} (see also Example \ref{example:first-expl} \label{item:ext-Dyn}). In the  literature on subfactors they appear as principal graphs, e.g. \cite{GHJ}.

\paragraph{Fusion quiver}

As fusion quiver we define a quiver $Q$ with a rigid monoidal structure on the category of modules $\mod(Q)$, such that those modules for which all arrows act  trivially form a rigid monoidal subcategory $\Cat{A}_{Q}$.
It follows that $\Cat{A}_{Q}$ is a multi-fusion category.

As main result we show that a Drinfeld quiver is a  fusion quiver:
\begin{theorem-n}(Theorem \ref{thm:ten-quiv})
  Let $\Cat{A}$ be a multi-fusion category and $z \in \cent(\Cat{A})$. The Drinfeld quiver $Q_{z}$ has a rigid monoidal structure $\qoti$. The fusion $\qoti$ of  modules $(a,f)$ and $(b,g)$ of $\mod(Q_{z})$ is  in graphical terms the module
  $(a \otimes b, f \qoti g)$ with 
  \begin{equation}
  \label{eq:def-qoti-intro}
 f \qoti  g =
\begin{tikzpicture}[very thick,scale=1,color=blue!50!black, baseline]
    \draw[color=blue!50!black] (-2.0,-1.5) node[above] (z) {{$z$}};
  \draw[color=blue!50!black]   (-1.2,-1.5) node[above] (a) {{$a$}};
  \draw[color=blue!50!black]   (-0.6,-1.5) node[above] (b) {{$b$}};
\draw[color=blue!50!black] (-1.2,1) -- (-1.2,-1); 
\draw[color=blue!50!black] (-0.6,1) -- (-0.6, -1); 
\draw[color=green!50!black] (-2,-1) .. controls +(0,0.5) and +(-0.5,-0.5) .. (-1.2,0.25);  
\fill[color=blue!50!black] (-1.2,0.25) circle (2.9pt) node[right, xshift=-1.8] (meet2) {{$f$}};
\end{tikzpicture} \;
+ 
\begin{tikzpicture}[very thick,scale=1,color=blue!50!black, baseline]
    \draw[color=blue!50!black] (-1.9,-1.5) node[above] (z) {{$z$}};
  \draw[color=blue!50!black] (-1.2,-1.5) node[above] (a) {{$a$}};
  \draw[color=blue!50!black] (-0.6,-1.5) node[above] (b) {{$b$}};
\draw[color=blue!50!black] (-0.6,1) -- (-0.6,-1); 
\draw[color=blue!50!black] (-1.2,1) -- (-1.2,-1); 
\draw[color=white, line width=4pt] (-1.9,-1) .. controls +(0,0.5) and +(-0.5,-0.5) .. (-0.6,0.25);  
\draw[color=green!50!black] (-1.9,-1) .. controls +(0,0.5) and +(-0.5,-0.5) .. (-0.6,0.25);  
\fill[color=blue!50!black] (-0.6,0.25) circle (2.9pt) node[right, xshift=-1.8] (meet2) {{$g$}};
\end{tikzpicture}\, .
\end{equation}
The  locally nilpotent modules of $Q_{z}$ form a rigid monoidal subcategory. 
\end{theorem-n}

After the publication of this work, the author became aware that this monoidal product appears already as \cite[Eq. 2.2]{MR10}  in the literature on conformal field theory. Therefor full credit on the theorem belongs to these authors. 
We note that in \cite{ARW25} this monoidal category appears as a special case of a more general monoidal category of perturbed topological defects.\footnote{The author is grateful to I.~Runkel for explaining these results and further promising connections between the two approaches.}

In the sequel we obtain in the  graphical calculus a simple expression for these fusion product on quiver modules, which directly allows for generalizations---semisimplicity
of the category $\Cat{A}$ was not needed for this theorem.
However, even in the semisimple case one obtains many examples, including the Hopf quivers of \cite{CRHopf} as well as the doubles of extended Dynkin quivers,
Example \ref{example:first-expl}.

In total a fusion quiver provides with the forgetful functor $\mod(Q) \rightarrow \Cat{A}_{Q}$ a retract diagram of monoidal functors 
 \begin{equation}
    \label{eq:retract-intro}
    \begin{tikzcd}
      \Cat{A}_{Q} \ar{r}{} \ar[bend left]{rr}{\id} & \mod(Q) \ar{r}{} & \Cat{A}_{Q}.
    \end{tikzcd}
  \end{equation}
  Note that the category $\mod(Q)$ is non-semisimple in case $Q$ is a non-trivial quiver.

     Fusion quivers and Drinfeld fusion quivers naturally organise in bicategories  $\fusqui$ and $\Dfusqui$, respectively.  Moreover,  by considering the fusion structure on the modules, we obtain a natural 2-functor
\begin{equation}
  \Dfusqui\rightarrow \fusqui.
\end{equation}

\paragraph{Relations on fusion quivers}

Many algebras of interest arise as quotients of path algebras of quivers by relations.
In case of a Drinfeld quiver $Q_{z}$, we present a class of relations that respect the fusion product of modules.

The categorical point of view on quivers allows to consider graphically the following relations on $Q_{z}$:
Given a morphism $\varphi: y \rightarrow z^{\otimes n}$ in $\cent(\Cat{A})$, we say that a $Q_{z}$-module $ f: Q_{z}(a) \rightarrow a $ satisfies the relation $\varphi$, if:
 \begin{equation}
    \label{eq:varphi-rel-rel}
\begin{tikzpicture}[very thick,scale=1,color=blue!50!black, baseline]
\draw (0,-2) -- (0,1.5); 
\draw[color=green!50!black] (-0.7,-0.75) .. controls +(0,0.5) and +(-0.5,-0.5) .. (0,0.15); 
\draw[color=green!50!black] (-1.3,-0.75) .. controls +(0,0.5) and +(-0.5,-0.5) .. (0,0.8);
\draw[color=green!50!black] (-1.6,-0.75) .. controls +(0,0.5) and +(-0.5,-0.5) .. (0,1.1);
\draw (-1.15,-1) -- (-1.15,-2); 
\fill[color=blue!50!black] (0,0.15) circle (2.9pt) node[right] (meet2) {{$f$}};
\fill[color=blue!50!black] (0,0.8) circle (2.9pt) node[right] (meet2) {{$f$}};
\fill[color=blue!50!black] (0,1.1) circle (2.9pt) node[right] (meet2) {{$f$}};
  \draw (-1.15,-1) node[minimum height=0.5cm,minimum width=1cm,draw,fill=white] {{$\varphi$}};
%
\draw (0,-2.5) node[above] (X) {{$m$}};
\draw[color=blue!50!black] (-1.15,-2.5) node[above] (A1) {{$y$}};
\draw[color=green!50!black] (-0.17,-0.1) node[above] (A1) {{$\vdots$}};
\draw[color=green!50!black] (-0.22,-0.6) node[above] (A1) {{$z$}};
\draw[color=green!50!black] (-0.6,-0.21) node[above] (A1) {{$z$}};
\draw[color=green!50!black] (-1,0.18) node[above] (A1) {{$z$}};
\end{tikzpicture} 
\, =0.
  \end{equation}
  This can be phrased in terms of the usual notion of relations on quivers.

  We say that $\varphi$ is a $q$-relation, if there is a primitive $n$-th root of unity $q$ such that whenever
  one postcomposes $\varphi$ with a braiding of two of the tensor factors in $z^{\otimes n}$, the resulting morphism is $q \cdot \varphi$. 
    \begin{theorem-n}(Theorem \ref{theorem:main-rel})
    Let $Q_{z}$ be a Drinfeld quiver and $\varphi$ a $q$-relation. Then the category of modules of $Q_{z}$ that satisfy the relation $\varphi$ is a rigid monoidal subcategory of $\mod(Q_{z})$. 
  \end{theorem-n}
  
  This construction is  applied to reproduce the rigid monoidal category of modules over the Taft Hopf algebra.
  By considering length two $q$-relations of the type $\varphi: \unit \rightarrow z \otimes z$, one can obtain certain preprojective algebras as corresponding quotients of the path algebra. 

For a  pivotal braided multi-fusion category $\Cat{A}$ we consider two situations in detail. First, if
  $z \in \Cat{A}$ is self-dual with Frobenius-Schur-indicator $-1$ and twist $1$, then the quiver $Q_{z}$ has
  a length two $q$-relation whose corresponding rigid monoidal category of modules corresponds to the modules over a  preprojective algebra in the classical sense. 
  
  In the second case we apply this construction, to equip the modules over a twisted version of the preprojective algebra over certain quivers with a rigid  monoidal structure,
    see Corollary \ref{corollary:q-preproj}.

 The construction of rigid monoidal structures for preprojective algebra  is new to our knowledge. This can be regarded as a first step in the application of our methods to construct new families of rigid monoidal categories. Note that for all non-trivial quivers the categories  of modules are non-semisimple, thus the Drinfeld quivers $Q_{z}$ provide a family of non-semisiple, non-finite hereditary rigid monoidal categories, see Remark \ref{remark:TypesDrinfeld}.

We conclude with an outline of the organization of the paper. 
To make it accessible for readers of both communities, we give a rather long introduction to quivers in Section \ref{sec:quiv-semis-categ}, as well as including definitions and results from tensor categories in Appendix \ref{sec:categorical-notions}.  
In Section \ref{sec:quiv-semis-categ} we  furthermore relate in   Proposition \ref{proposition:modules-endof}  quivers to endofunctors and
present in Definition \ref{definition:quiver-of-fun-mor-bicat} the bicategory of quivers.
In Section \ref{sec:fusi-categ-their} we consider actions of monoidal categories on quivers  in  
Definition \ref{definition:Action-quiv}.  In Proposition \ref{proposition:gen-quiv-A}, those are related to module categories and in 
in Remark \ref{remark:action-moduli} we consider the induced actions on the moduli spaces of semisimple quiver modules. 
In Section \ref{sec:tensor-quivers} we present  the main Definition \ref{definition:fus-quiver} of a fusion quiver. In Subsection \ref{sec:fusi-quiv-struct} we show in Theorem \ref{thm:ten-quiv} that Drinfeld quivers are fusion quivers and in 
Example \ref{example:first-expl} we consider fusion structures on  extended Dynkin quivers. 
We list the known classification of Drinfeld quivers up to rank two in Subsection \ref{sec:fusion-quivers-low}. 
In Theorem \ref{theorem:Moduli-structures} we consider the induced structure on the  moduli spaces. 

In Section \ref{sec:relat-fusi-quiv} we study homogenous relations for Drinfeld quivers that are compatible with the rigid monoidal structure with main Theorem \ref{theorem:main-rel}.

In  Section \ref{sec:prepr-algebra} we apply these constructions to  modules over preprojective algebra
in two cases, 
see Theorem \ref{theorem:prepro-selfdual} and 
Theorem \ref{theorem:double-cent-rigid-mon}. In  Definition \ref{definition:double-quiv-prepro}
we propose a generalization of categories of  preprojective modules to general endofunctors. 

Throughout this article we  work over a fixed algebraically closed field $\Bbbk$.

\vspace{1em}
{\small \noindent{\sc Acknowledgments: }
The author thanks the anonymous referee for helpful suggestions and Catherine Meusburger and  the Emmi-Noether-Seminar Erlangen for useful discussions and comments.
}

\section{Quivers and semisimple  categories}
\label{sec:quiv-semis-categ}
We build on the  basic correspondence between endofunctors of finite semisimple categories and finite quivers.
This correspondence allows us to define a natural bicategory of quivers and enables to describe modules over a  quiver in simple diagrammatic terms.

\subsection{Finite quivers}
\label{sec:finite-quivers}

We refer to \cite{reineke08} for a detailed introduction to quivers and their modules. 
A \emph{quiver} $Q=(I, Q^1; s,t)$ consists of  sets $I$  of vertices and  $Q^1$ of arrows  together with  source and target maps
$s,t: Q^1 \rightarrow I$. We write $\alpha: i \rightarrow j$ for $\alpha \in Q^{1}$ with source $i$ and target $j$.  The set of arrows from vertex $i$ to vertex $j$ is denoted by  $Q_{ij}$.
The quiver is called  non-trivial, if $Q^{1}$ is non-empty. 
We assume  that all  quivers are finite, i.e. both the vertices and arrows form finite sets. The \emph{rank} of a quiver is the cardinality of the set of vertices. 

Thus, a quiver is just a finite directed graph, and it is typically drawn as a graph in the plane to visualize the incidence relations, for example:
\begin{equation}
  \label{eq:pic0-quiv}
  Q=
  \begin{tikzpicture}[thick,scale=1,color=blue!50!black, baseline]
    \fill[color=blue!50!black] (-1.0,0)   circle (2.9pt) node[right, xshift=-1.8] (a) {};
    \fill[color=blue!50!black]   (1,0)   circle (2.9pt) node[right, xshift=-1.8] (b) {};
    \draw[->,color=blue!50!black] (-1,0)+(0.1,0.1) .. controls +(0.5,0.5) and +(-0.5,0.5) .. (1-0.1,0+0.1);  
    \draw[<-,color=blue!50!black] (-1+0.1,0-0.1) .. controls +(0.5,-0.5) and +(-0.5,-0.5) .. (1-0.1,0-0.1);  
    \draw[->,color=blue!50!black] (-1.1,0.1) .. controls +(-1,1) and +(-1,-1) .. (-1.1,-0.1);  
\end{tikzpicture} \; .
\end{equation}

\paragraph{Quiver modules}
Quivers are mainly of interest for their categories of representations, here called modules:
For a quiver $Q$, a \emph{module $(V,f$) over $Q$} or a \emph{Q-module} is a collection of   vector spaces $(V_i)_{i \in I}$ together with a collection of linear maps
$(f_{\alpha}: V_i \rightarrow V_j)_{\alpha: i \rightarrow j}$, indexed by the arrows $\alpha$ of $Q$. The module $V$ is called finite dimensional, if all $V_{i}$ are finite dimensional.
For a collection of vector spaces $(V_{i})_{i \in I}$ there is a distinguished module, where all linear maps $f_{\alpha}=0$ are the zero morphisms. Such a module is called \emph{vertex semisimple} in the sequel.

The  \emph{dimension vector} of a finite dimensional  module $V$ is the tuple  $x=(x_{i})_{i \in I}\in \N^{n}$ with $x_i=\dim_{\Bbbk}(V_i)$ and $n=|I|$.
A \emph{morphism of $Q$-modules} $h: (V,f) \rightarrow (W,g)$
is a collection of linear maps $h_i: V_i \rightarrow W_i$ for all $i \in I$, such that
\begin{equation}
  \label{eq:quiver-morph}
  \begin{tikzcd}
    V_i \ar{d}{h_i} \ar{r}{f_\alpha}& V_{j}  \ar{d}{h_j} \\
    W_{i} \ar{r}{g_{\alpha}}  & W_j 
  \end{tikzcd}
\end{equation}
commutes for all $\alpha \in Q_{ij}$ and  all vertices $i,j$. 
The corresponding linear category of  $Q$-modules is denoted $\Mod(Q)$ with the full subcategory
$\mod(Q)$ of finite dimensional modules. The full subcategory of $\mod(Q)$ of modules with fixed dimension vector $x$ is denoted
$\mod_x(Q)$. The vertex semisimple modules $\mod(Q)^{\vses}$   form a subcategory of $\mod(Q)$, however in general those are not the only semisimple modules over $Q$.

Let us reformulate these structures in a way that is suited for generalizations. Note that a  module over $Q$ is in particular an  $I$-graded vector space: Let $\Vect=\Vect_{\Bbbk}$ be the category of all (not necessarily finite-dimensional) $\Bbbk$-vector spaces.  The linear category  $\Vect^{I}$ of $I$-graded vector spaces is just the $|I|$-fold product of $\Vect$, its objects are tuples $(V_{i})_{i\in \I}$ with $V_{i} \in \Vect$.
By definition of $\mod(Q)$ for $Q$ a quiver with vertex set $I$ there is a forgetful functor
\begin{equation}
  \label{eq:forget-Q}
  U: \Mod(Q) \longrightarrow \Vect^{I},
\end{equation}
and if $\vect$ denotes the category of finite dimensional vector spaces, $U$ restricts to a functor $U: \mod(Q) \rightarrow \vect$. 
The $I$-graded vector space $V=(V_{i})_{i \in I}$ underlying a module over $Q$ is called  the \emph{dimension object} $V \in \Vect^{I}$ of the module.
It follows moreover that there is a canonical equivalence of categories $\Mod(Q)^{\vses} \cong \Vect^{I}$ and that from the inclusion of $\Mod(Q)^{\vses} \rightarrow \Mod(Q)$ we obtain a commuting diagram of functors
 \begin{equation}
    \label{eq:retract}
    \begin{tikzcd}
      \Vect^{I} \ar{r}{} \ar[bend left]{rr}{\id} & \Mod(Q) \ar{r}{U} & \Vect^{I}.
    \end{tikzcd}
  \end{equation}
Regarding the set $Q_{ij}$ as the basis of the vector space $\Bbbk Q_{ij}$, a module over $Q$ is the same as collection of linear maps $\Bbbk Q_{ij} \rightarrow  \Hom_{\Vect}(V_{i},V_{j})$.
Thus a $Q$-module is an element in the   \emph{space of modules over $Q$ with dimension object   $V \in \Vect^{I}$}:
\begin{equation}
  \label{eq:lin-rep-space}
  \Mod_{V}(Q)=\oplus_{i,j \in I}        (\Bbbk Q_{ij})^{*} \otimes \Hom_{\Vect}(V_{i},V_{j}).
\end{equation}

For a dimension object $V$ the group $\Gl(V)= \prod_{i \in I} \Gl(V_{i})$ acts on $\mod_{V}(Q)$ by conjugation:
Given a module $(f_{\alpha}: V_{i} \rightarrow V_{j})_{\alpha:i\rightarrow j}$, the module
$g \act f$ for $(g_{i})_{i \in I} \in \Gl(V)$ has on the arrow $\alpha: i \rightarrow j$ the linear map
\begin{equation}
  \label{eq:action-gauge}
 (g \act f)_{\alpha}= g_{j} \circ f_{\alpha} \circ g_{i}^{-1}. 
\end{equation}

Finally we note that given quivers $Q_{1}$ and $Q_{2}$ on the same vertex set $I$, their \emph{sum} $Q_{1} \oplus Q_{2}$ has vertices $I$ and as arrows the union
$(Q_{1} \oplus Q_{2})_{ij}=(Q_{1})_{ij} \cup (Q_{2})_{ij}$ of the arrows of $Q_{1}$ and $Q_{2}$.
\paragraph{Linear quivers}
From the point of view of the quiver modules, it is only a minor step to forget the distinct basis of $\Bbbk Q_{ij}$ and to consider the setting of 
 a \emph{linear quivers}  $Q$ with  vertices $I$: $Q$ consists of  a collection of finite dimensional vector spaces $Q(i,j)$ for all $i,j \in I$,
the \emph{space of arrows from $i$ to $j$}, whose non-zero elements will be called arrows. For an $I$-graded vector space $V=(V_{i})_{i \in I}$,
the space of modules over $Q$ on $V$ is then
\begin{equation}
  \label{eq:lin-rep-space-lin-quiv}
  \Mod_{V}(Q)=\oplus_{i,j \in I}        Q(i,j)^{*} \otimes \Hom_{\Vect}(V_{i},V_{j}).
\end{equation}
Thus, a finite quiver $Q$ induces a linear quiver with arrow space $Q(i,j)= \Bbbk Q_{ij}$. Conversely, upon picking bases we obtain a finite quiver from a linear quiver.
A collection of bases for all $Q(i,j)$ we call a \emph{basis of a linear quiver}. Under this correspondence, the two categories  of modules are canonically equivalent
and the sum of quivers corresponds to the taking the direct sum of the arrow spaces in the setting of linear quivers. In the sequel, by quiver we mean linear quiver and
a linear quiver with basis will be called \emph{combinatorial quiver}. 

The setting of linear quivers allows for the following  morphisms between quivers with the same vertices. 
\begin{definition}
  \label{definition:id-morph-lin-Quiv}
  Let $Q=(I,Q^{1};s,t)$ and $P=(I,P^{1};s',t')$ be two linear quivers  with the same vertices. 
  A \emph{morphism $\varphi:Q\rightarrow P$ over the identity} from $Q$ to $P$ consists of a collection of linear maps $\varphi_{i,j}:  P(i,j) \rightarrow Q(i,j) $. 
\end{definition}
We choose  the given direction of the linear maps $\varphi_{i,j}$, because by Equation \eqref{eq:lin-rep-space-lin-quiv}, a morphism $\varphi: Q \rightarrow P$ induces a linear map
\begin{equation}
  \label{eq:lin-map-mod}
  \begin{tikzcd}[row sep=large, column sep=large]
    \Mod_{V}(Q)) \ar{r}{\oplus_{i,j} \varphi_{i,j}^{*} \otimes \id} &   \Mod_{V}(P).
  \end{tikzcd}
    \end{equation}
It will become clear after Definition \ref{definition:bicat-quiv}, why we call these morphisms over the identity.

Clearly, quivers of rank $I$ and  morphisms over the identity  form a category $\Quiv_{I}$.

The paths of a quiver carry two distinct known algebraic structures, which are compared for instance in \cite{Simson}.
\begin{definition}
  \label{definition:path-structures}
  Let $Q$ be a (linear) quiver with vertices $I$ and $i,j \in I$. 
  \begin{definitionlist}
  \item The \emph{space  of paths of length $n$} in $Q$  from $i$ to $j$ is the vector space
    \begin{equation}
      \label{eq:paths-space}
      Q^{n}(i,j)= \oplus_{i_{1}, \ldots, i_{n-1} \in I}     Q(i_{n-1},j)\otimes \ldots \otimes     Q(i_{1},i_{2})  \otimes  Q(i,i_{1}), 
    \end{equation}
    whose elements are called \emph{paths of length $n$}. 
An elementary tensor  in this space is called a \emph{composed path} of length $n$ and denoted 
$\alpha_{n} \circ \ldots \circ \alpha_{1}$ with $t(\alpha_{i})=s(\alpha_{i+1})$ for all $i=1, \ldots, n-1$.
In case we are given a basis of each $Q(i,j)$ (i.e. the  quiver is a combinatorial quiver), the composed paths which consist of  a composition of elements in the bases, are called \emph{basic paths}. 
The paths of length $0$   form the vector space $\Bbbk I$ with canonical basis $e_{i}$ with $i \in I$.
    
  We set $Q^{n}= \oplus_{i,j \in I}Q^{n}(i,j)$ and denote $\Bbbk Q= \oplus_{n \in \N} Q^{n}$.
  \item The \emph{path algebra} $\PQ$ of $Q$ is the vector space $\PQ= \Bbbk Q$ with multiplication defined by  the composition  of paths and unit $u= \sum_{i \in I} e_{i}$. 
  \item The \emph{path coalgebra} $\CQ$ of $Q$
    is as vector space $\CQ=\Bbbk Q$ and
   the  coproduct is defined on  a composed path $\omega \in \CQ$  by
\begin{equation}
  \label{eq:path-coprod}
  \Delta(\omega)= \sum_{\alpha,\beta \in \Bbbk Q, \omega=\beta \circ \alpha}\beta \otimes \alpha \in \CQ \otimes \CQ
\end{equation}
and extended linearly to $\CQ$. The counit $\epsilon:  \CQ  \rightarrow \Bbbk$ of $\CQ$ is the linear map defined by $\epsilon(e_{i})=1$ for all $e_{i} \in \CQ^{0}=\Bbbk I$ extended by zero to $\CQ$. 
  \end{definitionlist}
\end{definition}
It is straightforward to see that $\PQ$ is an algebra and $\CQ$ a coalgebra. 
Note that $\Bbbk Q$ is infinite dimensional if and only if there exists an \emph{oriented cycle} in $Q$, i.e.
a loop (a path from a vertex to itself) of positive length. It is a standard result, that the category of finite dimensional left $\PQ$ modules, ${}_{\PQ}\mod$, is equivalent to the
category $\mod(Q)$ of quiver modules. 
For a basic path $\omega:i  \rightarrow j$ and a module $(V,f)$ over $Q$, the linear map $f(\omega)$ is thereby given as the composite of the linear maps associated to the paths  of length one in $\omega$. 

To describe the comodules over $\CQ$ we consider  for a linear quiver $Q$ the \emph{dual quiver} $Q^{*}$, which is the linear quiver with vertices $I$ and arrows from $i$ to $j$ defined

with the dual vector space as $Q^{*}(i,j)=Q(i,j)^{*}$.

 The \emph{canonical pairing} of $\CQ^{*}$ and $\PQ$  is the bilinear map
\begin{equation}
  \label{eq:quiver-pairing}
  <-,->: \CQ^{*} \times \PQ \rightarrow \Bbbk, 
 \end{equation}
 which for composed paths $\omega_{n} \circ \ldots \circ \omega_{1} \in  \CQ^{*}$ and $\lambda_{n} \circ \ldots
\circ  \lambda_{1} \in \PQ$ is defined as
 \begin{equation}
   \label{eq:pairing-composed-paths}
   < \omega_{n} \circ \ldots \circ \omega_{1},\lambda_{n} \circ \ldots \circ  \lambda_{1} >= <\omega_{n},\lambda_{n}>\cdot \ldots \cdot <\omega_{1}, \lambda_{1}>,
 \end{equation}
using on the right hand side the canonical $<-,->: \PQ(i,j)^{*} \times \PQ(i,j)\rightarrow \Bbbk $. Finally, Equation \eqref{eq:pairing-composed-paths} is extended linearly. In case that the lengths of paths or the sources  and targets of arrows does not match, the pairing is zero by definition. 

\begin{definition}
Let $Q$ be a linear quiver. A module $(V,f) \in \mod(Q)$ is called \emph{locally nilpotent}, if for all $v \in V$, $f_{\omega}(v)=0$ for almost all paths $\omega \in \PQ$.
The full subcategory of locally nilpotent modules over $Q$ is denoted $\modln(Q)$.
  \end{definition}

And we obtain,  which is shown for combinatorial quivers in   \cite[Prop.6.1]{CKQSplit}, see also \cite[Thm. 3.14]{Simson}. 

\begin{proposition}
  \label{proposition:path-coalg-comod}
  Let $Q$ be a linear quiver with dual quiver $Q^{*}$.
   The locally nilpotent (not necessarily finite dimensional) $Q$-modules are equivalent to the right $\CQ^{*}$-modules, that is
    \begin{equation}
      \label{eq:43}
      {}_{\PQ}\Modln \cong \Comod_{\CQ^{*}}. 
    \end{equation}
Restricting to finite dimensional modules provides an equivalence $    {}_{\PQ}\modln \cong \comod_{\CQ^{*}}. $
\end{proposition}
\begin{proof}
  
  We adopt the proof of  \cite[Prop.6.1]{CKQSplit} to the linear setting:
   Fix a basis $\{\lambda_{i,j;\alpha}\}_{\alpha}$ of $Q(i,j)$ for all $i,j \in I$. The corresponding basic paths are a basis for the spaces of paths $Q^{n}(i,j)$   and we obtain dual bases $\{\lambda_{i,j}^{\alpha}\}_{\alpha}$  for $Q^{*}(i,j)$, such that each path $\omega \in \PQ$ defines a dual path $\omega^{\vee} \in \CQ^{*}$ with respect to this basis. 
    
    For $(V,\Delta_{V})$ a locally nilpotent module, define the coaction
    \begin{equation}
      \label{eq:44}
      \Delta_{V}: V \rightarrow V \otik \CQ^*, \quad v \mapsto \sum_{\text{basic paths } \omega}f(\omega)(v) \otimes \omega^{\vee},
    \end{equation}
   By the locally nilpotence this is well-defined and  it follows, that $\Delta_{V}$ is independent of the chosen bases and indeed a coaction of $\CQ^{*}$.

    Conversely, given $(V,\Delta_{V})$ in $\Comod_{\CQ^*}$,  we define  the linear map $f(\lambda): V \rightarrow V$  for an arrow $\lambda \in Q(i,j)$  as
     \begin{equation}
      \label{eq:45}
      f(\lambda)(v)= (\id \otimes <-,\lambda>) \circ \Delta_{V}(v), 
    \end{equation}
    using the covector $<-,\lambda>: \CQ^{*} \rightarrow \Bbbk$ for $v \in V$.
    As in  \cite[Prop.6.1]{CKQSplit} are the two constructions  functorial and mutually inverse and provide the claimed equivalences of categories. It is clear, that the equivalence is compatible with the restriction to finite dimensional modules. 
\end{proof}
\paragraph{Relations for quivers}
Quotients of quivers by relations give rise to an important class of algebras: A \emph{relation} $\omega$ for $Q$ is a linear combination of paths, i.e. $\omega \in \Bbbk Q$. The corresponding path algebra with relation $\omega$ is $\PQ/\Omega_{\omega}$, where $\Omega_{\omega} \subset \PQ$ is the both-sided ideal generated by $\omega$. The modules
in $\mod(\PQ/\Omega_{\omega})$ are canonically identified with those quiver modules $(V,f) \in \mod(Q)$, such that
$f(\omega): V \rightarrow V$ is the zero map.

It is well-known that up to Morita equivalence, all finite dimensional algebras can be described by quivers with relations.

\paragraph{Moduli spaces of $Q$-modules}

For a quiver $Q$ there exist geometric tools to study its modules.
We only very briefly review this important subject

    following \cite{reineke08}. In Section \ref{sec:fusion-moduli-spaces} we will apply our algebraic results on categories of quivers to these geometric objects.

    The representation space $\mod_{m}(Q)$ of $Q$ with fixed dimension vector $m \in \N^{I}$  is in particular an affine space with  the action    \ref{eq:action-gauge}
of the reductive algebraic group $\Gl_{m}$.

Given an algebraic variety $X$ with action of an algebraic group $G$ on $X$, the \emph{categorical quotient} of
$X$ is (if it exists), a variety 
$M$ together with a $G$-invariant morphism $\pi: X \rightarrow M$, such that for all
$G$-invariant morphisms $f: X \rightarrow Y$ there exists a unique  morphism $\widehat{f}: M \rightarrow Y$, such that $f= \widehat{f} \circ \pi$, i.e.
\begin{equation}
  \begin{tikzcd}
  X \ar{r}{f} \ar{d}{\pi}  & Y \,. \\
  M \ar[dashed]{ur}[right, yshift=-3pt]{\widehat{f}} &   
  \end{tikzcd}
\end{equation}
In case of $X=\mod_{m}(Q)$ and $G=\Gl_{m}$, the categorical quotient exists, see e.g. \cite[Sec. 3.3]{reineke08} and is given by  the moduli space   $\Mses_{m}(Q)$ of semisimple $Q$-modules with dimension vector $m$.
That is,  $\Mses_{m}(Q)$ is a variety whose points correspond to the $G$-orbits of semisimple $Q$-modules with dimension vector $m$. We collect this as 
\begin{definition}
  \label{definition:ses-moduli}
  Let $Q$ be a quiver. For  $m \in \N^{I}$ the categorical quotient of the action of $G$ on $\mod_{m}(Q)$ is called the \emph{moduli space of $Q$-modules with dimension vector $m$}, denoted  $\Mses_{m}(Q)$. The product over all dimension vectors is the \emph{moduli space of semisimple $Q$-modules} $\Mses(Q)=\sqcup_{m \in \N^{I}}\Mses_{m}(Q)$.
\end{definition}

Basic properties of the  varieties $\Mses_{m}(Q)$ are studied in  \cite{BCSes}: The coordinate rings are generated by traces along oriented cycles in $Q$ and the varieties admit explicit finite stratifications into locally closed smooth irreducible subvarieties.

However, there are  more refined versions of moduli spaces for quivers, see \cite{reineke08} for an introduction.

\subsection{Semisimple categories and quivers}
In Appendix \ref{sec:mono-categ-their} we summarize standard concepts for linear categories
       following \cite{EGNObook}.
All categories and functors will be $\Bbbk$-linear. For finite semisimple categories we choose a  set $I$ labeling representatives for the isomorphism classes of the simple objects.
The cardinality of $I$ is called the rank of the category.

We first relate (linear) quivers to (linear) endofunctors to obtain a useful characterization of quivers. 
By  Lemma \ref{lemm:ses-cat-fun}, a  linear functor on a semisimple category is   determined by the values on $I$.
Recall \eqref{eq:repres-vect}  the product $V \otik m \in \Cat{M}$ of a vector space $V \in \vect$ with an object $m \in \Cat{M}$.

\begin{definition}
  \label{definition:quiv-of-fun}
  Let $\Cat{M}$ be a finite semisimple category with labeling set $I$ and corresponding simples  $(x_{i})_{i \in I}$  of $\Cat{M}$.
  \begin{definitionlist}
  \item Let  $R: \Cat{M} \rightarrow \Cat{M}$ a linear functor. The \emph{quiver $Q(R)$ of $R$}  has as
  vertices $I$ and as arrow spaces  from $i$ to $j$ for all $i,j \in I$
  \begin{equation}
    \label{eq:arrow-space}
    Q(R)(i,j)=\Hom_{\Cat{M}}(R(x_{i}),x_{j})^{*}
  \end{equation}
\item   Let $Q$  be  a linear quiver  on the vertex set $I$. The \emph{endofunctor $R(Q): \Cat{M} \rightarrow \Cat{M}$ of $Q$} is defined by
  \begin{equation}
    \label{eq:Fun-from-Q}
   R(Q)(x_{i})= \oplus_{j \in I}Q(i,j) \otik x_{j}. 
  \end{equation}
  \end{definitionlist}
\end{definition}
With the natural isomorphism of Equation \eqref{eq:nat-iso-dual} we obtain the isomorphism
 \begin{equation}
    \label{eq:arrow-space-iso}
    Q(R)(i,j)=\Hom_{\Cat{M}}(R(x_{i}),x_{j})^{*} \cong \Hom_{\Cat{M}}(x_{j},R(x_{i})). 
  \end{equation}
  for all arrow spaces of $Q(R)$. Note that this isomorphism depends on the choice of the pairing \eqref{eq:pairing-ses},
  for that reason we prefer to define the arrow space as in Equation \eqref{eq:arrow-space}.

\begin{proposition}
  \label{proposition:quiv-as-endof} 
  The assignments of Definition \ref{definition:quiv-of-fun} extend for all semisimple linear categories $\Cat{M}$ of rank $I$ to an equivalence of categories  
  \begin{equation}
    \label{eq:10}
     \Quiv_{I}\cong \End(\Cat{M})^{\opp}.
  \end{equation} 
\end{proposition}
\begin{proof}
  First we show that the maps $R$ and $Q$ from  Definition \ref{definition:quiv-of-fun} naturally extend to functors
  \begin{equation}
    \label{eq:Fun-R-Q}
    Q:     \End(\Cat{M}) \longrightarrow  \Quiv_{I} \quad \text{and} \quad R:  \Quiv_{I} \longrightarrow   \End(\Cat{M}).
  \end{equation}
  Indeed, if $\varphi:R \rightarrow R' $ is a natural transformation, the collection of linear maps
  \begin{equation}
    \label{eq:4}
   Q(\varphi)(i,j)=\Hom(\varphi_{x_{i}},x_{j})^{*}: \Hom_{\Cat{M}}(R(x_{i}),x_{j})^{*} \rightarrow \Hom_{\Cat{M}}(R'(x_{i}),x_{j})^{*} 
  \end{equation}
  defines a morphism $Q(\varphi): Q(R') \rightarrow Q(R)$ in $\Quiv_{I}$ and thus the functor $Q$.
  On the other hand, a collection of linear maps $\eta(i,j): Q(i,j) \rightarrow Q'(i,j)$, defines for all $x_{i}$ a morphism
  \begin{equation}
    \label{eq:3}
    \oplus_{j \in I} \eta(i,j) \otimes \id_{x_{j}}: \oplus_{j} Q(i,j) \otimes x_{j} \longrightarrow \oplus_{j} Q'(i,j) \otimes x_{j},
  \end{equation}
  and thus by Lemma \ref{lemm:ses-cat-fun} a natural transformation $R(\eta): R(Q) \rightarrow R(Q')$ and so the functor $R$.

We next show that for all linear functors $R: \Cat{M} \rightarrow \Cat{M}$   on the vertex set $I$, 
  the functor $R(Q(R))$ is canonically equivalent to $R$. 
 For a simple object $x_{i} \in \Cat{M}$,
  \begin{equation}
    \label{eq:58}
    R(Q(R))(x_{i}) = \oplus_{s}Q(R)(i,s) \otimes x_{s} = \oplus_{s} \Hom_{\Cat{M}}(R(x_{i}), x_{s})^{*} \otimes x_{s}
 \cong R(x_{i}),
\end{equation}
using Lemma \ref{lemma:facts-ses}, \refitem{item:coend-ses}.
This provides by Lemma \ref{lemm:ses-cat-fun} a natural isomorphism $R(Q(R)) \cong R$.

Now let $Q$ be a linear quiver on $I$. 
We then  compute
\begin{equation}
  \label{eq:61}
  \begin{split}
  Q(R(Q))(i,j) =& \Hom_{\Cat{M}}(R(Q)(x_{i}),x_{j})^{*}= \Hom_{\Cat{M}}(\oplus_{s} Q(i,s) \otimes x_{s},x_{j})^{*}\\
  \cong & \oplus_{s \in I}Q(i,s)^{**} \otimes \Hom_{\Cat{M}}(x_{s},x_{j})^{*} \cong Q(i,j),    
  \end{split}
\end{equation}
using the canonical isomorphism $V \cong V^{**}$ for a finite-dimensional $V \in \vect$ and the fact that $\Hom_{\Cat{M}}(x_{s},x_{j}) \cong \Bbbk$ if $s=j$ and zero else. Thus we obtain for all $i,j \in I$ a canonical isomorphism 
    $Q(i,j) \cong Q(R(Q))(i,j)$  and thus an isomorphism  $Q \cong Q(R(Q))$ of quivers. 
   It is straightforward to see that the two isomorphisms are natural and thus provide the claimed equivalence of categories.  
\end{proof}

As justified by Proposition \ref{proposition:quiv-as-endof}, from now on a quiver $Q$ is for us the same as a linear endofunctor $Q: \Cat{M} \rightarrow \Cat{M}$ 
of a finite semisimple category $\Cat{M}$, denoted also by $Q$.
As a consequence, the space of paths from $i$ to $j$ of length $n$ in $Q$ is
\begin{equation}
  \label{eq:paths-gen-Q}
  \PQ^{n}(i,j)=\Hom_{\Cat{M}}(Q^{n}(x_{i}),x_{j})^{*} \stackrel{~\eqref{eq:nat-iso-dual}}{\cong} \Hom_{\Cat{M}}(x_{j},Q^{n}(x_{i})). 
\end{equation}
The sum of linear quivers corresponds to the direct sum of the corresponding endofunctors.
Proposition \ref{proposition:quiv-as-endof} allows for the following definition:
\begin{definition}
  \label{definition:gen-quiv}
A \emph{generalized quiver} is a endofunctor $Q: \Cat{M} \rightarrow \Cat{M}$ of a linear category $\Cat{M}$. 
\end{definition}

In sight of the correspondence between quivers and endofunctors, we define modules for endofunctors as follows. 
\begin{definition}
  \label{definition:module-endo}
  Let $R: \Cat{M} \rightarrow \Cat{M}$ be a endofunctor of a category $\Cat{M}$. 
  A \emph{module over $R$} is  a  pair $(m,f)$ of an object  $m \in \Cat{M}$ and a morphism $f: R(m) \rightarrow m$. The object $m$ is  called  the \emph{dimension object} of the module.
  A \emph{morphism of modules} $h: (m,f) \rightarrow (n,g)$ is
  a morphism $h: m \rightarrow n$ such that $g \circ R(h)= h \circ f$. This defines the category $\mod(R)$ of modules over $R$  and the subcategories $\mod_{m}(R)$ of modules with fixed dimension object $m \in \Cat{M}$.

  A module $f: R(m) \rightarrow m$ is called \emph{locally nilpotent}, if there exists a $n \in \N$ with $f^{\ast n}=0$, where
   \begin{equation}
         \label{eq:f-ast}
f^{\ast n}=f \circ R(f) \circ \ldots  \circ  R^{n-1}(f): R^{n}(m) \rightarrow m. 
\end{equation}
The subcategory of locally nilpotent modules is denoted $\modln(R)$. 
 \end{definition}
 In the literature, an object $ (m,f) \in\mod(R)$ is also called an $R$-algebra (see \cite{BW-cat}), but  
 we find this terminology confusing. Again we have a forgetful functor $U: \mod(R) \rightarrow \Cat{M}$
and an inclusion of $\Cat{M}$ in $\mod(R)$ by the zero morphism. In total we obtain a diagram 
 as in  Equation \eqref{eq:retract}  
  \begin{equation}
    \label{eq:retract-gen}
    \begin{tikzcd}
      \Cat{M} \ar{r}{} \ar[bend left]{rr}{\id} & \mod(R) \ar{r}{U} & \Cat{M}.
    \end{tikzcd}
  \end{equation}
 We next show how to explicitly relate $R$-modules according to Definition \ref{definition:module-endo} for an endofunctor $R$ of a finite semisimple category  to modules over the associated quiver $Q(R)$. 

 Let  $R: \Cat{M} \rightarrow \Cat{M}$ be a linear endofunctor of the  semisimple category $\Cat{M}$ and  $f: R(m) \rightarrow m$ a module in $\mod_{m}(R)$ with $m \in \Cat{M}$. Consider the quiver $Q=Q(R)$.
 
The module over $Q$ corresponding to $f$ has  the vector spaces $m_{i}=\Hom_{\Cat{M}}(x_{i},m)$ assigned to
the vertices $i\in I_{\Cat{M}}$ of $Q$.

For a vector $v \in m_{i}$ for an $i \in I$ and a morphism $\lambda \in Q(i,j)=\Hom_{\Cat{M}}(R(x_{i}),x_{j})^{*}$ we define the linear map
$f(\lambda): m_{i}\rightarrow m_{j}$ as follows: First compose
\begin{equation}
  \label{eq:19}
  f \circ R(v) \in \Hom_{\Cat{M}}(R(x_{i}),m) \stackrel{~\eqref{eq:ses-otik}}{\cong}  \oplus_{j \in I} \Hom(R(x_{i}),x_{j}) \otik \Hom(x_{j},m).
\end{equation}
Composing then  with $\lambda \otik \id$ provides an element
\begin{equation}
  \label{eq:alg-f}
  f(\lambda)(v)= (\lambda \otimes \id) \circ f \circ R(v) \in \Hom_{\Cat{M}}(x_{j},m)=m_{j},
\end{equation}
as required. 

For later computations it is useful to reformulate this construction using graphical calculus for categories, see Appendix \ref{sec:categorical-notions} for a short review.

The element $v \in m_{i}$ is expressed graphically as $
\begin{tikzpicture}[very thick,scale=1,color=blue!50!black, baseline]
\draw (0,-0.75) -- (0,1); 
%
%
\fill[color=blue!50!black] (0,0.2) circle (2.9pt) node[right] (v) {{$v$}};
\draw[color=blue!50!black] (-0.1,-0.4) node[right] (A1) {{$i$}};
\draw[color=blue!50!black] (-0.1,0.6) node[right] (A1) {{$m$}};
\end{tikzpicture} 
,
$
and under the identification $\Hom_{\Cat{M}}(R(x_{i}),x_{j})^{*} \stackrel{~\eqref{eq:nat-iso-dual}}{\cong} \Hom_{\Cat{M}}(x_{j},R(x_{i}))$, for $\lambda \in \Hom_{\Cat{M}}(x_{j},R(x_{i}))$
the linear map \eqref{eq:alg-f} has the expression
    \begin{equation}
      \label{eq:f-on-arrow}
      f(\lambda)(v)=
\begin{tikzpicture}[very thick,scale=1,color=blue!50!black, baseline]
\draw (0,-1.5) -- (0,1.5); 
\draw[color=green!50!black] (-1,0) .. controls +(0,0.5) and +(-0.5,-0.5) .. (0,1);
\draw[color=green!50!black] (-1,0) .. controls +(0,-0.5) and +(-0.5,0.5) .. (0,-1);
%
%
\fill[color=blue!50!black] (0,1) circle (2.9pt) node[right] (meet1) {{$f$}};
\fill[color=blue!50!black] (0,0.2) circle (2.9pt) node[right] (v) {{$v$}};
\fill[color=blue!50!black] (0,-1) circle (2.9pt) node[right] (meet2) {{$\lambda$}};
\draw (0,-2.1) node[above] (X) {{$j$}};
\draw[color=green!50!black] (-0.7,-0.56) node[left] (A1) {{$R$}};
\draw[color=blue!50!black] (-0.1,-0.4) node[right] (A1) {{$i$}};
\draw[color=blue!50!black] (-0.1,0.6) node[right] (A1) {{$m$}};
\draw[color=blue!50!black] (-0.1,1.35) node[right] (A1) {{$m$}};
\end{tikzpicture} 
\, . 
    \end{equation}

More generally on an arrow $\omega$ of length $n$, i.e. a morphism $\omega: x_{j} \rightarrow R^{n}(x_{i})$  for some $i,j \in I$, the corresponding morphism
     $f(\omega): m_{i} \rightarrow  m_{j}$,  is defined  on
       $ v\in m_{i}$ as
       \begin{equation}
         \label{eq:f-on-paths}
         f(\omega)(v)=
\begin{tikzpicture}[very thick,scale=1,color=blue!50!black, baseline]
\draw (0,-2.5) -- (0,1.5); 
\draw[color=green!50!black] (-0.5,-0.75) .. controls +(0,0.5) and +(-0.5,-0.5) .. (0,0.15); 
\draw[color=green!50!black] (-1,-0.75) .. controls +(0,0.5) and +(-0.5,-0.5) .. (0,0.8);
\draw[color=green!50!black] (-1.3,-0.75) .. controls +(0,0.5) and +(-0.5,-0.5) .. (0,1.1);
\draw[color=green!50!black] (-0.5,-0.75) .. controls +(0,-0.5) and +(0,0.5) .. (-0.2,-1.5); 
\draw[color=green!50!black] (-1,-0.75) .. controls +(0,-0.5) and +(0,0.5) .. (-0.8,-1.5);
\draw[color=green!50!black] (-1.3,-0.75) .. controls +(0,-0.5) and +(0,0.5) .. (-1,-1.5);
%
\fill[color=blue!50!black] (0,0.15) circle (2.9pt) node[right] (meet2) {{$f$}};
\fill[color=blue!50!black] (0,0.8) circle (2.9pt) node[right] (meet2) {{$f$}};
\fill[color=blue!50!black] (0,1.1) circle (2.9pt) node[right] (meet2) {{$f$}};
\fill[color=blue!50!black] (0,-0.5) circle (2.9pt) node[right] (meet2) {{$v$}};
\draw (-0.5,-1.7) node[minimum height=0.5cm,minimum width=1cm,draw,fill=white] {{$\omega$}};
%
\draw (-0.1,-2.2) node[right] (X) {{$j$}};
\draw (-0.1,-1) node[right] (X) {{$i$}};
\draw[color=green!50!black] (-0.18,-0.1) node[above] (A1) {{$\vdots$}};
\draw[color=green!50!black] (-0.27,-0.75) node[above] (A1) {{$R$}};
\draw[color=green!50!black] (-0.43,-0.27) node[above] (A1) {{$R$}};
\draw[color=green!50!black] (-1,0) node[above] (A1) {{$R$}};
\end{tikzpicture} 
\, . 
       \end{equation} 
       With the notation of Equation \eqref{eq:f-ast},  $f(\omega)(v)= f^{\ast n} \circ R^{n}(v) \circ \omega$.
       \begin{proposition}
         \label{proposition:modules-endof}
   For an endofunctor $R \in \End(\Cat{M})$, there are canonical equivalences of   categories of modules
   \begin{equation}
     \label{eq:11}
   \mod(R)  \cong  \mod(Q(R)), \quad \text{and} \quad \modln(R) \cong \modln(Q(R)).    
   \end{equation}
 \end{proposition}
 \begin{proof}
   For the first part we have by semisimplicity for an object $m \in \Cat{M}$ the isomorphism
  \begin{equation}
    \label{eq:23}
    \mod_{m}(R) \cong \oplus_{i,j} \Hom_{\Cat{M}}(R(x_{i}),x_{j}) \otik \Hom_{\vect}(m_{i},m_{j}).
  \end{equation}
  On the other hand,  also
  \begin{equation}
    \label{eq:24}
    \begin{split}
      \mod_{m}(Q(R))= & \oplus_{i,j}   Q(R)(i,j)^{*}  \otik \Hom_{\vect}(m_{i},m_{j}) \\
      \cong & \oplus_{i,j} \Hom_{\Cat{M}}(R(x_{i}),x_{j})^{**} \otik \Hom_{\vect}(m_{i},m_{j}).  
    \end{split}
     \end{equation}
     With the canonical isomorphism $V^{**} \cong V$ for a finite dimensional vector space $V$,  we obtain the  isomorphism $\mod_{m}(R) \cong \mod_{m}(Q(R))$ for all $m \in \Cat{M}$. Clearly these isomorphisms are functorial and provide an equivalence of categories  $\mod(R) \cong \mod(Q(R))$.
     Expressed on elements, this isomorphism reduces to the one defined by  Equation \eqref{eq:f-on-arrow}.

  Suppose now that $f \in \mod(Q)$ is locally nilpotent. Thus we find $n \in \N$ such that for all paths $\omega \in \PQ_{n}$, $f(\omega)=0$. By Equation \eqref{eq:f-on-paths}, this implies that
  $f^{\ast n} \circ R(Q)^{n}(v) \circ \omega=0$ for all $i,j$ and  $\omega \in \Hom_{\Cat{M}}(x_{j},R(Q)^{n}(x_{i}))$ and all $v \in \Hom_{\Cat{M}}(x_{i},m)$. But this is by semisimplicity equivalent  to
  $f ^{\ast n} \circ g=0$ for all $g \in \Hom_{\Cat{M}}(x_{j}, R(Q)^{n}(m))$ and all $j \in I$. Again by semisimplicity, we conclude $f^{\ast n}=0$ in $\mod(R(Q))$.

  Conversely, if $f^{\ast n}=0$ in $\mod(R)$, it follows directly from Equation  \eqref{eq:f-on-paths}, that $f(\omega)=0$ for all $\omega \in \PQ_{n}$, and thus $f$ is locally nilpotent as a $Q(R)$-module.  
\end{proof}

Under the equivalence of the proposition, the action \eqref{eq:action-gauge} of $g \in \Gl(m)$ on $f \in \mod_{m}(R)=\Hom_{\Cat{M}}(R(m),m)$ is given as
\begin{equation}
  \label{eq:Gl-action-Fun}
  g \act f= g \circ f \circ R(g^{-1}). 
\end{equation}
\begin{remark}
  It is natural to ask, whether for an endofunctor $R: \Cat{M} \rightarrow \Cat{M}$ there is an analogue of the path algebra of a quiver. That is, is there an algebra $A_{R}$ in $\Vect$, such that the category of left modules over $A_{R}$ is equivalent to the category $\mod(R)$ as in Definition \ref{definition:module-endo}?
  We believe that such an algebra can be constructed explicitly from a progenerator $P$ for $\Cat{M}$, namely
  \begin{equation}
    \label{eq:gen-alg}
    A_{R}= \oplus_{n \in \N} \End_{\Cat{M}}(P, R^{n}(P)),
  \end{equation}
  with multiplication given by the composition
  \begin{equation*}
  \End(P, R^{n}(P)) \otik \End(P,R^{m}(P))\rightarrow \End(P,R^{n+m}(P))    
  \end{equation*}
  with $(f,g) \mapsto  f \circ R^{n}(g)$. This defines a unital associative algebra.
  In case of a finite semisimple $\Cat{M}$, taking $P= \oplus_{i \in I}x_{i}$, this reproduces the path algebra $\PQ(R)$. A proof in the general case is beyond the scope of this paper. 
\end{remark}

\paragraph{A bicategory of endofunctors}

We want to extend Definition \ref{definition:id-morph-lin-Quiv} to allow for  morphisms between quivers with different vertices. By Proposition \ref{proposition:quiv-as-endof}
we identify quivers with endofunctors of finite semisimple categories. 
For general endofunctors we define:

\begin{definition}
  \label{definition:quiver-of-fun-mor-bicat}
  Let $\Cat{M}, \Cat{M}'$ be  linear categories  and $R: \Cat{M} \rightarrow \Cat{M}$ and $R': \Cat{M}' \rightarrow \Cat{M}'$ be linear  endofunctors.
  \begin{definitionlist}
  \item A \emph{morphism} of pairs $(\Cat{M}, R) \rightarrow (\Cat{M}', R')$ is  a tuple $(F, \eta)$ of a functor  $F: \Cat{M} \rightarrow \Cat{M}'$ together  with     a natural transformation  $\eta:R'F \rightarrow FR$.
  \item Let $(F, \eta)$ and $(F', \eta')$ two morphisms of pairs  $(\Cat{M}, R) \rightarrow (\Cat{M}', R')$. A \emph{2-morphism}
    $\Gamma: (F, \eta) \rightarrow (F', \eta')$ is a natural transformation $\Gamma: F \rightarrow F'$, such that
    \begin{equation}
      \label{eq:comm-diag-2-morph}
      \begin{tikzcd}
        R'F    \ar{d}{R'\Gamma}\ar{r}{\eta}  &  FR  \ar{d}{\Gamma R} \\
         R'F'   \ar{r}{\eta'}  &   F'R
      \end{tikzcd}
    \end{equation}
    commutes. 
  \end{definitionlist}
\end{definition}
For $(F,\eta): R \longrightarrow R'$ and $(G, \rho): R' \rightarrow \widetilde{R}$ morphisms of pairs, their composition  is given as the composite functor $GF$ with natural transformation
\begin{equation}
  \label{eq:composition-quiv-1-morph}
  \begin{tikzcd}
    \widetilde{R} G F \ar{r}{\rho F} & G R' F \ar{r}{G \eta} & GF R. 
  \end{tikzcd}
\end{equation}
This definition is entirely abstract, it uses  only that linear categories, linear functors and linear
natural transformation  form  the bicategory $\Lincat$  with the sub bicategory $\lincatses$ of finite semisimple 2-categories.
Indeed, Definition \ref{definition:quiver-of-fun-mor-bicat} as well as the composition \eqref{eq:composition-quiv-1-morph} make sense for any bicategory $\Cat{B}$ and gives the
\emph{bicategory $\Omega(\Cat{B})$ of endomorphisms of $\Cat{B}$}: As objects of  $\Omega(\Cat{B})$ we consider the pairs $(b,f)$ of an object $b \in\Cat{B}$ together with a 1-morphism $f: b \rightarrow b$. As morphisms and 2-morphism between pairs we consider the  structures that are analogous to Definition \ref{definition:quiver-of-fun-mor-bicat} in this abstract setting.
The symbol $\Omega(\Cat{B})$ is chosen due to the similarity of the construction with the loop object construction
in algebraic topology.

\begin{proposition}
  \label{proposition:Bicat-Omega}
  Let $\Cat{B}$ be a bicategory. The endomorphisms of $\Cat{B}$ form a bicategory $\Omega(\Cat{B})$.  
\end{proposition}
\begin{proof}
  To show that $\Omega(\Cat{B})$  is a bicategory we show that  it is a bicategory of 2-functors: Let $\Cat{X}$ be the strict bicategory on one 1-morphism, that is
  $\Cat{X}$ has one object, 1-morphisms are $\{X^{n}\}_{n \in \N}$ and only identity 2-morphisms. The composition of 1-morphisms is the obvious $X^{n} \circ X^{m}=X^{n+m}$.
  Now consider the bicategory of $\Fun^{\mathrm{lax}}(\Cat{X},\Cat{B})$ of strict 2-functors, lax 2-transformations and modifications, see Subsection \ref{sec:bicategories}. 
  A strict 2-functor $F: \Cat{X} \rightarrow \Cat{B}$ is the same as an object $b \in \Cat{B}$ together with a 1-morphism $R: b\rightarrow b$, which is the image of $X$.
  Similarly, a lax 2-transformation is the same as a morphism in $\Omega(\Cat{B})$ and a modification the same as a 2-morphism in $\Omega(\Cat{B})$.
  As also the compositions agree, we conclude that $\Fun^{\mathrm{lax}}(\Cat{X},\Cat{B})$ and $\Omega(\Cat{B})$ are equivalent, in particular is  $\Omega(\Cat{B})$ a bicategory.  
\end{proof}

For later use we record, that there is an evident forgetful 2-functor
\begin{equation}
  \label{eq:forget-2-fun}
  \Omega(\Cat{B}) \rightarrow \Cat{B},
\end{equation}
which forgets the additional 1-morphism on the objects and similar for the 1- and 2-morphisms.

We finally arrive at
\begin{definition}
  \label{definition:bicat-quiv}
  The \emph{bicategory of finite quivers} is $\Quiv=\Omega(\lincatses)$, while the category of generalized quivers is $\Omega(\Lincat)$.  
\end{definition}
Obviously,  $\Omega(\lincatses)$ is a sub bicategory of $\Omega(\Lincat)$. A morphism of pairs $(F, \eta): Q \rightarrow Q'$ in $\Quiv$ with $F=\id$ reduces to a morphism in $\Quiv_{I}$ from Definition \ref{definition:id-morph-lin-Quiv}.  
\begin{remark}
  \label{remark:morph-give-fun-Q}
  Let $(F, \eta): Q \rightarrow Q'$ be a 1-morphism of quivers $Q \in \End(\Cat{M})$ and $Q'\in \End(\Cat{M}')$.
  It induces a functor
  \begin{equation}
    \label{eq:eta-fun}
    \widehat{\eta}: \mod(Q) \rightarrow \mod(Q'),
  \end{equation}
  with components $\widehat{\eta}_{m}: \mod_{m}(Q) \rightarrow \mod_{F(m)}(Q')$ given as the composite
  \begin{equation}
    \label{eq:13}
    \widehat{\eta}_{m}: \Hom_{\Cat{M}}(Qm,m) \stackrel{F(-)}{\longrightarrow} \Hom_{\Cat{M}}(FQ(m),F(m))
    \stackrel{\eta_{m}^{*}}{\longrightarrow} \Hom_{\Cat{M}}(Q'F(m),F(m))=\mod_{F(m)}(Q'). 
  \end{equation}
  On a morphism, $h: (m,f) \rightarrow (m',f')$ in $\mod(Q)$, it is defined as $\widehat{\eta}(h)=F(h)$. By naturality of $\eta$, this is indeed a morphism from $\widehat{\eta}(m,f)$ to $\widehat{\eta}(m',f')$.

  Note that the objects of $\mod(Q)$ form a vector space and the functor $\widehat{\eta}$ has the additional property of being linear on the objects. However, not all functors between quiver modules with this property come from a morphism
  in $\Omega(\Lincat)$, an example is the action \eqref{eq:Gl-action-Fun} of the groups $\Gl(m)$: Picking for all objects $m \in \Cat{M}$ an object $g_{m}\in \Gl(m)$, provides a functor
  $\Ad_{g}: \mod(Q) \rightarrow \mod(Q)$ sending a module $(m,f)$ to $(m,g \act f) \in \mod_{m}(Q)$ and a morphism $h: (m,f) \rightarrow (m',f')$ to the morphism $g_{m'}\circ h \circ g_{m}$. It is easy to see that this defines a
  functor with is also linear on objects, but i.g. there exists no morphism $(F, \eta): Q \rightarrow Q'$  in $\Quiv$ with $\widehat{\eta}=\Ad_{g}$. 
\end{remark}

\section{Fusion categories and their actions on quivers}
\label{sec:fusi-categ-their}
We now consider actions of monoidal categories on categories of quiver modules. See Appendix \ref{sec:mono-categ-their} for a short review on module categories and related structures.

By Definition \ref{definition:gen-quiv},
a generalized  quiver is   a linear endofunctor $Q: \Cat{M} \rightarrow \Cat{M}$ on a  linear category $\Cat{M}$. By Proposition \ref{proposition:Bicat-Omega}, the endomorphism category of an object $(\Cat{M}, Q)$ of $\Omega(\Lincat)$, denoted $\End_{\Omega}(\Cat{M},Q)$ is a monoidal category. Thus we define
\begin{definition}
  \label{definition:Action-quiv}
  Let $\Cat{A}$ be a monoidal category, $(\Cat{M},Q)$ a generalized quiver. 
  A \emph{left action of $\Cat{A}$ on $Q$} is a strong monoidal functor  
  \begin{equation}
    \label{eq:left-action-Q}
    L: \Cat{A} \rightarrow \End_{\Omega}(\Cat{M},Q).  
  \end{equation}

  Analogously, a right action of $\Cat{A}$ on $Q$, and a bimodule action of monoidal categories $\Cat{A}$, $\Cat{B}$
  on $Q$ are defined as strong monoidal functors from $\Cat{A}^{\rev}$ and $\Cat{A} \times \Cat{B}^{\rev} $ on $\End_{\Omega}(\Cat{M}, Q)$, respectively. Here $\Cat{B}^{\rev}$ is the category $\Cat{B}$ with reversed monoidal product. 
\end{definition}
A quiver $Q$  with an $\Cat{A}$-action on $Q$ is also called a \emph{$\Cat{A}$-left module quiver} and analogously we speak of  $\Cat{A}$-right module quivers and $(\Cat{A},\Cat{B})$-bimodule quivers.

To define morphisms between $\Cat{A}$-left module quivers, we recall that a monoidal category $\Cat{A}$ can be equivalently regarded as a bicategory with one object, see Appendix \ref{sec:bicategories} for a review.

In this interpretation, a strong monoidal functor \eqref{eq:left-action-Q} is the  same as a 2-functor $\Cat{A} \rightarrow \Omega(\Lincat)$ with sends the single object of the bicategory $\Cat{A}$
to the object $(\Cat{M},Q)$ of $\Omega(\Lincat)$. This allows to define the \emph{bicategory of $\Cat{A}$-left module quivers} as
\begin{equation}
  \label{eq:15}
  {}_{\Cat{A}}\Quiv= \twoFunol(\Cat{A}, \Omega(\Lincat)).
\end{equation}

We proceed to unpack the definition of a $\Cat{A}$-left module quiver, using  the  bicategory  ${}_{\Cat{A}}\Modoplax$
of $\Cat{A}$-module  categories with oplax module functors, see Definition \ref{definition:bicat-mod-cat}.
\begin{proposition}
  \label{proposition:gen-quiv-A}
  Let $\Cat{A}$ be a monoidal category and $(\Cat{M},Q)$ a generalized quiver.
  \begin{propositionlist}
    \item\label{item:act-oplax}  An action of $\Cat{A}$ on $Q$ is the same as a $\Cat{A}$-left module category structure on $\Cat{M}$ together with the structure of an oplax 
  $\Cat{A}$-module functor $Q: \AM  \rightarrow \AM$ on $Q$. Likewise, a right (bimodule)-action correspond to right-module  (bimodule)  category $\Cat{M}$ and right-module (bimodule) structure on $Q$. 
\item This  correspondence 
  extends to an  equivalence of bicategories  ${}_{\Cat{A}}\Quiv \cong \Omega({}_{\Cat{A}}\Modoplax)$ and analogously for right-module and bimodule quivers.
  \end{propositionlist}
 \end{proposition}
 \begin{proof}
   For the first part we unpack the structures of an  action of $\Cat{A}$ on $Q$ given by a strong monoidal functor
    $L: \Cat{A} \rightarrow \End_{\Omega}(\Cat{M},Q)$. For each object $a \in \Cat{A}$ it provides a 1-morphism $(L_{a}, \eta_{a}): (\Cat{M},Q) \rightarrow (\Cat{M},Q)$.
  Under the forgetful 2-functor
  $\Omega(\Lincat) \rightarrow \Lincat$ from Equation \eqref{eq:forget-2-fun}, $L$ defines a strong monoidal functor
  $\Cat{A} \rightarrow \End_{\Lincat}(\Cat{M})$, i.e. the structure of a left $\Cat{A}$-module category on $\Cat{M}$.
  The additional natural transformations $\eta_{a}: Q L_{a} \rightarrow  L_{a}Q $ for each object $a \in \Cat{A}$
  provides natural morphisms $ Q(a \act m) \rightarrow a \act Q(m )$ for all $m \in \Cat{M}$.
  Together with the  coherence morphisms of the monoidal functor $L$, these morphisms equip $Q$ with the  structure of an oplax module functor.

 For the second part we use of the bicategory $\Cat{X}$ of the proof of Proposition \ref{proposition:Bicat-Omega}, by which we need to show that there is a canonical equivalence of bicategories
   \begin{equation}
     \label{eq:21}
     \twoFunol(\Cat{A},\twoFunls(\Cat{X},\Cat{B})) \cong \twoFunls(\Cat{X},\twoFunol(\Cat{A},\Cat{B})),
   \end{equation}
   for $\Cat{B}=\Lincat$. However, such an equivalence exists for any bicategory $\Cat{B}$, as can be seen by a straightforward, but lengthy unpackaging of the structures similar to part \refitem{item:act-oplax}.
\end{proof}

By the preceeding proof we obtain for a monoidal category $\Cat{A}$ a forgetful 2-functor
\begin{equation}
  \label{eq:forget-A-mod}
  U:  {}_{\Cat{A}}\Quiv \rightarrow {}_{\Cat{A}}\Modoplax,
\end{equation}
thus, a $\Cat{A}$-module quiver $(\Cat{M},Q)$ has  an  \emph{underlying module category $\AM$}. For $Q_{1}$ and $Q_{2}$ two left $\Cat{A}$-module quivers on $\Cat{M}$ and $\Cat{N}$, respectively, a morphism of $\Cat{A}$-module quivers from $Q_{1}$ to $Q_{2}$ is a module functor $F: \Cat{M} \rightarrow \Cat{N}$ together  with a module natural transformation $F Q_{1} \rightarrow Q_{2} F$.

A  module category $\AM$ is called indecomposable,  if it  is not  equivalent to  a module category
$\AMp \oplus \AMpp$ for non-trivial  module categories $\AMp$ and $\AMpp$.
\begin{definition}
  For module quivers $Q,Q'$, their \emph{sum} is the module quiver $Q \oplus Q'$ corresponding to the sum of the module functors. A module quiver $(\Cat{M},Q)$ is called \emph{indecomposable}, if $\Cat{M}$ is an indecomposable module category and $Q$  is indecomposable as an object in
  $\End_{\Cat{A}}(\Cat{M})$. 
\end{definition}
In the semisimple case we obtain, see Subsection \ref{sec:mono-categ-their} for the notion of a multi-fusion category. 

\begin{proposition}
  Let $\Cat{A}$ be a multi-fusion category.
  There are finitely many indecomposable $\Cat{A}$-module quivers and each $\Cat{A}$-module quiver on a indecomposable module category is a  sum of these. 
\end{proposition}
\begin{proof}
  For a given semisimple 
  module category $\AM$, the category $\End_{\Cat{A}}(\AM)$ is semisimple with finitely many indecomposable $\Cat{A}$-endofunctors of $\AM$, see \cite{Mue1}.
  In total  there are only finitely many indecomposable $\Cat{A}$-module categories by  \cite[Corollary 2.35]{ENOfus} and the statements follow. 
\end{proof}

As a corollary  of Proposition \ref{proposition:gen-quiv-A}  we obtain for a monoidal category $\Cat{A}$ the following description of the  bimodule quivers with underlying bimodule category $\Cat{A}$ seen  as a bimodule category $\AAA$ over itself.
Let $\cent(\Cat{A})$ be the Drinfeld center of $\Cat{A}$. For each $z \in \cent(\Cat{A})$,
 the endofunctor
  \begin{equation}
    \label{eq:Qz}
    Q_{z}= z \otimes -: \Cat{A} \rightarrow \Cat{A}
  \end{equation}
  is a bimodule functor and as recalled in  Equation \ref{eq:Drinf-mod-fun}, all bimodule endofunctors of $\AAA$
  are of this type. Thus we have
\begin{corollary}
  \label{corollary:center-bimod}
  Let $\Cat{A}$ be a monoidal category and $Q: \Cat{A} \rightarrow \Cat{A}$ a  $\Cat{A}$-bimodule quivers with underlying bimodule category $\AAA$. There is an object $z \in \cent(\Cat{A})$, unique up to isomorphism, and an isomorphism of bimodule quivers $Q \cong Q_{z}$. 
Conversely, for each $z \in \cent(\Cat{A})$, the functor $Q_{z}: \Cat{A} \rightarrow \Cat{A}$ is a bimodule quiver. \end{corollary}

Those quivers  play an important role in the sequel and deserve
\begin{definition}
  \label{definition:Drinfeld-quivers}
  Let $\Cat{A}$ be a monoidal category and  $z \in \cent(\Cat{A})$. The generalized quiver $Q_{z}$ given by the endofunctor 
  \begin{equation}
    \label{eq:Qz}
    Q_{z}= z \otimes -: \Cat{A} \rightarrow \Cat{A},
  \end{equation}
is called \emph{generalized Drinfeld quiver of $z$}. In case that $\Cat{A}$ is a multi-fusion category
 $Q_{z}$ is called the \emph{Drinfeld quiver of $z$}. 
\end{definition} 
By construction, a Drinfeld quiver $Q_{z}$ has as set of vertices a labeling set $I$ for the simples of $\Cat{A}$ and arrows $Q_{z}(i,j)=\Hom_{\Cat{A}}(z \otimes x_{i},x_{j})^{*} \cong \Hom_{\Cat{A}}(x_{j}, z \otimes x_{i})$, using the natural isomorphism \eqref{eq:nat-iso-dual}. The sum of Drinfeld quivers $Q_{z}\oplus Q_{y}$ with $z,y \in \cent(\Cat{A})$ is equivalent to the Drinfeld quiver $Q_{z \oplus y}$ for the sum of the  objects.

\begin{remark}
  \label{remark:action-moduli}
  Let $Q: \Cat{M} \rightarrow \Cat{M}$ be  a linear quiver on a module category $\AM$. It follows from
  Remark \ref{remark:morph-give-fun-Q}, that a left action of $\Cat{A}$ on $Q$ provides a collection of linear maps
  \begin{equation}
    \label{eq:action-quiv}
    \varphi_{x;m}: \mod_{m}(Q) \rightarrow \mod_{x\act m}(Q),
  \end{equation}
  for all $x \in \Cat{A}$ and $m \in \Cat{M}$, 
  which are equivariant with respect to the action of $\Gl(m)$ on both representation spaces. Furthermore, the linear maps are coherent: $\varphi_{\unit,m}$ has to be the identity on $\mod_{m}(Q)$ and 
  \begin{equation}
    \label{eq:coher-varphi}
    \begin{tikzcd}
      \mod_{m}(Q) \ar{r}{\varphi_{x;m}} & \mod_{x \act m}(Q) \ar{r}{\varphi_{y,x \act m}} &
      \mod_{y \act (x \act m)}(Q)\cong \mod_{(y \otimes x) \act a }(Q)
    \end{tikzcd}
  \end{equation}
  coincides with $\varphi_{y \otimes x;m}$. Here we used the coherence structure of the
  module category in the last isomorphism.

  The action descends to an action of the Grothendieck ring $\Gr(\Cat{A})$ (see Subsection \ref{sec:mono-categ-their})
on  the categorical quotient $\Mses(Q)$ of Definition \ref{definition:ses-moduli}
  \begin{equation}
    \label{eq:actGr}
    \Gr(\Cat{A}) \times \Mses(Q) \rightarrow \Mses(Q), \quad (a,f) \mapsto a \act f.
  \end{equation}
  The action is homogenous in the sense, that for $f \in \Mses_{m}(Q)$ and $a \in \Gr(\Cat{A})$, we have $a \act f \in \Mses_{am}(Q)$. We thus obtain the structure of a $\Gr(\Cat{A})$-module on the semisimple moduli space $\Mses$. 
\end{remark}

\begin{example}
  \label{example:A-actions}
  \begin{examplelist}
    \item We consider $G$-actions on quivers, where $G$ is a finite group.  A $G$-action  on a quiver $Q$ in our setting is the same as a $\vect_{G}$-action on $Q$, where $\vect_{G}$ is the fusion category of finite-dimensional $G$-graded vector spaces, see Example \ref{example:vectG}.
  The  $\vect_{G}$-module categories and their module endofunctors are classified in
\cite{OsFin}:
  The  indecomposable $\vect_{G}$-module categories correspond to pairs $(H,\Psi)$ of a subgroup $H \subseteq G$ with a class $ \Psi \in H^{2}(H, \Bbbk^{\times})$.

  A pair $(H,\Psi)$ defines an algebra $A(H, \Psi) \in \vect_{G}$, which is the group algebra of $H$ with multiplication twisted by $\Psi$ and the module category $\Cat{M}(H,\Psi)$ corresponding to $(H, \Psi)$ is the category of right $A(H, \Psi)$-modules in $\vect_{G}$ with module action given by the monoidal product of $\vect_{G}$ from the left.
  
  The simple objects in the  category $\End_{\vect_{G}}(\Cat{M})$ given by  the simple  $A(H, \Psi)$-bimodules in $\vect_{G}$ \cite[Prop.3.1]{OsFin}. Thus by Remark \ref{remark:action-moduli}, for any such bimodule $W$ with corresponding quiver $Q_{W}$, the group $G$ acts on $\mod(Q_{W})$ and on the moduli space $\Mses(Q_{W})$.
\item \label{item:Sl2-module}Consider for  $q \in \Bbbk^{\times}$ the following braided monoidal category $\Cat{C}_{q}$: If $q= \pm 1 $ or $q$ is not a root of unity,  $\Cat{C}_{q}=\mod(\SL_{q}(2))$ are the finite dimensional $\SL_{q}(2)$-modules.
In the other  cases, $\Cat{C}_{q}$ is the semisimple subquotient of $\mod(\SL_{q}(2))$. 
In each case, $\Cat{C}_{q}$ is generated by the  fundamental representation $V$. For any connected quiver $Q$ with vertices $I$, there exists by \cite{EOGraphs} a $q \in \Bbbk^{\times}$ and an action of $\Cat{C}_{q}$ on $\Cat{M}=\vect^{I}$, such that the quiver $Q_{V}=V \act -: \Cat{M} \rightarrow \Cat{M}$ is isomorphic to $Q$.
  Since $\Cat{C}_{q}$ is braided, it follows that any connected $Q$ has the structure of a  $\Cat{C}_{q}$-module quiver for some $q$. 
\end{examplelist}
\end{example}

\section{Fusion quivers}

\label{sec:tensor-quivers}

We consider quivers with a rigid monoidal structure on their categories of modules. 
For a  rigid  monoidal  
category   $\Cat{A}$ and an object  $z \in \cent(\Cat{A})$ we show that the Drinfeld quiver $Q_{z}$ from Definition \ref{definition:Drinfeld-quivers} naturally possesses  such a structure.
We discuss explicit examples of such quivers $Q_{z}$ and examine induced structures on the moduli space  of semisimple modules..

We propose that the following defines an interesting class of quivers:
\begin{definition}
  \label{definition:fus-quiver}
  A \emph{generalized fusion quiver} is a  generalized quiver, i.e. a linear functor $Q: \Cat{A} \rightarrow \Cat{A}$ on a linear category $\Cat{A}$,
  together with  a linear monoidal structure $\qoti$ on $\mod(Q)$, such that
  \begin{definitionlist}
  \item The monoidal category $(\mod(Q), \qoti)$ is rigid monoidal,
  \item the semisimple subcategory $\Cat{A}_{Q}=\mod(Q)^{\vses}$ of vertex semisimple modules (see Subsection \ref{sec:finite-quivers}) is a rigid monoidal subcategory.
  \end{definitionlist}
 A \emph{fusion quiver} is a generalized fusion quiver, such that $\Cat{M}$ is finitely semisimple.
    A (generalized) fusion quiver is called \emph{non-trivial}, if $Q$ is not the zero morphism in $\End(\Cat{M})$. 
   \end{definition}
   We call the rigid monoidal category $\Cat{A}_{Q}$ the underlying monoidal category of $Q$.
By \eqref{eq:retract-gen}, as a linear category, $\Cat{A}_{Q}$ is identified with $\Cat{A}$. 
   If $Q$ is a fusion quiver, $\Cat{A}_{Q}$ is a multi-fusion category. 
 A fusion quiver is non-trivial if and only if the corresponding combinatorial quiver is non-trivial, i.e. has a non-empty set of arrows.

 \begin{definition}
   Let $Q: \Cat{A} \rightarrow \Cat{A}$ and $Q': \Cat{A}' \rightarrow \Cat{A}'$ be generalized fusion quivers. A \emph{1-morphism of fusion quivers} is a strong monoidal functor
   $F: \mod(Q) \rightarrow \mod(Q')$ which restricts to a monoidal functor between the underlying monoidal categories of vertex semisimple modules.
 \end{definition}
 It is clear, that together with monoidal natural transformations as 2-morphisms, fusion quivers and 1- and 2-morphisms of fusion quivers form a bicategory $\fusqui$. 
\subsection{Fusion quiver structures on Drinfeld quivers}
\label{sec:fusi-quiv-struct}

We provide a  canonical construction of generalized fusion quivers for the generalized Drinfeld quivers of Definition \ref{definition:Drinfeld-quivers}:

Let $\Cat{A}$ be a  monoidal category and $z \in \cent(\Cat{A})$ with linear  functor $Q_{z}= z \otimes -: \Cat{A} \rightarrow \Cat{A}$.  
We abbreviate the category $\mod(Q_{z})=\Cat{C}_{z}$, which by Definition \ref{definition:module-endo} has  objects tuples $(a,f)$ where $a \in \Cat{A}$ and $f: z \otimes  a \rightarrow a$ is a morphism in $\Cat{A}$. 
If the object $a \in \Cat{A}$ is clear from the context, we denote the object just as $f$.

To define the fusion of two $Q_{z}$-modules, we use the diagrammatic calculus, see Subsection \ref{sec:mono-categ-their}.
For two objects $(a,f)$ and $(b,g)$ of $\Cat{C}_{z}$ we define the object $(a \otimes b, f \qoti g)$ in $\Cat{C}_{z}$ as the object
$a \otimes b$ in $\Cat{A}$ together with the morphism
\begin{equation}
  \label{eq:qotie}
  f \qoti g: z \otimes ( a \otimes b) \rightarrow a \otimes b
\end{equation}
defined as
\begin{equation}
  \label{eq:def-qoti}
 f \qoti  g =
\begin{tikzpicture}[very thick,scale=1,color=blue!50!black, baseline]
    \draw[color=blue!50!black] (-2.0,-1.5) node[above] (z) {{$z$}};
  \draw[color=blue!50!black]   (-1.2,-1.5) node[above] (a) {{$a$}};
  \draw[color=blue!50!black]   (-0.6,-1.5) node[above] (b) {{$b$}};
\draw[color=blue!50!black] (-1.2,1) -- (-1.2,-1); 
\draw[color=blue!50!black] (-0.6,1) -- (-0.6, -1); 
\draw[color=green!50!black] (-2,-1) .. controls +(0,0.5) and +(-0.5,-0.5) .. (-1.2,0.25);  
\fill[color=blue!50!black] (-1.2,0.25) circle (2.9pt) node[right, xshift=-1.8] (meet2) {{$f$}};
\end{tikzpicture} \;
+ 
\begin{tikzpicture}[very thick,scale=1,color=blue!50!black, baseline]
    \draw[color=blue!50!black] (-1.9,-1.5) node[above] (z) {{$z$}};
  \draw[color=blue!50!black] (-1.2,-1.5) node[above] (a) {{$a$}};
  \draw[color=blue!50!black] (-0.6,-1.5) node[above] (b) {{$b$}};
\draw[color=blue!50!black] (-0.6,1) -- (-0.6,-1); 
\draw[color=blue!50!black] (-1.2,1) -- (-1.2,-1); 
\draw[color=white, line width=4pt] (-1.9,-1) .. controls +(0,0.5) and +(-0.5,-0.5) .. (-0.6,0.25);  
\draw[color=green!50!black] (-1.9,-1) .. controls +(0,0.5) and +(-0.5,-0.5) .. (-0.6,0.25);  
\fill[color=blue!50!black] (-0.6,0.25) circle (2.9pt) node[right, xshift=-1.8] (meet2) {{$g$}};
\end{tikzpicture} 
\end{equation}
using the braiding of $z$ in the second term. As usual in the diagrammatic calculus, the associator of $\Cat{A}$ is not depicted, but the
evaluation of a diagram makes sense also for non-strict monoidal categories: For a choice of bracketing of the external strands, a diagram evaluates to a unique morphism.

From Corollary \ref{corollary:center-bimod}, $\Cat{C}_{z}$ is  a bimodule category over $\Cat{A}$ with  actions
denoted $b \act f$ and $f \ract b$ for $(a,f) \in \Cat{C}_{z}$ and $b \in \Cat{A}$, and we conclude, that by definition of $\qoti$,
\begin{equation}
  \label{eq:qoti-bimod}
  f \qoti g= f \ract b + a \act  g.
\end{equation}
For morphisms $h_{1}:(a,f) \rightarrow (a',f')$ and $h_{2}:(b,g) \rightarrow (b',g')$
in $\Cat{C}_{z}$, we set
\begin{equation}
  \label{eq:qoti-morph}
  h_{1} \qoti h_{2}= h_{1} \otimes h_{2}: a \otimes b  \rightarrow a' \otimes b'. 
\end{equation}

\begin{lemma}
  The morphism $h_{1} \qoti h_{2}$ in $\Cat{A}$ is a morphism from $f \qoti g$ to $f' \qoti g'$ in $\Cat{C}_{z}$. 
\end{lemma}
\begin{proof}
  This follows most easily in the graphical calculus:
  \begin{equation}
    \label{eq:Pic2}
\begin{tikzpicture}[very thick,scale=1,color=blue!50!black, baseline]
    \draw[color=blue!50!black] (-2.0,-1.5) node[above] (z) {{$z$}};
  \draw[color=blue!50!black]   (-1.2,-1.5) node[above] (a) {{$a$}};
  \draw[color=blue!50!black]   (-0.6,-1.5) node[above] (b) {{$b$}};
\draw[color=blue!50!black] (-1.2,1) -- (-1.2,-1); 
\draw[color=blue!50!black] (-0.6,1) -- (-0.6, -1); 
\draw[color=green!50!black] (-2,-1) .. controls +(0,0.5) and +(-0.5,-0.5) .. (-1.2,0.25);  
\fill[color=blue!50!black] (-1.2,0.25) circle (2.9pt) node[right, xshift=-1.8] (meet2) {{$f$}};
\fill[color=blue!50!black] (-1.2,-0.5) circle (2.9pt) node[right, xshift=-1.8] (mor1) {{$h_{1}$}};
\fill[color=blue!50!black] (-0.6,-0.5) circle (2.9pt) node[right, xshift=-1.8] (mor2) {{$h_{2}$}};
\end{tikzpicture} \;
+ 
\begin{tikzpicture}[very thick,scale=1,color=blue!50!black, baseline]
    \draw[color=blue!50!black] (-1.9,-1.5) node[above] (z) {{$z$}};
  \draw[color=blue!50!black] (-1.2,-1.5) node[above] (a) {{$a$}};
  \draw[color=blue!50!black] (-0.6,-1.5) node[above] (b) {{$b$}};
\draw[color=blue!50!black] (-0.6,1) -- (-0.6,-1); 
\draw[color=blue!50!black] (-1.2,1) -- (-1.2,-1); 
\draw[color=white, line width=4pt]  (-1.9,-1) .. controls +(0,0.5) and +(-0.5,-0.5) .. (-0.6,0.25);  
\draw[color=green!50!black] (-1.9,-1) .. controls +(0,0.5) and +(-0.5,-0.5) .. (-0.6,0.25);  
\fill[color=blue!50!black] (-0.6,0.25) circle (2.9pt) node[right, xshift=-1.8] (meet2) {{$g$}};
\fill[color=blue!50!black] (-1.2,-0.5) circle (2.9pt) node[right, xshift=-1.8] (mor1) {{$h_{1}$}};
\fill[color=blue!50!black] (-0.6,-0.5) circle (2.9pt) node[right, xshift=-1.8] (mor2) {{$h_{2}$}};
\end{tikzpicture} 
=
\begin{tikzpicture}[very thick,scale=1,color=blue!50!black, baseline]
    \draw[color=blue!50!black] (-2.0,-1.5) node[above] (z) {{$z$}};
  \draw[color=blue!50!black]   (-1.2,-1.5) node[above] (a) {{$a$}};
  \draw[color=blue!50!black]   (-0.6,-1.5) node[above] (b) {{$b$}};
\draw[color=blue!50!black] (-1.2,1) -- (-1.2,-1); 
\draw[color=blue!50!black] (-0.6,1) -- (-0.6, -1); 
\draw[color=green!50!black] (-2,-1) .. controls +(0,0.5) and +(-0.5,-0.5) .. (-1.2,0.25);  
\fill[color=blue!50!black] (-1.2,0.25) circle (2.9pt) node[right, xshift=-1.8] (meet2) {{$f$}};
\fill[color=blue!50!black] (-1.2,0.6) circle (2.9pt) node[right, xshift=-1.8] (mor1) {{$h_{1}$}};
\fill[color=blue!50!black] (-0.6,0.6) circle (2.9pt) node[right, xshift=-1.8] (mor2) {{$h_{2}$}};
\end{tikzpicture} \;
+ 
\begin{tikzpicture}[very thick,scale=1,color=blue!50!black, baseline]
    \draw[color=blue!50!black] (-1.9,-1.5) node[above] (z) {{$z$}};
  \draw[color=blue!50!black] (-1.2,-1.5) node[above] (a) {{$a$}};
  \draw[color=blue!50!black] (-0.6,-1.5) node[above] (b) {{$b$}};
\draw[color=blue!50!black] (-0.6,1) -- (-0.6,-1); 
\draw[color=blue!50!black] (-1.2,1) -- (-1.2,-1); 
\draw[color=white, line width=4pt] (-1.9,-1) .. controls +(0,0.5) and +(-0.5,-0.5) .. (-0.6,0.25);  
\draw[color=green!50!black] (-1.9,-1) .. controls +(0,0.5) and +(-0.5,-0.5) .. (-0.6,0.25);  
\fill[color=blue!50!black] (-0.6,0.25) circle (2.9pt) node[right, xshift=-1.8] (meet2) {{$g$}};
\fill[color=blue!50!black] (-1.2,0.6) circle (2.9pt) node[right, xshift=-1.8] (mor1) {{$h_{1}$}};
\fill[color=blue!50!black] (-0.6,0.6) circle (2.9pt) node[right, xshift=-1.8] (mor2) {{$h_{2}$}};
\end{tikzpicture} 
  \end{equation}
  using in the  first  term, that $h_{1}$ is a morphism in $\Cat{C}_{z}$ and in the second term, that $h_{2}$ is a morphism in $\Cat{C}_{z}$ as well as the naturality of the  braiding. 
\end{proof}

Recall that  the locally nilpotent modules over $Q_{z}$, denoted $\Clnz=\modln(Q_{z})$, correspond to comodules over the  path coalgebra by Proposition \ref{proposition:path-coalg-comod}. 
\begin{theorem}[\cite{MR10}]
  \label{thm:ten-quiv}
  Let $\Cat{A}$ be a linear monoidal category and $z \in \cent(\Cat{A})$.
  \begin{theoremlist}
    \item The category $\Cat{C}_{z}$ of  modules over the quiver $Q_{z}$ with the product
  \begin{equation}
    \label{eq:mon-prod}
    \qoti: \Cat{C}_{z} \times \Cat{C}_{z} \longrightarrow \Cat{C}_{z}
  \end{equation}
  is a linear  monoidal category with unit $(\unit, 0)$, i.e. the $Q_{z}$-module with object $\unit \in \Cat{A}$ together with the zero morphism. Moreover is the monoidal product compatible with the bimodule structure over $\Cat{A}$ in the sense that \eqref{eq:mon-prod} is a balanced  bimodule functor.
  The category $\Clnz$ is a  monoidal subcategory of $\Cat{C}_{z}$.
\item In case that $\Cat{A}$ is rigid,  $\Cat{C}_{z}$ is a generalized fusion quiver with
  $\Clnz$ a  rigid monoidal  subcategory. If $\Cat{A}$ is a multi-fusion category, $\Cat{C}_{z}$ is a fusion quiver. 
  \end{theoremlist}
\end{theorem}
\begin{proof}
  The first part of the result is proven in \cite[Sec.2]{MR10}, see statement in the Introduction, without the perspective on quivers. We repeat the proof in our notation and prove the second part:
  
  We first show that $\qoti$ defines a monoidal structure. We  assume for simplicity  that $\Cat{A}$ is a strict monoidal category and show that then also $\Cat{C}_{z}$ is strict monoidal. In case $\Cat{A}$ is not strict, it is straightforward to see, that the associators and unitors of $\Cat{A}$ provide the corresponding coherence morphisms for $\Cat{C}_{z}$ and $\Clnz$.
  
  The object $(\unit, 0)$ is the unit  object, since with the $0$-morphism in 
  Equation  \eqref{eq:qoti-bimod} we deduce $f \qoti (\unit,0)=f$ and $(\unit, 0) \qoti f=f$ for all objects $(a,f) \in \Cat{C}_{z}$.
  Given three objects $(a,f),(b,g), (c,h)$ in $\Cat{C}_{z}$, we compute
  \begin{equation}
    \label{eq:Pic3}
      (f \qoti  g) \qoti h =
\begin{tikzpicture}[very thick,scale=1,color=blue!50!black, baseline]
    \draw[color=blue!50!black] (-2.0,-1.5) node[above] (z) {{$z$}};
  \draw[color=blue!50!black]   (-1.2,-1.5) node[above] (a) {{$a$}};
  \draw[color=blue!50!black]   (-0.6,-1.5) node[above] (b) {{$b$}};
    \draw[color=blue!50!black]   (0,-1.5) node[above] (c) {{$c$}};

\draw[color=blue!50!black] (-1.2,1) -- (-1.2,-1); 
\draw[color=blue!50!black] (-0.6,1) -- (-0.6, -1); 
\draw[color=blue!50!black] (0,1) -- (0, -1); 
\draw[color=green!50!black] (-2,-1) .. controls +(0,0.5) and +(-0.5,-0.5) .. (-1.2,0.25);  
\fill[color=blue!50!black] (-1.2,0.25) circle (2.9pt) node[right, xshift=-1.8] (meet2) {{$f$}};
\end{tikzpicture} \;
+ 
\begin{tikzpicture}[very thick,scale=1,color=blue!50!black, baseline]
    \draw[color=blue!50!black] (-1.9,-1.5) node[above] (z) {{$z$}};
  \draw[color=blue!50!black] (-1.2,-1.5) node[above] (a) {{$a$}};
  \draw[color=blue!50!black] (-0.6,-1.5) node[above] (b) {{$b$}};
      \draw[color=blue!50!black]   (0,-1.5) node[above] (c) {{$c$}};
\draw[color=blue!50!black] (-0.6,1) -- (-0.6,-1); 
\draw[color=blue!50!black] (-1.2,1) -- (-1.2,-1); 
\draw[color=blue!50!black] (0,1) -- (0, -1); 
\draw[color=white, line width=4pt] (-1.9,-1) .. controls +(0,0.5) and +(-0.5,-0.5) .. (-0.6,0.25);  
\draw[color=green!50!black] (-1.9,-1) .. controls +(0,0.5) and +(-0.5,-0.5) .. (-0.6,0.25);  
\fill[color=blue!50!black] (-0.6,0.25) circle (2.9pt) node[right, xshift=-1.8] (meet2) {{$g$}};
\end{tikzpicture} 
+
\begin{tikzpicture}[very thick,scale=1,color=blue!50!black, baseline]
    \draw[color=blue!50!black] (-2.0,-1.5) node[above] (z) {{$z$}};
  \draw[color=blue!50!black]   (-1.2,-1.5) node[above] (a) {{$a$}};
  \draw[color=blue!50!black]   (-0.6,-1.5) node[above] (b) {{$b$}};
    \draw[color=blue!50!black]   (0,-1.5) node[above] (c) {{$c$}};

\draw[color=blue!50!black] (-1.2,1) -- (-1.2,-1); 
\draw[color=blue!50!black] (-0.6,1) -- (-0.6, -1); 
\draw[color=blue!50!black] (0,1) -- (0, -1); 
\draw[color=white, line width=4pt] (-1.9,-1) .. controls +(0,0.5) and +(-0.5,-0.5) .. (0,0.25);  
\draw[color=green!50!black] (-1.9,-1) .. controls +(0,0.5) and +(-0.5,-0.5) .. (0,0.25);  
\fill[color=blue!50!black] (0,0.25) circle (2.9pt) node[right, xshift=-1.8] (meet2) {{$h$}};
\end{tikzpicture}
=      f \qoti  (g \qoti h )\; ,
  \end{equation}
  showing that $\qoti$ defines a strict monoidal structure on $\Cat{C}_{z}$.

  Next we show that $\Cat{C}_{z}$ is rigid: The right dual of $(a,f)$ is $(a^{*},f^{\dagger})$, where $f^{\dagger}: a^{*} \rightarrow z \otimes a^{*}$ is defined as
  \begin{equation}
    \label{eq:Pic4}
      f^{\dagger}= -
\begin{tikzpicture}[very thick,scale=1,color=blue!50!black, baseline=-0.2cm]
  \draw (0,-0.5) -- (0,0.5);
  \draw[color=blue!50!black] (0,-0.5) .. controls +(0,-0.5) and +(0,-0.5) .. (0.5,-0.5); 
  \draw[color=blue!50!black] (0,0.5) .. controls +(0,0.5) and +(0,0.5) .. (-0.5,0.5);
  \draw (0.5,-0.5) -- (0.5,1);
      \draw (-0.5,0.5) -- (-0.5,-1);
\draw[color=white, line width=4pt] (-1,-1) .. controls +(0,0.5) and +(-0.5,-0.5) .. (0,0);  
\draw[color=green!50!black] (-1,-1) .. controls +(0,0.5) and +(-0.5,-0.5) .. (0,0);  
\fill[color=blue!50!black] (0,0) circle (2.9pt) node[right] (meet2) {{$f$}};
\draw (-0.5,-1.5) node[above] (X) {{$a^*$}};
\draw[color=green!50!black] (-1,-1.5) node[above] (A1) {{$z$}};
\end{tikzpicture} 
\, . 
  \end{equation}
 We show that the evaluation and coevaluation morphisms $\ev{a}: a^{*} \otimes a \rightarrow \unit$ and
  $\coev{a}: \unit \rightarrow a \otimes a^{*}$ of $a$ in $\Cat{A}$ are morphisms in $\Cat{C}_{z}$ and thus equip $(a^{*},f^{\dagger})$ with the structure of a right dual:
  For the evaluation, $\ev{a}: f^{\dagger} \qoti f \rightarrow 0$ (here we consider the  unit object $(\unit,0)$ on the right  hand side) to be a  morphism, we need to show that
  \begin{equation}
    \label{eq:Pic5}
\begin{tikzpicture}[very thick,scale=1,color=blue!50!black, baseline=-0.2cm]
  \draw (0,-0.5) -- (0,0.5);
  \draw[color=blue!50!black] (0,0.5) .. controls +(0,0.5) and +(0,0.5) .. (-0.5,0.5);
  \draw (-0.5,0.5) -- (-0.5,-1);
        \draw (0,-0.5) -- (0,-1);
\draw[color=white, line width=4pt] (-1,-1) .. controls +(0,0.5) and +(-0.5,-0.5) .. (-0.5,0);  
\draw[color=green!50!black] (-1,-1) .. controls +(0,0.5) and +(-0.5,-0.5) .. (-0.5,0);  
\fill[color=blue!50!black] (-0.5,0) circle (2.9pt) node[left] (meet2) {{$f^{\dagger}$}};
\draw (-0.5,-1.5) node[above] (X) {{$a^*$}};
\draw (0,-1.5) node[above] (X) {{$a$}};
\draw[color=green!50!black] (-1,-1.5) node[above] (A1) {{$z$}};
\end{tikzpicture} 
+
\begin{tikzpicture}[very thick,scale=1,color=blue!50!black, baseline=-0.2cm]
  \draw (0,-0.5) -- (0,0.5);
  \draw[color=blue!50!black] (0,0.5) .. controls +(0,0.5) and +(0,0.5) .. (-0.5,0.5);
  \draw (-0.5,0.5) -- (-0.5,-1);
  \draw (0,-0.5) -- (0,-1);
\draw[color=white, line width=4pt] (-1,-1) .. controls +(0,0.5) and +(-0.5,-0.5) .. (0,0);  
\draw[color=green!50!black] (-1,-1) .. controls +(0,0.5) and +(-0.5,-0.5) .. (0,0);  
\fill[color=blue!50!black] (0,0) circle (2.9pt) node[right] (meet2) {{$f$}};
\draw (-0.5,-1.5) node[above] (X) {{$a^*$}};
\draw (0,-1.5) node[above] (X) {{$a$}};
\draw[color=green!50!black] (-1,-1.5) node[above] (A1) {{$z$}};
\end{tikzpicture}
=0 \, .
  \end{equation}
  This follows readily from the triangle identity \eqref{eq:triangle}.  For the coevaluation, it follows analogously that
  \begin{equation}
    \label{eq:Pic6}
\begin{tikzpicture}[very thick,scale=1,color=blue!50!black, baseline=-0.2cm]
  \draw (0,-0.5) -- (0,0.5);
  \draw[color=blue!50!black] (0,-0.5) .. controls +(0,-0.5) and +(0,-0.5) .. (-0.5,-0.5);
  \draw (-0.5,-0.5) -- (-0.5,1);
        \draw (0,0.5) -- (0,1);
\draw[color=white, line width=4pt] (-1,-1) .. controls +(0,0.5) and +(-0.5,-0.5) .. (-0.5,0);  
\draw[color=green!50!black] (-1,-1) .. controls +(0,0.5) and +(-0.5,-0.5) .. (-0.5,0);  
\fill[color=blue!50!black] (-0.5,0) circle (2.9pt) node[left] (meet2) {{$f$}};
\draw (-0.8,0.2) node[above] (X) {{$a$}};
\draw (0.3,0.2) node[above] (X) {{$a^{*}$}};
\draw[color=green!50!black] (-1,-1.5) node[above] (A1) {{$z$}};
\end{tikzpicture} 
+
\begin{tikzpicture}[very thick,scale=1,color=blue!50!black, baseline=-0.2cm]
\draw (0,-0.5) -- (0,0.5);
  \draw[color=blue!50!black] (0,-0.5) .. controls +(0,-0.5) and +(0,-0.5) .. (-0.5,-0.5);
  \draw (-0.5,-0.5) -- (-0.5,1);
        \draw (0,0.5) -- (0,1);
\draw[color=white, line width=4pt] (-1,-1) .. controls +(0,0.5) and +(-0.5,-0.5) .. (0,0);  
\draw[color=green!50!black] (-1,-1) .. controls +(0,0.5) and +(-0.5,-0.5) .. (0,0);  
\fill[color=blue!50!black] (0,0) circle (2.9pt) node[right] (meet2) {{$f^{\dagger}$}};
\draw (-0.8,0.2) node[above] (X) {{$a$}};
\draw (0.3,0.2) node[above] (X) {{$a^{*}$}};
\draw[color=green!50!black] (-1,-1.5) node[above] (A1) {{$z$}};
\end{tikzpicture}
=0
\, ,
  \end{equation}
  thus $\coev{a}: 0 \rightarrow f \qoti f^{\dagger}$ is a morphism in $\Cat{C}_{z}$. Analogously it follows, that the left dual of $f$ is $(\tensor[^*]{a}{},  \tensor[^\dagger]{f}{})$
  with 
  \begin{equation}
    \label{eq:Pic7}
      \tensor[^\dagger]{f}{}= -
\begin{tikzpicture}[very thick,scale=1,color=blue!50!black, baseline=-0.2cm]
  \draw (0,-0.5) -- (0,0.5);
  \draw[color=blue!50!black,yscale=-1] (0,-0.5) .. controls +(0,-0.5) and +(0,-0.5) .. (0.5,-0.5); 
  \draw[color=blue!50!black,yscale=-1] (0,0.5) .. controls +(0,0.5) and +(0,0.5) .. (-0.5,0.5);
  \draw[yscale=-1] (0.5,-0.5) -- (0.5,1);
      \draw[yscale=-1] (-0.5,0.5) -- (-0.5,-1);
\draw[color=white, line width=4pt] (-1,-1) .. controls +(0,0.5) and +(-0.5,-0.5) .. (0,0);  
\draw[color=green!50!black] (-1,-1) .. controls +(0,0.5) and +(-0.5,-0.5) .. (0,0);  
\fill[color=blue!50!black] (0,0) circle (2.9pt) node[right] (meet2) {{$f$}};
\draw (0.5,-1.65) node[above] (X) {{$ \tensor[^*]{a}{}$}};
\draw[color=green!50!black] (-1,-1.5) node[above] (A1) {{$z$}};
\end{tikzpicture} 
\, .
\end{equation}
  We conclude, that $\qoti$ gives $\Cat{C}_{z}$ the structure of a rigid monoidal category.

  It remains to show that the monoidal structure is compatible with the bimodule structure.
  This is apparent from Equation \eqref{eq:qoti-bimod}: For all objects $(a,f), (b,g) \in \Cat{C}_{z}$ and all $X \in \Cat{A}$, we have
  $(c \act f) \qoti g=c \act (f \qoti g)$ as well as $f\qoti (g \ract c)= (f \qoti g) \ract c$ and $(f  \ract c) \qoti g= f \qoti (c \act g)$. Thus $\qoti$ is a balanced bimodule functor.

  Let us finally show that $\Clnz$ is a rigid monoidal subcategory. If $f,g \in \Clnz$, there exist $n,m \in \N$  with  $f^{\ast n}=0=g^{\ast m}$. By Definition \eqref{eq:def-qoti} of $f \qoti g$ we conclude that $(f \qoti g)^{\ast (n+m)}=0$, since  each term in the sum contains either $\geq n$ times the morphism $f$ or $\geq m$ times the morphism $g$. Moreover, it is easy to see, that also  $(f^{\dagger})^{\ast n}=0$ for $f \in \Clnz$, thus the statement follows. 
\end{proof}
The monoidal structure on $\mod(Q)$ for the  fusion quiver $Q=Q_{z}$  is called the \emph{Drinfeld fusion quiver structure} on $Q_{z}$.

We consider the following morphisms between Drinfeld fusion quivers:
\begin{definition}
  \label{definition:Dfusqui}
  Let $(\Cat{A}, Q_{z})$ and $(\Cat{A}', Q_{z'})$ be generalized Drinfeld fusion quivers with underlying monoidal categories $\Cat{A}$ and $\Cat{A}'$. A \emph{morphisms of generalized Drinfeld fusion quivers} is a  pair $(F, \rho)$ of a strong monoidal functor $F: \Cat{A} \rightarrow \Cat{A}'$ together with a morphism $\rho: z' \rightarrow F(z)$ in $\cent(\Cat{A}')$.  For two morphisms $(F_{1},\rho_{1})$ and $(F_{2},\rho_{2})$ from $(\Cat{A},Q_{z})$ to $(\Cat{A}', Q_{z'})$, a \emph{2-morphism} $\eta:(F_{1},\rho_{1}) \rightarrow (F_{2},\rho_{2})$ is a monoidal natural transformation $\eta: F_{1} \rightarrow F_{2}$, such that $\eta_{z}\circ \rho_{1}=\rho_{2}$. 
\end{definition}
In particular, for $\Cat{A}=\Cat{A}'$, a morphism $\rho: z' \rightarrow z$ defines a morphism of generalized Drinfeld fusion quivers $(\id, \rho)$ from $(\Cat{A}, Q_{z})$ to $(\Cat{A}', Q_{z'})$.

As a consequence, two Drinfeld fusion quivers  $(\Cat{A}, Q_{z})$ and $(\Cat{A}', Q_{z'})$ are equivalent, if there exists a monoidal equivalence  $F: \Cat{A} \rightarrow \Cat{A}'$ and an isomorphism $ z' \cong F(z)$ in $\cent(\Cat{A}')$. 

It follows readily, that this defines a bicategory $\Dfusqui$ of Drinfeld fusion quivers and that by considering the fusion structure on the modules, we obtain a 2-functor
\begin{equation}
  \label{eq:36}
  \Dfusqui\rightarrow \fusqui.
\end{equation}
We consider it as interesting open problem to study the  properties of this 2-functor (is it full, faithful, an equivalence?).

\begin{example}
  \label{example:first-expl}
  \begin{examplelist}
  \item \label{item:gen-Jordan} Let $\Cat{A}$ be a multi-fusion category. For the  choice of $\unit \in \cent(\Cat{A})$ we obtain the fusion quiver $Q_{\unit}=\sqcup_{i \in I}J_{i}$, which is the union of the Jordan quivers $J_{i}$ at the vertex $i$, where $J_{i}$ consists  of
    a  single loop at vertex $i$. An object $(a,f) \in \Cat{C}_{\unit}$ consists of an endomorphism $f: a \rightarrow a$ in $\Cat{A}$ and the monoidal product with $(b,g)$ in $\Cat{C}_{\unit}$ is
    $f \qoti g= f \otimes \id_{b} + \id_{a} \otimes g$.

  With the object  $z=\unit^{\oplus n} \in \cent(\Cat{A})$  we  see that the quiver $Q_{z}$ which consists of $n$ loops at all vertices is also a Drinfeld quiver. 
\item
  \label{item:cyclic-grp}Consider the group $G= \Z /n \Z $ with generator $g$.
For the fusion category $\Cat{A}= \vect_{G}$  with simple objects $e_{g^{i}}$ for $i \in \{1, \ldots, n\}$, see Example \ref{example:vectG} and  $n$-th root of unity $q$ the object $z=e_{g}$ has the  structure of an object in $\cent(\Cat{A})$ with 
with half-braiding given by $c_{z, e_{g^{i}}}: z \otimes e_{g^{i}} \stackrel{q^{i} \cdot \id}{\longrightarrow} e_{g^{i}} \otimes z$.
The corresponding Drinfeld quiver $Q_{z}$  is an oriented cycle on $n$ vertices.
\item  \label{item:ext-Dyn} The double of  an extended Dynkin quiver is a  Drinfeld quiver: Let $G \subset \Sl_{2}(\C)$ be a finite subgroup.
  The category $\Rep(G)$ of finite dimensional $G$-modules is symmetric monoidal and thus each object of $\Rep(G)$
  has a canonical structure of an object in $\cent(\Rep(G))$. For the $G$-representation $V$ which is the fundamental representation of $\Sl_{2}(\C)$ restricted to $G$, the quiver $Q_{V}$ is thus a Drinfeld quiver. 
  By the Mc Kay-correspondence it corresponds to the doubles of   extended Dynkin quivers:
  For an extended Dynkin quiver $D$ with unoriented edges, consider its double with each edge corresponding to two arrows with opposite orientations. The McKay-correspondence  \cite{McKay} states, that the resulting quivers are in bijection with finite subgroups of  $\Sl_{2}(\C)$. 
  \end{examplelist}
\end{example}

By Theorem \ref{thm:ten-quiv} we associate to every object $z \in \cent(\Cat{A})$ a rigid monoidal category $\Cat{C}_{z}$. 
The following functoriality of the  categories $\Cat{C}_{z}$ shows that the construction can be regarded as a categorification of  $\cent(\Cat{A})$, as objects $z$ get mapped to monoidal categories $\Cat{C}_{z}$ and morphisms to monoidal functors.
\begin{lemma} 
  \label{lemma:2-functor-Q}
Let  $(F, \rho)$ be a morphism of generalized Drinfeld fusion quivers from 
$(\Cat{A}, Q_{z})$ to $(\Cat{A}', Q_{z'})$.  It induces a strict monoidal functor
 \begin{equation}
    \label{eq:Morph-Functor}
    (F,\rho): \Cat{C}_{z} \rightarrow \Cat{C}_{z'} 
  \end{equation}
   which sends an object $(a,f)$ to the object $(Fa, (\rho \otimes \id) \circ F(f))$ in $\Cat{C}_{z'}$ and is given by application of $F$ on morphisms.
\end{lemma}
\begin{proof}
  It is clear, that Equation \ref{eq:Morph-Functor} defines  a well-defined functor. It directly follows from the definition, using that $\rho$ is a morphism in the center, that
  $(F, \rho)(f_{1} \qoti f_{2}) \cong (F,\rho)(f_{1}) \qoti (F, \rho)(f_{2})$ for all objects $(a,f_{1}),(b,f_{2}) \in \Cat{C}_{z}$ with the monoidal constraint of $F$ as isomorphism, showing that $(F, \rho)$ is  monoidal.  
\end{proof} 

\begin{remark}
  \label{remark:exotic-asso}
  \begin{remarklist}
  \item At the moment of writing, the author does not know examples of fusion quivers, that are not Drinfeld quivers. Even for a Drinfeld quiver $Q_{z}$  it seems likely that in general  there exist exotic associators for
    the bifunctor $\qoti$, i.e. associators on $\mod(Q_{z})$ that do not arise from the associators of the underlying monoidal category $\Cat{A}$.
    As example consider for $\Bbbk=\C$ the Drinfeld
     quiver
\begin{equation}
  \label{eq:rk1}
J=J_{1}=
\begin{tikzpicture}[thick,scale=1,color=blue!50!black, baseline]
    \fill[color=blue!50!black] (-1.0,0)   circle (2.9pt) node[right, xshift=-1.8] (a) {};
        \draw[->,color=blue!50!black] (-0.9,0.1) .. controls +(1,1) and +(-1,1) .. (-1.1,0.1);  
\end{tikzpicture} \;  
\end{equation}
of Example \ref{example:first-expl}\refitem{item:gen-Jordan} for $\Cat{A}=\vect_{\C}$.
A module $(V,f)$ of $J$ consist of an finite-dimensional vector space $V$ and an endomorphism $f$ of $V$.
The fusion product is given by $(V,f) \qoti (W,g)=(V\otik W, f \otik 1 + 1 \otik g)$. Given a third endomorphism
$(X,h)$, define the isomorphism
$\E^{f \otik g \otik h}: V \otik W \otik X \rightarrow V \otik W \otik X$, regarded as an isomorphism (suppressing the canonical associator of $\vect$)
\begin{equation}
  \label{eq:30}
  \omega_{f,g,h}=\E^{f \otik g \otik h}: (f \qoti g) \qoti h \rightarrow f \qoti (g \qoti h). 
\end{equation}
Indeed, $ \omega_{f,g,h}$ is a morphism in $\mod(J)$, since the exponential commutes with $f \otik 1 \otik 1$, $1 \otik g \otik 1$ and $1 \otik 1 \otik h$. Moreover, it is natural since with $\rho: (V,f) \rightarrow (V',f')$ a morphism, $\rho \circ f = f' \circ \rho$ implies $(\rho \otik 1 \otik 1 )  \E^{f \otik g \otik h} = \E^{f \otik g \otik h} (\rho \otik 1 \otik 1)$ and analogously for the other tensor factors. Moreover, $\omega_{f,g,h}$ satisfies the pentagon for a forth module $(Y,l)$ (omitting $\otik$-symbols)
\begin{equation}
  \label{eq:penta-omega}
  (1_{V} \otik \E^{g \otik h \otik l}) \circ \E^{f \otik (g \qoti h) \otik l} \circ (\E^{f \otik g \otik h} \otik 1_{Y})=
  \E^{f \otik g \otik (h \qoti l)} \circ \E^{(f \qoti g) \otik h \otik l},
\end{equation}
as is seen by a direct computation using $\E^{f \qoti g}= \E^{f \otik 1+ 1 \otik g}= \E^{f \otik 1} \E^{1 \otik g}$.
Since the monoidal unit is $\unit=(\C,0)$, $\omega_{f,\unit,g}=\id$, and with the unitor of $\vect$, the triangle axiom for $\omega$ is satisfied. Thus, $(\mod(J), \qoti, \omega)$ is a monoidal category.

However, it is monoidally equivalent to the standard Drinfeld fusion quiver $(\mod(J), \qoti)$:
For modules $(V,f)$ and $(W,g)$, the isomorphism
\begin{equation}
  \label{eq:2-coboundary}
  \varphi_{f,g}=\E^{- \frac{1}{2}f^{2} \otik g}: V \otik W \rightarrow V \otik W
\end{equation}
is a  natural isomorphism, that satisfies by direct computation, for a third module $(X,h)$:
\begin{equation}
  \label{eq:35}
(1_{V} \otik \varphi_{g,h}) \circ \varphi_{f, g \qoti h} \circ \varphi^{-1}_{f \qoti g, h} \circ (\varphi^{-1}_{f,g} \otik 1_{X})= \omega_{f,g,h}. 
\end{equation}
Thus it follows that the functor $(F, \varphi): (\mod(J), \qoti) \rightarrow (\mod(J), \qoti, \omega)$ with $F=\id$ and
$F(f \qoti g) \stackrel{\varphi_{f,g} \cdot \id}{\longrightarrow} F(f) \qoti F(g)$ is a monoidal equivalence. 
  \item The following is an alternative construction of monoidal structures on quiver modules:
    Starting with a monoidal category and an oplax monoidal functor $Q: \Cat{A} \rightarrow \Cat{A}$ regarded as generalized quiver, the category $\mod(Q)$ is monoidal with monoidal product of $Q$-modules $(a,f)$ and $(b,g)$ with $f: Q(a) \rightarrow a$ and $g: Q(b) \rightarrow b$ morphisms in $\Cat{A}$ given by the composite
    \begin{equation}
      \label{eq:8}
      f \widetilde{\otimes} g: Q(a \otimes b) \longrightarrow Q(a) \otimes Q(b) \stackrel{f \otimes g}{\longrightarrow} a \otimes b, 
    \end{equation}
    using the oplax monoidal structure of $Q$.  However, the category $(\mod(Q), \widetilde{\otimes})$ is in general not rigid: In case that $Q$ is even a strong monoidal functor, it is straightforward to see that an
    object $(a,f) \in \mod(Q)$ has a left and right dual, if and only if $f: Q (a) \rightarrow a $ is an isomorphism in $\Cat{A}$.
  \end{remarklist}
\end{remark}

Which quivers do admit the structure of a fusion quiver? An obstruction for a quiver to be a Drinfeld fusion quiver is the following.
Recall that vertex $i$ of a $Q$ is called a \emph{source} if  no edge of $Q$ ends in $i$ and a  \emph{sink}, if no edge starts in $i$.
\begin{proposition}
  \label{proposition:no-s-s}
  Let $Q=Q_{z}$ be a non-trivial Drinfeld  quiver. Then there are no sources and sinks in $Q$. 
\end{proposition}
\begin{proof}
  The category $\Cat{A}_{Q}=\mod(Q)^{\vses}$ is a semisimple rigid monoidal category. Thus for $z \in \cent(\Cat{A})$ not the zero object, the  product $z \otimes x_{i} \in \Cat{A}$ is not the zero object, so there is for all $i\in I$ some $j \in I$, such that $Q_{z}(i,j)\neq 0$.  Similarly, by considering the product with $z^{*}$, we obtain that for all $i \in I$ there is a $j \in I$ such that $Q_{z}(j,i) \neq 0$, thus there exist not sources or sinks in $Q_{z}$.
\end{proof}
\begin{remark}
  \label{remark:TypesDrinfeld}
  If we consider the Drinfeld quivers $(\Cat{A},z)$ as a method to construct rigid monoidal categories, note that the category $\Cat{C}_{z}$ is non-semisimple and non-finite: Since $Q_{z}$ is not the trivial quiver, it is non-semisimple and since necessarily it has cycles by  Proposition \ref{proposition:no-s-s}, the path algebra is infinite-dimensional and thus $\Cat{C}_{z}$ is non-finite. Moreover we know that $\Cat{C}_{z}$ is hereditary as for all categories of quiver modules the extension groups above degree one vanish, see  \cite{reineke08}. 
\end{remark}
\subsection{Fusion quivers in low rank}
\label{sec:fusion-quivers-low}
 The classification of Drinfeld quivers of a given rank amounts  by Definition \ref{definition:Dfusqui} to the classification of the simple objects in the Drinfeld doubles of all multi-fusion categories of the given rank.
 For low ranks these are known and we list the possible quivers that arise by rank of the underlying multi-fusion category $\Cat{A}$. 

\paragraph{Rank one}

In rank one, there only one fusion category, $\vect$, which has $\cent(\vect)=\vect$ with unique simple object $z= \unit$.
The corresponding indecomposable fusion quiver is  the Jordan quiver
\begin{equation}
  \label{eq:rk1}
J=J_{1}=
\begin{tikzpicture}[thick,scale=1,color=blue!50!black, baseline]
    \fill[color=blue!50!black] (-1.0,0)   circle (2.9pt) node[right, xshift=-1.8] (a) {};
        \draw[->,color=blue!50!black] (-0.9,0.1) .. controls +(1,1) and +(-1,1) .. (-1.1,0.1);  
\end{tikzpicture} \;  
\end{equation}
of Example \ref{example:first-expl}\refitem{item:gen-Jordan} for $\Cat{A}=\vect$.  
The direct sums of $J$ generate all quivers of rank one. Thus, all rank one quivers are fusion quivers.

\paragraph{Rank two}

For rank two, there are four inequivalent fusion categories  $\Cat{A}_{k}$, $k\in \{1, \ldots,4\}$ by
\cite{ostrik2002fusion}.  The categories $\Cat{A}_{1}$ and $\Cat{A}_{2}$ have apart from the unit $\unit$ one object $x$ with the monoidal product determined by  $x^{2}= \unit$ in both cases.
    The categories $\Cat{A}_{3}$ and $\Cat{A}_{4}$ have apart from the unit one object $y$ and the fusion rules $y^{2}=  \unit + y$. The definition of the associators will not be relevant here. 

 The Drinfeld center $\cent(\Cat{A}_{i})$ has in each case four simple  objects, denoted $x_{ij}$ with $i,j \in \{1, \ldots, 4\}$ with $x_{ij} \in \cent(\Cat{A}_{i})$ simple objects for  all $j$.
 From the classification of \cite{ostrik2002fusion} we immediately obtain the following classification of rank two Drinfeld quivers:
 \begin{proposition}
   \label{proposition:rank-two}
   There are 16 inequivalent indecomposable Drinfeld quivers of rank 2 with four possible underlying quivers:
\begin{equation}
\begin{array}{rcl}
   Q^{1}=&
  \begin{tikzpicture}[thick,scale=1,color=blue!50!black, baseline]
    \fill[color=blue!50!black] (-1.0,0)   circle (2.9pt) node[right, xshift=-1.8] (a) {};
    \fill[color=blue!50!black]   (1,0)   circle (2.9pt) node[right, xshift=-1.8] (b) {};
    \draw[->,color=blue!50!black] (-0.9,0.1) .. controls +(1,1) and +(-1,1) .. (-1.1,0.1);  
        \draw[->,color=blue!50!black] (1.1,0.1) .. controls +(1,1) and +(-1,1) .. (0.9,0.1);  
    %
    %
    \draw (-1,0) node[below] (X) {{$\unit$}};
  \end{tikzpicture}  ,\quad
            Q^{2}&=
  \begin{tikzpicture}[thick,scale=1,color=blue!50!black, baseline]
    \fill[color=blue!50!black] (-1.0,0)   circle (2.9pt) node[right, xshift=-1.8] (a) {};
    \fill[color=blue!50!black]   (1,0)   circle (2.9pt) node[right, xshift=-1.8] (b) {};
    \draw[->,color=blue!50!black] (-1,0)+(0.1,0.1) .. controls +(0.5,0.5) and +(-0.5,0.5) .. (1-0.1,0+0.1);  
        \draw[<-,color=blue!50!black] (-1+0.1,0-0.1) .. controls +(0.5,-0.5) and +(-0.5,-0.5) .. (1-0.1,0-0.1);  
       \draw (-1,0) node[below] (X) {{$\unit$}};
     \end{tikzpicture} \; , \\
   Q^{3}=&
  \begin{tikzpicture}[thick,scale=1,color=blue!50!black, baseline]
    \fill[color=blue!50!black] (-1.0,0)   circle (2.9pt) node[right, xshift=-1.8] (a) {};
    \fill[color=blue!50!black]   (1,0)   circle (2.9pt) node[right, xshift=-1.8] (b) {};
    \draw[->,color=blue!50!black] (-1,0)+(0.1,0.1) .. controls +(0.5,0.5) and +(-0.5,0.5) .. (1-0.1,0+0.1);  
        \draw[<-,color=blue!50!black] (-1+0.1,0-0.1) .. controls +(0.5,-0.5) and +(-0.5,-0.5) .. (1-0.1,0-0.1);  
               \draw[->,color=blue!50!black] (1.1,0.1) .. controls +(1,1) and +(1,-1) .. (1.1,-0.1);  
       \draw (-1,0) node[below] (X) {{$\unit$}};
     \end{tikzpicture} \; ,\quad \quad\quad
    Q^{4}&=
  \begin{tikzpicture}[thick,scale=1,color=blue!50!black, baseline]
    \fill[color=blue!50!black] (-1.0,0)   circle (2.9pt) node[right, xshift=-1.8] (a) {};
    \fill[color=blue!50!black]   (1,0)   circle (2.9pt) node[right, xshift=-1.8] (b) {};
    \draw[->,color=blue!50!black] (-1,0)+(0.1,0.1) .. controls +(0.5,0.5) and +(-0.5,0.5) .. (1-0.1,0+0.1);  
        \draw[<-,color=blue!50!black] (-1+0.1,0-0.1) .. controls +(0.5,-0.5) and +(-0.5,-0.5) .. (1-0.1,0-0.1);  
        \draw[->,color=blue!50!black] (-1.1,0.1) .. controls +(-1,1) and +(-1,-1) .. (-1.1,-0.1);  
        \draw[->,color=blue!50!black] (1.1,0.1) .. controls +(1,1) and +(1,-1) .. (1.1,-0.1);  
        \draw (-1,0) node[below] (X) {{$\unit$}};
        \draw (1.8,0.3) node[above] (X) {{$2$}};
\end{tikzpicture} \; ,
\end{array}
\end{equation}
where the ``2'' in $Q^{4}$ indicates  a double arrow. 
    \end{proposition}
    \begin{proof}
      Let $U$ be the respective forgetful functor $U:\cent(\Cat{A}_{i}) \rightarrow \Cat{A}_{i}$ in each case. 
     
      In the centers of $\Cat{A}_{1}$ and $\Cat{A}_{2}$ there are two simples $x_{1j}$ and $x_{2j}$ with $j=\{1,2\}$ such that $U(x_{1j})=\unit$ and $U(x_{2j})=\unit$, and thus the corresponding quivers are $Q^{1}$.
      The  forgetful functor applied to the remaining two simple objects $x_{1j}$ and $x_{2j}$ with $j=3,4$ is the object $x$ in the corresponding categories $\Cat{A}_{1}$ and $\Cat{A}_{2}$, thus the corresponding quivers are $Q^{2}$.
      
  In case of $\Cat{A}_{3}$ and $\Cat{A}_{4}$, again the corresponding centers have  unit objects $x_{i1}$ with $i=3,4$ with corresponding quiver $Q^{1}$ and two objects $x_{i2}$ and $x_{i3}$ for $i=3,4$ with underlying object $y$, which gives the quiver $Q^{3}$.
Finally, the last two objects in the centers have underlying object $y^{2}$, which corresponds to the quiver $Q^{4}$. 
    \end{proof}

    Note that except from $Q^{4}$, all quivers in the above list have thereby different monoidal structures on their modules. Since $Q^{4}\cong Q^{1} \oplus Q^{3}$, its modules can also be equipped with different rigid monoidal structures, using the sum of the Drinfeld quivers. 

    Not all quivers of rank two admit the structure of a Drinfeld quiver. First of all, it is evident, that only symmetric quivers are possible. Even among these, not all are occur: Since it is only possible to take the sum of Drinfeld quivers  with the same underlying monoidal category, the quiver
    \begin{equation}
  \label{eq:3} Q=
  \begin{tikzpicture}[thick,scale=1,color=blue!50!black, baseline]
    \fill[color=blue!50!black] (-1.0,0)   circle (2.9pt) node[right, xshift=-1.8] (a) {};
    \fill[color=blue!50!black]   (1,0)   circle (2.9pt) node[right, xshift=-1.8] (b) {};
    \draw[->,color=blue!50!black] (-1,0)+(0.1,0.1) .. controls +(0.5,0.5) and +(-0.5,0.5) .. (1-0.1,0+0.1);  
        \draw[<-,color=blue!50!black] (-1+0.1,0-0.1) .. controls +(0.5,-0.5) and +(-0.5,-0.5) .. (1-0.1,0-0.1);  
        \draw[->,color=blue!50!black] (1.1,0.1) .. controls +(1,1) and +(1,-1) .. (1.1,-0.1);  
        \draw (-1,0) node[below] (X) {{$\unit$}};
        \draw (0,0.4) node[above] (X) {{$2$}};
         \draw (0,-1) node[above] (X) {{$2$}};
\end{tikzpicture} \;
\end{equation}
does not admit the structure of a Drinfeld quiver, since this would require to take the sum of $Q^{2}$ and $Q^{3}$, which have different underlying monoidal categories.
It is unclear, whether $Q$ admits the structure of  a fusion quiver at all, see Remarks  \ref{remark:exotic-asso} and \ref{remark:use-obst}.

\subsection{Fusion moduli spaces}
\label{sec:fusion-moduli-spaces}
We examine the induced structures on the semisimple quiver moduli space for a fusion quiver $Q$ with underlying multi-fusion category $\Cat{A}$.
Recall from Subsection \ref{sec:finite-quivers}
the moduli space $\Mses_{a}(Q)$ of semisimple quiver modules with dimension vector $a \in \N^{I} $.  The dimension vectors $a \in \N^{I}$ of quiver modules correspond canonically to elements in the   Grothendieck ring $\Gr(\Cat{A})$ of $\Cat{A}$. We will from now on consider  $a \in \Gr(\Cat{A})$ as dimension vector, 
it  at least one basis coefficient of $a$ is negative, we define the set $\Mses_{a}$ to be the empty set. Thus, the  moduli space has a $\Gr(\Cat{A})$ grading   $\Mses=\sqcup_{a\in \Gr(\Cat{A})}\Mses_{a}$. 

From Theorem \ref{thm:ten-quiv} we obtain directly
\begin{corollary}
  Let $Q$ be a fusion quiver on the multi-fusion category $\Cat{A}$. For all dimension vectors $a,b \in \Gr(\Cat{A})$, there is a unique morphism
  $\mu_{a,b}: \Mses_{a} \times \Mses_{b} \rightarrow \Mses_{ab}$, such that
  \begin{equation}
    \label{eq:27}
    \begin{tikzcd}
      \mod_{x_{a}}(Q)  \times \mod_{x_{b}}(Q) \ar{r}{\qoti} \ar{d}{}  & \mod_{x_{a} \otimes x_{b}}(Q) \ar{d}{}\\
      \Mses_{a} \times \Mses_{b} \ar{r}{\mu_{a,b}} & \Mses_{a  b}
    \end{tikzcd}
  \end{equation}
  commutes, where the downward pointing arrows are the morphisms from the categorical quotient and $x_{a},x_{b} \in \Cat{A}$ are objects corresponding to $a,b \in \Gr(\Cat{A})$, respectively. 
  The maps $\{\mu_{a,b}\}$ are associative, i.e. for all $a,b,c \in \Cat{A}$
  \begin{equation}
    \label{eq:asso}
    \mu_{a  b,c} (\mu_{a,b}\times \id) = \mu_{a,  c}(\id \times \mu_{b,c}).
  \end{equation}
\end{corollary}

We conclude, that $\Mses=\sqcup_{a\in \Gr(\Cat{A})}\Mses_{a}$ is an $\Gr(\Cat{A})$-graded ring. 
In the sequel we write $\mu_{a,b}(f, g)=fg$ for $f \in \mod_{a}(Q)$
and $g\in \mod_{b}(Q)$. 

The right duality  $f \mapsto f^{\dagger}$ for $f \in \mod(Q)$ descends to a morphism
\begin{equation}
  \label{eq:29}
  D: \Mses_{a} \rightarrow \Mses_{a^{*}}.
\end{equation}
It follows, that $D$ is an anti-automorphism of $\Mses$, i.e.
with $D(fg)=D(g)D(f)$ with inverse given by the left duality.

\begin{definition}
  The  \emph{ trace $\epsilon_{a}$}  on $\Mses_{a}$ is the morphism
  \begin{equation}
    \label{eq:28}
    \epsilon_{a}: \Mses_{a} \longrightarrow \N,
  \end{equation}
 induced by $\widetilde{\epsilon}(f)= \dim_{\Bbbk}(\Hom_{\mod(Q)}(\unit, f))$
 for $f\in \mod_{x_{a}}(Q)$.
It defines the \emph{pairing}, 
 \begin{equation}
    \label{eq:kappa}
    \kappa_{a, b}: \Mses_{a} \times \Mses_{b}\longrightarrow \N,
  \end{equation}
 as $\kappa_{a,b}= \epsilon_{a b} \circ \mu_{a,b}$. 
\end{definition}

We readily obtain the following list of properties of these structures.
\begin{theorem}
  \label{theorem:Moduli-structures}
  Let $Q$ be a fusion quiver with vertices in $I$, with semisimple moduli spaces  $\Mses$ and underlying multi-fusion category $\Cat{A}$. 
  The $\Gr(\Cat{A})$-graded moduli space $\Mses=\sqcup_{a\in \Gr(\Cat{A})}\Mses_{a}$ 
  has a
  \begin{enumerate}
  \item graded associative and unital multiplication
    \begin{equation}
      \label{eq:2}
      \mu: \Mses \times \Mses \rightarrow \Mses, 
    \end{equation}
  \item graded anti-homomorphism $D: \Mses \rightarrow \Mses$, that is $D(gf)=D(f)D(g)$,
    \item trace $\epsilon: \Mses \rightarrow \N$,  with  induced pairing $\kappa$ of \eqref{eq:kappa}
    \end{enumerate}
    such that
    \begin{enumerate}
 \addtocounter{enumi}{3}
       \item \label{item:eps-Dsquare} $\epsilon(fg)= \epsilon(gD^{2}(f))$ for all $f, g \in \Mses$,
        \item $\epsilon(D^{2}(f))=\epsilon(f)$ for all $f \in \Mses$,
      \item $\epsilon(fD(f)) >0$ for all $f \in \Mses$. 
    \end{enumerate}
\end{theorem} 
\begin{proof}
  The associative multiplication, the anti-homomorphism $D$ and the trace $\epsilon$ are obtained from the monoidal structure, the duality, and the $\Hom$-space with the unit as above. 
  Property \refitem{item:eps-Dsquare} now follows directly from the rigidity  in $\mod(Q)$:
  \begin{equation}
    \label{eq:32}
    \begin{split}
      \Hom_{\mod(Q)}(\unit, f \qoti g)  &\cong    \Hom_{\mod(Q)}(f^{\dagger}, g) \\
                      &\cong     \Hom_{\mod(Q)}(\unit, g \qoti f^{\dagger \dagger}). 
    \end{split}
  \end{equation}
  By specializing to $g=1$ we obtain the second property. The last property follows from $\Hom_{\mod(Q)}(\unit, f \qoti f^{\dagger}) \cong \Hom_{\mod(Q)}(f,f)$,
since  the identity is a non-zero morphism. 
\end{proof}
Since for a multi-fusion category $\Cat{A}$ we have for all $z \in \cent(\Cat{A})$  the fusion  quiver $Q_{z}$, we obtain moduli spaces $\Mses(Q_{z})=\sqcup_{a\in \Gr(\Cat{A})}\Mses(Q_{z})_{a}$ for  all $z \in \cent(\Cat{A})$.
By Lemma \ref{lemma:2-functor-Q}, these are functorial in $z$: A morphism $f: z\rightarrow y$ in $\cent(\Cat{A})$ provides a morphism
$f^{*}: \Mses(Q_{y})\rightarrow \Mses(Q_{z})$, which is easily seen to respect the multiplications. 
\begin{remark}
\label{remark:use-obst}
  It is unclear to the author at the time of writing, how restrictive these structures are for a given quiver $Q$. Can one show that certain  quivers do not admit the structure of a fusion quiver by properties of their semisimple moduli spaces?
\end{remark}

\section{Relations on Drinfeld quivers}
\label{sec:relat-fusi-quiv}
We study relations on the path algebras of Drinfeld quivers, that are compatible with the monoidal product on the
category of modules. These so-called $q$-relations  provide a construction  of quotients of path algebras   whose category of representations has a rigid monoidal structure.

\subsection{Homogenous relations on quivers by morphisms}
In Subsection \ref{sec:finite-quivers} the notion of a quiver with relations is recalled.
One is often interested in the subclass of \emph{homogenous relations}  that are compatible with the grading given by lengths,  see  e.g.  \cite[Thm. 3.14]{Simson}, and are generated by paths  of a fixed length.  

First we discuss a conceptual  approach to homogenous relations for a generalized quiver:
Let $\Cat{A}$ be a monoidal category and $\AM$
a module category. For $x \in \Cat{A}$, consider the generalized quiver $Q_{x}=x \act - : \Cat{M} \rightarrow \Cat{M}$. Let   $\varphi: y \rightarrow x^{\otimes n}$ be a morphism in $\Cat{A}$ for some object $y \in \Cat{A}$ and $n \in \N$ and we define
\begin{definition}
  \label{definition:mod-varphi}
 A module $(m,f) \in \mod(Q_{x})$ \emph{satisfies the relation $\varphi$}, if $(\varphi \act \id_{m}) \circ f^{\ast n}=0$, i.e. 
  \begin{equation}
    \label{eq:varphi-rel}
\begin{tikzpicture}[very thick,scale=1,color=blue!50!black, baseline]
\draw (0,-2) -- (0,1.5); 
\draw[color=green!50!black] (-0.7,-0.75) .. controls +(0,0.5) and +(-0.5,-0.5) .. (0,0.15); 
\draw[color=green!50!black] (-1.3,-0.75) .. controls +(0,0.5) and +(-0.5,-0.5) .. (0,0.8);
\draw[color=green!50!black] (-1.6,-0.75) .. controls +(0,0.5) and +(-0.5,-0.5) .. (0,1.1);
\draw (-1.15,-1) -- (-1.15,-2); 
\fill[color=blue!50!black] (0,0.15) circle (2.9pt) node[right] (meet2) {{$f$}};
\fill[color=blue!50!black] (0,0.8) circle (2.9pt) node[right] (meet2) {{$f$}};
\fill[color=blue!50!black] (0,1.1) circle (2.9pt) node[right] (meet2) {{$f$}};
  \draw (-1.15,-1) node[minimum height=0.5cm,minimum width=1cm,draw,fill=white] {{$\varphi$}};
%
\draw (0,-2.5) node[above] (X) {{$m$}};
\draw[color=blue!50!black] (-1.15,-2.5) node[above] (A1) {{$y$}};
\draw[color=green!50!black] (-0.17,-0.1) node[above] (A1) {{$\vdots$}};
\draw[color=green!50!black] (-0.22,-0.6) node[above] (A1) {{$x$}};
\draw[color=green!50!black] (-0.6,-0.21) node[above] (A1) {{$x$}};
\draw[color=green!50!black] (-1,0.18) node[above] (A1) {{$x$}};
\end{tikzpicture} 
\, =0.
  \end{equation}
  We denote by $\mod(Q_{x}, \varphi)$ be the full subcategory of $\mod(Q_{x})$ on the objects that satisfy the relation $\varphi$. 
\end{definition}
For the left hand side of Equation \eqref{eq:varphi-rel} we introduce for $(m,f) \in \Cat{C}_{z}$ the notation $\varphi(f):= (\varphi \act \id_{m}) \circ f^{\ast n} \in \Hom_{\Cat{A}}(y \act m, m)$, such that
$f \in \mod(Q_{x}, \varphi) \iff \varphi(f)=0$.
In particular, all modules $(m,0)$ satisfy the $\varphi$-relation and thus $\mod(Q_{x}, \varphi)$ contains $\Cat{M}$
as a full subcategory.

We assume for the remainder of this subsection, that $\Cat{A}$ is semisimple and $\Cat{M}$ is finite semisimple with labeling set $I$. Thus, the endofunctor $Q_{x}$ corresponds to a combinatorial quiver with vertices $I$ as in 
Definition \ref{definition:quiv-of-fun}.

Recall from Section \ref{sec:finite-quivers}, that a relation $\Omega_{\omega} \subset \PQ$ is an ideal generated by paths $(\omega_{\lambda})_{\lambda} \in \PQ$. It is called \emph{homogenous}, if all $\omega_{\lambda} \in \PQ^{n}$ for some $n \in \N$.
We want  to show that the  category $\mod(Q_{x}, \varphi)$ is the category of modules over the  quotient of the path algebra $\PQ_{x}$ by an ideal of homogenous relations. As in Equation \ref{eq:paths-gen-Q} we identify the space of paths
$\PQ(i,j)^{n} = \Hom_{\Cat{M}}(m_{j}, x^{\otimes n} \act m_{i})$.
\begin{definition}
  Let $Q=Q_{x}= x \act -: \Cat{M} \rightarrow \Cat{M}$ be a quiver  as above with labeling set $I$ for the representatives of the simples $m_{i} \in \Cat{M}$, $i \in I$ and let 
    $\varphi: y \rightarrow x^{\otimes n}$  be a morphism in $\Cat{A}$. 
  For all morphisms $\lambda: m_{j} \rightarrow y \act m_{i}$ in $\Cat{M}$ for $i,j \in I$, 
  define the path of length $n$
  \begin{equation}
    \label{eq:38}
  \varphi_{\lambda}=(\varphi \act m_{i}) \circ \lambda \in  \PQ(i,j)^{n}.
\end{equation}
The \emph{ideal $\Omega_{\varphi}$ of $\varphi$} is the ideal of $\PQ$ generated by all paths $\varphi_{\lambda_{i,j,\alpha}}$ for all $i,j \in I$ and all $\alpha$, where  $\{\lambda_{i,j, \alpha}\}_{\alpha}\subset
\Hom_{\Cat{M}}(m_{j}, x \act m_{i})$  form a basis for each $i,j \in I$. 
\end{definition}
  
Using semisimplicity of $\Cat{M}$, we get explicit basis coefficients

\begin{equation}
  \varphi_{\lambda_{i,j,\alpha}}=\;
  \begin{tikzpicture}[very thick,scale=1,color=blue!50!black, baseline]
\draw (0,-1.5) -- (0,1.5); 
\draw (0,-0.8) -- (-1.3,0.5); 
\draw[color=green!50!black] (-1.3,0.5)-- (-1.3,1.4); 
\draw[color=green!50!black] (-1.1,0.5)-- (-1.1,1.4); 
\draw[color=green!50!black] (-0.6,0.5)-- (-0.6,1.4); 
\draw (-1,0.4) node[minimum height=0.5cm,minimum width=1cm,draw,fill=white] {{$\varphi$}};
%
\draw [decorate,
    decoration = {brace}] (-1.4,1.5) --  (-0.5,1.5);
\draw (0,-2) node[above] (X) {{$j$}};
\draw[color=green!50!black] (-0.95,1.5) node[above] (A1) {{$n$}};
\draw[color=blue!50!black] (0,-0.8) node[right] (A1) {{$\lambda_{i,j,\alpha}$}};
\draw[color=blue!50!black] (0.2,0.5) node[above] (A1) {{$i$}};
\end{tikzpicture}
=
\sum_{i_{l},\alpha_{l}} \varphi_{i,i_{1}, \ldots, i_{n},j}^{\alpha;\alpha_{1},\ldots,\alpha_{n}}\;
 \begin{tikzpicture}[very thick,scale=1,color=blue!50!black, baseline]
\draw (0,-1.5) -- (0,2); 
\draw[color=green!50!black] (0,-1) -- (-0.8,0);  
\draw[color=green!50!black] (0,-0.2) -- (-0.8,0.8);  
\draw[color=green!50!black] (0,1.3) -- (-0.6,2);  

%
\draw (0,-2) node[above] (X) {{$j$}};
\draw[color=blue!50!black] (0,-1) node[right] (A1) {{$\alpha_{1}$}};
\draw[color=blue!50!black] (0,-0.6) node[right] (A1) {{$i_{1}$}};
\draw[color=blue!50!black] (0,-0.2) node[right] (A1) {{$\alpha_{2}$}};
\draw[color=blue!50!black] (0,0.2) node[right] (A1) {{$i_{2}$}};
\draw[color=green!50!black] (0,0.8) node[left] (A1) {{$\vdots$}};
\draw[color=blue!50!black] (0,1.3) node[right] (A1) {{$\alpha_{n}$}};
\draw[color=blue!50!black] (0,1.7) node[right] (A1) {{$i$}};
\draw[color=green!50!black] (-0.4,-1.1) node[above] (A1) {{$x$}};
\draw[color=green!50!black] (-0.4,-0.3) node[above] (A1) {{$x$}};
\draw[color=green!50!black] (-0.4,1.25) node[above] (A1) {{$x$}};
\end{tikzpicture}
\end{equation}
  with $\varphi_{i,j}^{\alpha;\alpha_{1},\ldots,\alpha_{n}}\in \Bbbk$ and abbreviating $\lambda_{i_{1},j, \alpha_{1}}$ by $\alpha_{1}$ on the right hand side and similar for the other basis vectors. 

  Using again semisimplicity,  the ideal $\Omega_{\varphi}$ of $\PQ$  consists explicity for all $i,j$ of all paths
  \begin{equation}
    \label{eq:6}
    \Omega_{\varphi}(i,j)=
  \begin{tikzpicture}[very thick,scale=1,color=blue!50!black, baseline]
\draw (0,-2) -- (0,2); 
\draw (-1.3,-0.8) -- (-1.3,0.5); 
\draw[color=green!50!black] (-2.5,-0.8)-- (-2.5,1.8); 
\draw[color=green!50!black] (-2.3,-0.8)-- (-2.3,1.8); 
\draw[color=green!50!black] (-1.9,-0.8)-- (-1.9,1.8); 
\draw[color=green!50!black] (-1.65,0.5)-- (-1.65,1.8); 
\draw[color=green!50!black] (-1.45,0.5)-- (-1.45,1.8); 
\draw[color=green!50!black] (-1.1,0.5)-- (-1.1,1.8); 
\draw[color=green!50!black] (-0.7,-0.8)-- (-0.7,1.8); 
\draw[color=green!50!black] (-0.5,-0.8)-- (-0.5,1.8); 
\draw[color=green!50!black] (-0.2,-0.8)-- (-0.2,1.8); 
\draw (-1.3,0.4) node[minimum height=0.5cm,minimum width=0.8cm,draw,fill=white] {{$\varphi$}};
\draw (-1.2,-1) node[minimum height=0.5cm,minimum width=3cm,draw,fill=white] {{}};
%
\draw [decorate,  decoration = {brace}] (-2.55,1.9) --  (-1.85,1.9);
\draw[color=green!50!black] (-2.2,2) node[above] (A1) {{$k$}};
\draw [decorate,  decoration = {brace}] (-1.7,1.9) --  (-1.05,1.9);
\draw[color=green!50!black] (-1.375,2) node[above] (A1) {{$n$}};
\draw [decorate,  decoration = {brace}] (-0.75,1.9) --  (-0.15,1.9);
\draw[color=green!50!black] (-0.45,2) node[above] (A1) {{$l$}};
\draw (0.2,-1.9) node[above] (X) {{$j$}};
\draw[color=blue!50!black] (0.2,0.5) node[above] (A1) {{$i$}};
\draw[color=green!50!black] (-2.7,-0.5) node[above] (A1) {{$x$}};
\draw[color=green!50!black] (-0.9,-0.5) node[above] (A1) {{$x$}};
\draw[color=blue!50!black] (-1.5,-0.5) node[above] (A1) {{$y$}};
\end{tikzpicture}
\end{equation}
  for all $k,l \in \N$.

  \begin{proposition}
    A module $(m,f)$ over $Q_{x}$ satisfies the $\varphi$-relation, if and only if on  all paths
    $\varphi_{\lambda_{i,j,\alpha}}\in \PQ^{n}$ for all $i,j, \alpha$, it satisfies  $f(\varphi_{\lambda_{i,j,\alpha}})=0$. This is the case if and only if,
    $f$ vanishes on all paths in $\Omega_{\varphi}$. 

    The category of modules over $\PQ/\Omega_{\varphi}$ is  canonically identified with
    \begin{equation}
      \label{eq:modules-relations}
      \mod(\PQ/\Omega_{\varphi})=\mod(Q_{x},\varphi). 
    \end{equation}
    \end{proposition}
\begin{proof}
  As in general for a relation $\Omega_{\varphi}$ on $Q$, a quiver module $f$ is in $\mod(Q,\varphi)$ if and only if
  $f(\omega)=0$ for all $\omega \in \Omega_{\varphi}$. This is by definition of $\Omega_{\varphi}$ the case if and only if
  $f(\varphi_{\lambda_{i,j,\alpha}})=0$ for all $i,j, \alpha$. By definition of $(\varphi_{\lambda_{i,j,\alpha}})$, this is the case if and only if
\begin{equation}
  \label{eq:30p}  
\begin{tikzpicture}[very thick,scale=1,color=blue!50!black, baseline]
\draw (0,-2.5) -- (0,1.5); 
\draw[color=green!50!black] (-0.7,-0.75) .. controls +(0,0.5) and +(-0.5,-0.5) .. (0,0.15); 
\draw[color=green!50!black] (-1.3,-0.75) .. controls +(0,0.5) and +(-0.5,-0.5) .. (0,0.8);
\draw[color=green!50!black] (-1.6,-0.75) .. controls +(0,0.5) and +(-0.5,-0.5) .. (0,1.1);
\draw[color=blue!50!black] (-1.15,-1) .. controls +(0,-0.5) and +(-0.5,0.5) .. (0,-2); 
\fill[color=blue!50!black] (0,0.15) circle (2.9pt) node[right] (meet2) {{$f$}};
\fill[color=blue!50!black] (0,0.8) circle (2.9pt) node[right] (meet2) {{$f$}};
\fill[color=blue!50!black] (0,1.1) circle (2.9pt) node[right] (meet2) {{$f$}};
  \draw (-1.15,-1) node[minimum height=0.5cm,minimum width=1cm,draw,fill=white] {{$\varphi$}};
%
  \draw (0,-2.35) node[right] (X) {{$j$}};
  \draw (0,-1.6) node[right] (X) {{$i$}};
\draw[color=blue!50!black] (-0.7,-2.1) node[above] (A1) {{$y$}};
\fill[color=blue!50!black] (0,-2) circle (2.9pt) node[right] (meet2) {{$\lambda_{\alpha}$}};
\fill[color=blue!50!black] (0,-1) circle (2.9pt) node[right] (meet2) {{$v$}};
\draw[color=green!50!black] (-0.17,-0.1) node[above] (A1) {{$\vdots$}};
\draw[color=green!50!black] (-0.22,-0.6) node[above] (A1) {{$x$}};
\draw[color=green!50!black] (-0.6,-0.21) node[above] (A1) {{$x$}};
\draw[color=green!50!black] (-1,0.18) node[above] (A1) {{$x$}};
\end{tikzpicture} 
\, =0 
\end{equation}

for all $i,j,\alpha$, which by semisimplicity is equivalent to $\varphi(f)=0$, which proves the claims. 
\end{proof}

\begin{remark}
  \begin{remarklist}
  \item  Not all homogenous relations on $Q_{x}$ can in general be obtained by this construction: by construction for each vertex $i$ there exists a vertex $j$ and a path from $i$ to $j$ in $\Omega_{\varphi}$.
  \item The construction can be generalized to the case that $Q: \Cat{M} \rightarrow \Cat{M}$ is an endofunctor,
    where we can take relations generated by natural transformations $\varphi: P \rightarrow Q^{n}$ for some endofunctor $P $ of $\Cat{M}$. 
  \end{remarklist}
\end{remark}
 \subsection{$q$-Relations and rigid monoidal subcategories of $\mod(Q_{z})$}
We now consider the case where $\Cat{A}$ is monoidal, $z \in \cent(\Cat{A})$ and $Q=Q_{z}$ is the associated Drinfeld quiver. We consider homogenous relations
as in Definition \ref{definition:mod-varphi} for $\mod(Q)$ that are compatible with the fusion product of Theorem \ref{thm:ten-quiv}. 

We therefor recall the usual  combinatorics of q-binomials: In the polynomial ring $\Z[q]$, one defines for $n \in \N$
the $q$-numbers $(n)_{q}= 1 + q + \ldots + q^{n-1}$,  the $q$-factorials
$(n)!_{q}=(n)_{q}\cdot \ldots \cdot (2)_{q}(1)_{q}$ and finally for $i \leq n$ the $q$-binomial
$\binom{n}{i}_{\!q}= \frac{(n)!_{q}}{(n-i)!_{q}(i)!_{q}}$, which happens to lie in $\Z[q]$.  
The importance  of these $q$-numbers is due to the following well-known
\begin{lemma}
  \label{lemma:q-comm-Alg}
  Let $A$ be an associative  algebra over $\Z[q]$ and $x,y \in A$ with $xy=qyx$. Then
  \begin{equation}
    \label{eq:q-comm}
    (x+y)^{n}=  \sum_{i=0}^{n}\binom{n}{i}_{\!\!\! q}y^{i}x^{n-i}.  
  \end{equation}
\end{lemma}
Suppose now that $q \in \C$ is a primitive $N$-th root of unity. Then $(N)_{q}=0$ and also $\binom{N}{k}_{\!q}=0$ for all
$0<k<N$, thus in this  case in the setting of the Lemma \ref{lemma:q-comm-Alg}, $(x+y)^{N}=x^{N}+y^{N}$.  We use these combinatorics and define
\begin{definition}
  \label{definition:q-state}
  Let $\Cat{A}$ be a monoidal category and $z \in \cent(\Cat{A})$ with generalized Drinfeld quiver $Q_{z}$. For $N \in \N$ and $q \in \C$ a primitive $N$-th root of unity, a \emph{$q$-relation of length $N$ for $Q_{z}$} is an object $y \in \cent(\Cat{A})$ together with  a morphism $\varphi: y \rightarrow z^{\otimes N}$ in $\cent(\Cat{A})$ such that
  \begin{equation}
    \label{eq:q-state}
    (\id_{z^{\otimes i}} \otimes  c_{z,z} \otimes \id_{z^{\otimes (N-i-2)}}) \circ \varphi= q \cdot \varphi.
  \end{equation}
  for all $i\in {0, \ldots, N-2}$. 
\end{definition}

It follows, that for $\varphi_{1}: y_{1} \rightarrow z_{1}^{\otimes N}$ and   $\varphi_{2}: y_{2} \rightarrow z_{2}^{\otimes N}$ two $q$-relations of the same length and the same root of unity $q$, also their sum
\begin{equation}
  \label{eq:37}
  \varphi_{1} \oplus \varphi_{2}: y_{1} \oplus y_{2} \rightarrow (z_{1} \oplus z_{2})^{\otimes N}
\end{equation}
is a $q$-relation of length $N$ with the same $q$ for the sum of the Drinfeld quivers $Q_{z_{1}\oplus z_{2}}$.

\begin{theorem}
  \label{theorem:main-rel}
Let $\Cat{A}$ be a monoidal category with generalized Drinfeld quiver $Q_{z}$ for 
   $z \in \cent(\Cat{A})$ and $\varphi$ a $q$-relation of length $N$ for $z$. 
  \begin{theoremlist}
  \item If $f,g \in \mod(Q_{z})$ satisfy the relation $\varphi$, also $f \qoti g$ satisfies the relation $\varphi$.
  \item If $f \in  \mod(Q_{z})$  satisfy the relation $\varphi$, also $f^{\dagger}$ and $\tensor[^{\dagger}]{f}{}$   satisfy the relation $\varphi$.
        \end{theoremlist}
   The categories $\mod(Q_{z},\varphi)$ and $\mod(Q_{z}, \varphi)^{\mathsf{ln}}$ of locally nilpotent modules in $\mod(Q_{z},\varphi)$ are   rigid  monoidal abelian subcategory of $\mod(Q_{z})$. 
 \end{theorem}
 \begin{proof}
   Assume  $(a,f),(b,g) \in \Cat{C}_{z}$ satisfy the $\varphi$-relation. The morphism $\varphi(f \qoti g)$ is by definition of $\qoti$ a sum of $2^{N}$-terms:
   \begin{equation}
     \label{eq:varphi-tech}
     \begin{split}
       &\varphi(f \qoti g)= \\
       &\sum_{s \in \{0,1\}^{N}} ((f^{s_{N}} \otimes g^{1-s_{N}})(\id \otimes c_{z,a} \otimes \id )^{1-s_{N}})  \ldots ((f^{s_{1}} \otimes g^{1-s_{1}})(\id \otimes c_{z,a} \otimes \id )^{1-s_{1}})  
     (\varphi \otimes \id_{a} \otimes \id_{b}),   
     \end{split}
       \end{equation}
   where for a morphism $h:v \rightarrow w$ in $\Cat{A}$ we use the convention, that $h^{1}=h$ and $h^{0}=\id_{v}$.
   Graphically for $N=3$, the term with $s=(s_{i})_{i=1,2,3}=(1,0,0)$ is represented as
     \begin{equation}
    \label{eq:Pic29}
\begin{tikzpicture}[very thick,scale=1,color=blue!50!black, baseline]
  \draw (0,-2) -- (0,1.5); 
  \draw[color=green!50!black] (-0.7,-0.75) .. controls +(0,0.5) and +(-0.5,-0.5) .. (0,0.15); 
  \draw[color=white, line width=4pt] (-1.3,-0.75) .. controls +(0,1) and +(-0.5,-0.5) .. (0.6,0.8); 
  \draw[color=green!50!black] (-1.3,-0.75) .. controls +(0,1) and +(-0.5,-0.5) .. (0.6,0.8);
  \draw[color=white, line width=4pt] (-1.6,-0.75) .. controls +(0,1) and +(-0.5,-0.5) .. (0.6,1.1); 
\draw[color=green!50!black] (-1.6,-0.75) .. controls +(0,1) and +(-0.5,-0.5) .. (0.6,1.1);
\draw (-1.15,-1) -- (-1.15,-2); 
  \draw (0.6,-2) -- (0.6,1.5); 
\fill[color=blue!50!black] (0,0.2) circle (2.9pt) node[right] (meet2) {{$f$}};
\fill[color=blue!50!black] (0.6,0.8) circle (2.9pt) node[right] (meet2) {{$g$}};
\fill[color=blue!50!black] (0.6,1.1) circle (2.9pt) node[right] (meet2) {{$g$}};
  \draw (-1.15,-1) node[minimum height=0.5cm,minimum width=1cm,draw,fill=white] {{$\varphi$}};
%
\draw (0,-2.5) node[above] (X) {{$a$}};
\draw[color=blue!50!black] (-1.15,-2.5) node[above] (A1) {{$y$}};
\draw[color=green!50!black] (-0.9,-0.75) node[above] (A1) {{$\dots$}};
\draw[color=green!50!black] (-0.22,-0.6) node[above] (A1) {{$z$}};
\draw[color=green!50!black] (-0.6,-0.25) node[above] (A1) {{$z$}};
\draw[color=green!50!black] (-1,0.18) node[above] (A1) {{$z$}};
\end{tikzpicture} 
=0 \, .
  \end{equation}
  To abbreviate, we denote the term corresponding to $s \in \{0,1\}^{N}$ by $<\varphi, g^{1-s_{N}} \ast f^{s_{N}}\ast \ldots \ast g^{1-s_{1}} \ast f^{s_{1}}>$, so that for instance the term in Figure \eqref{eq:Pic29} is denoted
  $<\varphi, g \ast g \ast f>$. With the braiding of $\cent(\Cat{A})$ we can permute $f$ and $g$ as follows. Assume that for a fixed $s$ there is $s_{i}=1$ and $s_{i+1}=0$. Then the corresponding term in
  $\varphi(f \qoti g)$ is, using \eqref{eq:q-state}:
 \begin{equation}
  \label{eq:4}  
\begin{tikzpicture}[very thick,scale=1,color=blue!50!black, baseline]
  \draw (0,-2) -- (0,1.5); 
  \draw[color=green!50!black] (-0.9,-0.75) .. controls +(0,0.5) and +(-0.5,-0.5) .. (0,0.2); 
  \draw[color=white, line width=4pt] (-1.3,-0.75) .. controls +(0,1) and +(-0.5,-0.5) .. (0.6,1); 
  \draw[color=green!50!black] (-1.3,-0.75) .. controls +(0,1) and +(-0.5,-0.5) .. (0.6,1);
\draw (-1.15,-1) -- (-1.15,-2); 
  \draw (0.6,-2) -- (0.6,1.5); 
\fill[color=blue!50!black] (0,0.2) circle (2.9pt) node[right] (meet2) {{$f$}};
\fill[color=blue!50!black] (0.6,1) circle (2.9pt) node[right] (meet2) {{$g$}};
  \draw (-1.15,-1) node[minimum height=0.5cm,minimum width=1cm,draw,fill=white] {{$\varphi$}};
%
  \draw (0,-2.5) node[above] (X) {{$a$}};
  \draw (0.6,-2.5) node[above] (X) {{$b$}};
\draw[color=blue!50!black] (-1.15,-2.5) node[above] (A1) {{$y$}};
\draw[color=green!50!black] (-1.5,-0.75) node[above] (A1) {{$\dots$}};
\draw[color=green!50!black] (-0.4,-0.75) node[above] (A1) {{$\dots$}};
\draw[color=green!50!black] (-0.22,-0.5) node[above] (A1) {{$z$}};
\draw[color=green!50!black] (-0.6,0.23) node[above] (A1) {{$z$}};
\end{tikzpicture} 
\, =
\begin{tikzpicture}[very thick,scale=1,color=blue!50!black, baseline]
  \draw (0,-2) -- (0,1.5); 
  \draw[color=green!50!black] (-0.9,-0.75) .. controls +(0,0.5) and +(-0.5,-0.5) .. (0,1); 
  \draw[color=white, line width=4pt] (-1.3,-0.75) .. controls +(0,1) and +(-0.5,-0.5) .. (0.6,1); 
  \draw[color=green!50!black] (-1.3,-0.75) .. controls +(0,1) and +(-0.5,-0.5) .. (0.6,1);
\draw (-1.15,-1) -- (-1.15,-2); 
  \draw (0.6,-2) -- (0.6,1.5); 
\fill[color=blue!50!black] (0,1) circle (2.9pt) node[right] (meet2) {{$f$}};
\fill[color=blue!50!black] (0.6,1) circle (2.9pt) node[right] (meet2) {{$g$}};
  \draw (-1.15,-1) node[minimum height=0.5cm,minimum width=1cm,draw,fill=white] {{$\varphi$}};
%
  \draw (0,-2.5) node[above] (X) {{$a$}};
    \draw (0.6,-2.5) node[above] (X) {{$b$}};
\draw[color=blue!50!black] (-1.15,-2.5) node[above] (A1) {{$y$}};
\draw[color=green!50!black] (-1.5,-0.75) node[above] (A1) {{$\dots$}};
\draw[color=green!50!black] (-0.4,-0.75) node[above] (A1) {{$\dots$}};
\draw[color=green!50!black] (-0.64,-0.5) node[above] (A1) {{$z$}};
\draw[color=green!50!black] (-0.6,0.23) node[above] (A1) {{$z$}};
\end{tikzpicture} 
= q \cdot
\begin{tikzpicture}[very thick,scale=1,color=blue!50!black, baseline]
  \draw (0,-2) -- (0,1.5); 
  \draw[color=white, line width=4pt](-0.9,-0.75) .. controls +(0,0.5) and +(-0.5,-0.5) .. (0.6,1); 
  \draw[color=green!50!black] (-0.9,-0.75) .. controls +(0,0.5) and +(-0.5,-0.5) .. (0.6,1); 
  \draw[color=green!50!black] (-1.3,-0.75) .. controls +(0,1) and +(-0.5,-0.5) .. (0,1);
\draw (-1.15,-1) -- (-1.15,-2); 
  \draw (0.6,-2) -- (0.6,1.5); 
\fill[color=blue!50!black] (0.6,1) circle (2.9pt) node[right] (meet2) {{$g$}};
\fill[color=blue!50!black] (0,1) circle (2.9pt) node[right] (meet2) {{$f$}};
  \draw (-1.15,-1) node[minimum height=0.5cm,minimum width=1cm,draw,fill=white] {{$\varphi$}};
%
  \draw (0,-2.5) node[above] (X) {{$a$}};
    \draw (0.6,-2.5) node[above] (X) {{$b$}};
\draw[color=blue!50!black] (-1.15,-2.5) node[above] (A1) {{$y$}};
\draw[color=green!50!black] (-1.5,-0.75) node[above] (A1) {{$\dots$}};
\draw[color=green!50!black] (-0.4,-0.75) node[above] (A1) {{$\dots$}};
\draw[color=green!50!black] (-0.3,-0.4) node[above] (A1) {{$z$}};
\draw[color=green!50!black] (-0.7,0.4) node[above] (A1) {{$z$}};
\end{tikzpicture} 
\end{equation}
  We obtain, as in Lemma \ref{lemma:q-comm-Alg},
  \begin{equation}
    \label{eq:22}
    \begin{split}
      \varphi(f \qoti g) &= \sum_{s \in \{0,1\}^{N}} <\varphi, g^{1-s_{N}} \ast f^{s_{N}}\ast \ldots \ast g^{1-s_{1}} \ast f^{s_{1}}> \\
      &= \sum_{i=0}^{N} \binom{N}{i}_{\!\!\!q} <\varphi, f^{\ast i} \ast g^{\ast (N-i)}> \\
      &=<\varphi, \id_{b} \ast f^{\ast N} > + <\varphi, g^{\ast N} \ast \id_{a}>,
    \end{split}
  \end{equation}
  using in the last step that $q$ is a primitive $N$-th root of unity. 
  The remaining terms are $<\varphi, \id_{b} \ast f^{\ast N} >=\varphi(f) \otimes \id_{b}=0$ and
  \begin{equation}
  \label{eq:4}  
  <\varphi, g^{\ast N} \ast \id_{a}>=
\begin{tikzpicture}[very thick,scale=1,color=blue!50!black, baseline]
\draw (0,-2) -- (0,1.5); 
\draw (-1.15,-1) -- (-1.15,-2); 
\draw (-0.5,-2) -- (-0.5,1.5); 
\draw[color=white, line width=4pt]  (-0.7,-0.75) .. controls +(0,0.5) and +(-0.5,-0.5) .. (0,0.15); 
\draw[color=green!50!black] (-0.7,-0.75) .. controls +(0,0.5) and +(-0.5,-0.5) .. (0,0.15); 
\draw[color=white, line width=4pt] (-1.3,-0.75) .. controls +(0,0.5) and +(-0.5,-0.5) .. (0,0.8);
\draw[color=green!50!black] (-1.3,-0.75) .. controls +(0,0.5) and +(-0.5,-0.5) .. (0,0.8);
\draw[color=white, line width=4pt]  (-1.6,-0.75) .. controls +(0,0.5) and +(-0.5,-0.5) .. (0,1.1); 
\draw[color=green!50!black] (-1.6,-0.75) .. controls +(0,0.5) and +(-0.5,-0.5) .. (0,1.1);
%
%
\fill[color=blue!50!black] (0,0.15) circle (2.9pt) node[right] (meet2) {{$g$}};
\fill[color=blue!50!black] (0,0.8) circle (2.9pt) node[right] (meet2) {{$g$}};
\fill[color=blue!50!black] (0,1.1) circle (2.9pt) node[right] (meet2) {{$g$}};
  \draw (-1.15,-1) node[minimum height=0.5cm,minimum width=1cm,draw,fill=white] {{$\varphi$}};
%
  \draw (0,-2.5) node[above] (X) {{$b$}};
  \draw (-0.5,-2.5) node[above] (X) {{$a$}};
\draw[color=blue!50!black] (-1.15,-2.5) node[above] (A1) {{$y$}};
\draw[color=green!50!black] (-0.2,-0.1) node[above] (A1) {{$\vdots$}};
\draw[color=green!50!black] (-0.22,-0.5) node[above] (A1) {{$z$}};
\draw[color=green!50!black] (-0.63,-0.23) node[above] (A1) {{$z$}};
\draw[color=green!50!black] (-1,0.18) node[above] (A1) {{$z$}};
\end{tikzpicture} 
=
\begin{tikzpicture}[very thick,scale=1,color=blue!50!black, baseline]
\draw (0,-2) -- (0,1.5); 
\draw (-0.5,-2)  .. controls +(0,0.2) and +(0.5,0) .. (-1.3,-1.5) 
 .. controls +(-1.5,0) and +(0,-0.5) .. (-1,1.5);
 \draw[color=white, line width=4pt] (-1.15,-1) -- (-1.15,-2); 
 \draw (-1.15,-1) -- (-1.15,-2); 
\draw[color=white, line width=4pt]  (-0.7,-0.75) .. controls +(0,0.5) and +(-0.5,-0.5) .. (0,0.15); 
\draw[color=green!50!black] (-0.7,-0.75) .. controls +(0,0.5) and +(-0.5,-0.5) .. (0,0.15); 
\draw[color=white, line width=4pt] (-1.3,-0.75) .. controls +(0,0.5) and +(-0.5,-0.5) .. (0,0.8);
\draw[color=green!50!black] (-1.3,-0.75) .. controls +(0,0.5) and +(-0.5,-0.5) .. (0,0.8);
\draw[color=white, line width=4pt]  (-1.6,-0.75) .. controls +(0,0.5) and +(-0.5,-0.5) .. (0,1.1); 
\draw[color=green!50!black] (-1.6,-0.75) .. controls +(0,0.5) and +(-0.5,-0.5) .. (0,1.1);
%
%
\fill[color=blue!50!black] (0,0.15) circle (2.9pt) node[right] (meet2) {{$g$}};
\fill[color=blue!50!black] (0,0.8) circle (2.9pt) node[right] (meet2) {{$g$}};
\fill[color=blue!50!black] (0,1.1) circle (2.9pt) node[right] (meet2) {{$g$}};
  \draw (-1.15,-1) node[minimum height=0.5cm,minimum width=1cm,draw,fill=white] {{$\varphi$}};
%
  \draw (0,-2.5) node[above] (X) {{$b$}};
  \draw (-0.5,-2.5) node[above] (X) {{$a$}};
\draw[color=blue!50!black] (-1.15,-2.5) node[above] (A1) {{$y$}};
\draw[color=green!50!black] (-0.2,-0.1) node[above] (A1) {{$\vdots$}};
\draw[color=green!50!black] (-0.22,-0.5) node[above] (A1) {{$z$}};
\draw[color=green!50!black] (-0.63,-0.23) node[above] (A1) {{$z$}};
\draw[color=green!50!black] (-1,0.18) node[above] (A1) {{$z$}};
\end{tikzpicture} 
=0 \; ,
\end{equation}
  using, that $\varphi$ is a morphism in $\cent(\Cat{A})$ in the second step.

  Suppose that $(a,f)$ satisfies  relation $\varphi$. To see that also $f^{\dagger}$ satisfies relation $\varphi$, we depict the argument for $N=2$ in the graphical calculus, the general case is given by a repeated use of this identity:
\begin{equation}
  \label{eq:comp-dagger-varphi}
\begin{tikzpicture}[very thick,scale=1,color=blue!50!black, baseline=-1.5cm]
  \draw (0,-0.5) -- (0,0.5);
  \draw[color=blue!50!black] (0,-0.5) .. controls +(0,-0.5) and +(0,-0.5) .. (0.5,-0.5); 
  \draw[color=blue!50!black] (0,0.5) .. controls +(0,0.5) and +(0,0.5) .. (-0.5,0.5);
  \draw (0.5,-0.5) -- (0.5,1);
  \draw (-0.5,0.5) -- (-0.5,-1);
   \draw (0-1,-0.5-2) -- (0-1,0.5-2);
  \draw[color=blue!50!black] (0-1,-0.5-2) .. controls +(0,-0.5) and +(0,-0.5) .. (0.5-1,-0.5-2); 
  \draw[color=blue!50!black] (0-1,0.5-2) .. controls +(0,0.5) and +(0,0.5) .. (-0.5-1,0.5-2);
  \draw (0.5-1,-0.5-2) -- (0.5-1,1-2);
  \draw (-0.5-1,0.5-2) -- (-0.5-1,-2-2);
  %
  \draw (-2.5,-3) -- (-2.5,-4);
  %
  \draw[color=white, line width=4pt]  (-2.8,-2.8) .. controls +(0,1) and +(-0.5,-0.5) .. (0,0);  
  \draw[color=green!50!black] (-2.8,-2.8) .. controls +(0,1) and +(-0.5,-0.5) .. (0,0); 
  \draw[color=white, line width=4pt]  (-2.2,-2.8) .. controls +(0,0.5) and +(-0.5,-0.5) .. (-1,-2);  
  \draw[color=green!50!black] (-2.2,-2.8) .. controls +(0,0.5) and +(-0.5,-0.5) .. (-1,-2);
  \draw (-2.5,-3) node[minimum height=0.5cm,minimum width=1cm,draw,fill=white] {{$\varphi$}};
  %
   \fill[color=blue!50!black] (0,0) circle (2.9pt) node[right] (meet2) {{$f$}};
   \fill[color=blue!50!black] (0-1,0-2) circle (2.9pt) node[right] (meet2) {{$f$}};
\draw (-0.5,-1.5) node[right] (X) {{$a$}};
\draw[color=green!50!black] (-2,-2.5) node[above] (A1) {{$z$}};
\draw[color=green!50!black] (-2.5,-2) node[above] (A1) {{$z$}};
\draw[color=blue!50!black] (-2.3,-4) node[above] (A1) {{$y$}};
\end{tikzpicture} 
\, =
\begin{tikzpicture}[very thick,scale=1,color=blue!50!black, baseline=-1.5cm]
    \draw (0-1,-0.5-2) -- (0-1,0.5-1);
  \draw[color=blue!50!black] (0-1,-0.5-2) .. controls +(0,-0.5) and +(0,-0.5) .. (0.5-1,-0.5-2); 
  \draw[color=blue!50!black] (0-1,0.5-1) .. controls +(0,0.5) and +(0,0.5) .. (-0.5-1,0.5-1);
  \draw (0.5-1,-0.5-2) -- (0.5-1,1);
  \draw (-0.5-1,0.5-1) -- (-0.5-1,-2-2);
  %
  \draw (-2.5,-3) -- (-2.5,-4);
%
 \draw[color=white, line width=4pt]  (-2.2,-2.8) .. controls +(0,1) and +(-0.5,-0.5) .. (0-1,-1); 
  \draw[color=green!50!black] (-2.2,-2.8) .. controls +(0,1) and +(-0.5,-0.5) .. (0-1,-1); 
  \draw[color=white, line width=4pt]   (-2.8,-2.8) .. controls +(0,0.5) and +(-0.5,-0.5) .. (-1,-2); 
  \draw[color=green!50!black] (-2.8,-2.8) .. controls +(0,0.5) and +(-0.5,-0.5) .. (-1,-2);
  \draw (-2.5,-3) node[minimum height=0.5cm,minimum width=1cm,draw,fill=white] {{$\varphi$}};
  %
   \fill[color=blue!50!black] (0-1,-1) circle (2.9pt) node[right] (meet2) {{$f$}};
   \fill[color=blue!50!black] (0-1,0-2) circle (2.9pt) node[right] (meet2) {{$f$}};
   \draw (-1.1,-1.5) node[right] (X) {{$a$}};
\draw[color=green!50!black] (-1.8,-1.7) node[above] (A1) {{$z$}};
\draw[color=green!50!black] (-2.5,-2.5) node[above] (A1) {{$z$}};
\draw[color=blue!50!black] (-2.3,-4) node[above] (A1) {{$y$}};;
\end{tikzpicture}
= q \cdot
\begin{tikzpicture}[very thick,scale=1,color=blue!50!black, baseline=-1.5cm]
  
   \draw (0-1,-0.5-2) -- (0-1,0.5-1);
  \draw[color=blue!50!black] (0-1,-0.5-2) .. controls +(0,-0.5) and +(0,-0.5) .. (0.5-1,-0.5-2); 
  \draw[color=blue!50!black] (0-1,0.5-1) .. controls +(0,0.5) and +(0,0.5) .. (-0.5-1,0.5-1);
  \draw (0.5-1,-0.5-2) -- (0.5-1,1);
  \draw (-0.5-1,0.5-1) -- (-0.5-1,-2-2);
  %
  \draw (-2.5,-3) -- (-2.5,-4);
  %
  \draw[color=white, line width=4pt]   (-2.2,-2.8) .. controls +(0,0.5) and +(-0.5,-0.5) .. (-1,-2); 
  \draw[color=green!50!black] (-2.2,-2.8) .. controls +(0,0.5) and +(-0.5,-0.5) .. (-1,-2);
 \draw[color=white, line width=4pt]  (-2.8,-2.8) .. controls +(0,1) and +(-0.5,-0.5) .. (0-1,-1); 
  \draw[color=green!50!black] (-2.8,-2.8) .. controls +(0,1) and +(-0.5,-0.5) .. (0-1,-1); 
  \draw (-2.5,-3) node[minimum height=0.5cm,minimum width=1cm,draw,fill=white] {{$\varphi$}};
  %
   \fill[color=blue!50!black] (0-1,-1) circle (2.9pt) node[right] (meet2) {{$f$}};
   \fill[color=blue!50!black] (0-1,0-2) circle (2.9pt) node[right] (meet2) {{$f$}};
   \draw (-1.1,-1.5) node[right] (X) {{$a$}};
\draw[color=green!50!black] (-2,-2.5) node[above] (A1) {{$z$}};
\draw[color=green!50!black] (-2.5,-2) node[above] (A1) {{$z$}};
\draw[color=blue!50!black] (-2.3,-4) node[above] (A1) {{$y$}};;
\end{tikzpicture}
=0
\end{equation}
Analogously, we conclude, that $\tensor[^{\dagger}]{f}{}$ satisfies relation $\varphi$.

Clearly the subcategory $\Cat{A}_{Q}$ of vertex semisimple modules $(a, 0)$ of $\mod(Q_{z})$ satisfies relation $\varphi$, in particular the monoidal unit $(\unit, 0)$. It follows, that $\mod(Q_{z}, \varphi)$ is a rigid monoidal subcategory of $\mod(Q_{z})$. The category $\mod(Q_{z}, \varphi)^{\mathsf{ln}}$ of additionally locally nilpotent modules is as
 the intersection of two rigid monoidal subcategories also a rigid monoidal subcategory. 
\end{proof}
Since the category $\mod(Q_{z}, \varphi)$ is equivalent to the category of finite dimensional $\PQ_{z}/\Omega_{\varphi}$-modules, we deduce from  Proposition \ref{proposition:self-inj}:
\begin{corollary}
  In the situation of Theorem \ref{theorem:main-rel}, the
  algebra  $\PQ_{z}/\Omega_{\varphi}$ is self-injective: The
  classes of finite dimensional projective and injective $\PQ_{z}/\Omega_{\varphi}$-modules agree. 
\end{corollary}
We proceed to consider first examples of $q$-relations. 
  \begin{example}
    \begin{examplelist}
      \item A length one $q$-relation for $Q_{z}$  with $z \in \cent(\Cat{A})$ has $q=1$ and is just a morphism $\varphi:y \rightarrow z$ in $\cent(\Cat{A})$. It is trivially a $q$-relation for $q=1$.
Since $\cent(\Cat{A})$ is semisimple, we can consider the image $\im(\varphi) \subset z$ and choose a splitting $z \cong \im (\varphi) \oplus \widetilde{z}$
in $\cent(\Cat{A})$. It follows easily that we can identify those modules over $Q_{z}$ that satisfy the relation $\varphi$
with the modules over $Q_{\widetilde{z}}$, i.e. there is an equivalence of categories
$\mod(Q_{z},\varphi) \cong \mod(Q_{\widetilde{z}})$.  In this sense length $1$ relations are not interesting.
\item For a length two $q$-relation and
  $\chara(\Bbbk) \neq 2$,
  we have $q=-1$ and it follows, that a corresponding $q$-relation $\varphi: y \rightarrow z^{\otimes 2}$ fulfills   $c_{z,z} \circ c_{z,z} \circ \varphi=\varphi$. Conversely for a
  morphism  $\psi: y \rightarrow z \otimes z$ in $\cent(\Cat{A})$ such that $c_{z,z} \circ c_{z,z} \circ \psi=\psi$,
  the morphism
  \begin{equation}
    \label{eq:25}
    \varphi= \frac{1}{2} (\psi- c_{z,z} \circ \psi)
  \end{equation}
  is a length $2$ $q$-relation.
    \end{examplelist}
 \end{example}
 The following particular length two $q$-relations for  $\Cat{A}$  a pivotal multi-fusion category will be related to preprojective algebras in Section
  \ref{sec:prepr-algebra} to preprojective algebras. 

  Recall the notion of twist and Frobenius-Schur indicators for $\Cat{A}$, see Section \ref{sec:mono-categ-their}.
\begin{lemma}
  \label{lemma:prepro-relation}
  Let $z \in \cent(\Cat{A})$ be a  selfdual simple object in a pivotal multi-fusion category $\Cat{A}$   with twist $\theta_{z} \in \Bbbk^{\times}$ 
  such that
  \begin{equation}
    \label{eq:41}
    \theta_{z} \cdot \FS(z)=-1. 
  \end{equation}
  Then $(\id \otimes \phi)\coev{z}: \unit \rightarrow z \otimes z$  is a length $2$ $q$-relation for any isomorphism
  $\phi: z^{*} \rightarrow z$. 
\end{lemma}
\begin{proof}
  This follows directly from
  \begin{equation}
  \label{eq:45}
  \begin{tikzpicture}[very thick,scale=1,color=blue!50!black, baseline]
    \draw[color=blue!50!black] (-0.5,-0.5) .. controls +(0,-0.5) and +(0,-0.5) .. (-1.2,-0.5);
    \draw[color=blue!50!black] (-0.5,-0.5)  -- (-0.5,0);
    \draw[color=blue!50!black] (-1.2,-0.5)  -- (-1.2,0);
    %
    \draw[color=blue!50!black] (-0.5,0) .. controls +(0,0.3) and +(0,-0.5) .. (-1.2,1.1);
     \draw[color=white, line width=4pt]  (-1.2,0) .. controls +(0,0.2) and +(0,-0.5) .. (-0.5,1.1); 
    \draw[color=blue!50!black] (-1.2,0) .. controls +(0,0.2) and +(0,-0.5) .. (-0.5,1.1);
    %
    %
    \fill[color=blue!50!black] (-0.5,-0.2) circle (2.9pt) node[right] (meet2) {{$\phi$}};
    \draw[color=blue!50!black] (-0.24,-0.8) node[above] (A1) {{$x^{*}$}};
        \draw[color=blue!50!black] (-0.35,0) node[above] (A1) {{$x$}};
\draw[color=blue!50!black] (-1.4,-0.8) node[above] (A1) {{$x$}};
  \end{tikzpicture}
  = 
  \begin{tikzpicture}[very thick,scale=1,color=blue!50!black, baseline]
     \draw[color=blue!50!black] (-0.5,-0.5) .. controls +(0,-0.5) and +(0,-0.5) .. (-1.2,-0.5);
    \draw[color=blue!50!black] (-0.5,-0.5)  -- (-0.5,0);
    \draw[color=blue!50!black] (-1.2,-0.5)  -- (-1.2,0);
    %
    \draw[color=blue!50!black] (-0.5,0) .. controls +(0,0.3) and +(0,-0.5) .. (-1.2,1.1);
     \draw[color=white, line width=4pt]  (-1.2,0) .. controls +(0,0.2) and +(0,-0.5) .. (-0.5,1.1); 
    \draw[color=blue!50!black] (-1.2,0) .. controls +(0,0.2) and +(0,-0.5) .. (-0.5,1.1);
    %
    \fill[color=blue!50!black] (-1.1,0.75) circle (2.9pt) node[left] (meet2) {{$\phi$}};
     \draw[color=blue!50!black] (-0.24,-0.8) node[above] (A1) {{$x^{*}$}};
        \draw[color=blue!50!black] (-1.4,0.8) node[above] (A1) {{$x$}};
             \draw[color=blue!50!black] (-1.4,-0.8) node[above] (A1) {{$x$}};
      \end{tikzpicture}
  \stackrel{~\eqref{eq:twist-def}}{=} \theta_{x}
  \begin{tikzpicture}[very thick,scale=1,color=blue!50!black, baseline]
     \draw[color=blue!50!black] (-0.5,-0.5) .. controls +(0,-0.5) and +(0,-0.5) .. (-1.2,-0.5);
    \draw[color=blue!50!black] (-0.5,-0.5)  -- (-0.5,1.1);
    \draw[color=blue!50!black] (-1.2,-0.5)  -- (-1.2,1.1);
    %
    %
    \fill[color=blue!50!black] (-1.2,0.3) circle (2.9pt) node[left] (meet2) {{$\phi$}};
    \draw[color=blue!50!black] (-0.24,-0.65) node[above] (A1) {{$x$}};
    \draw[color=blue!50!black] (-1.4,-0.65) node[above] (A1) {{$x^{*}$}};
    \draw[color=blue!50!black] (-1.4,0.6) node[above] (A1) {{$x$}};
  \end{tikzpicture}
  \stackrel{~\eqref{eq:FS-def}}{=}\theta_{x} \FS(x)
  \begin{tikzpicture}[very thick,scale=1,color=blue!50!black, baseline]
    \draw[color=blue!50!black] (-0.5,-0.5) .. controls +(0,-0.5) and +(0,-0.5) .. (-1.2,-0.5);
    \draw[color=blue!50!black] (-0.5,-0.5)  -- (-0.5,1.1);
    \draw[color=blue!50!black] (-1.2,-0.5)  -- (-1.2,1.1);
    %
    %
    %
    \fill[color=blue!50!black] (-0.5,0.1) circle (2.9pt) node[right] (meet2) {{$\phi$}};
    \draw[color=blue!50!black] (-0.24,-0.6) node[above] (A1) {{$x^{*}$}};
        \draw[color=blue!50!black] (-0.24,0.25) node[above] (A1) {{$x$}};
\draw[color=blue!50!black] (-1.4,-0.6) node[above] (A1) {{$x$}};
  \end{tikzpicture}
  \, .
\end{equation}
\end{proof}

We conclude the section with an example of a class of $q$-relations, where the resulting monoidal category is
the category of modules over a finite-dimensional Hopf algebra.

Consider the fusion quiver $Q_{z}$ from Example \ref{example:first-expl} \refitem{item:cyclic-grp}: For  $G=\Z/ n\Z$ with generator $g \in
G$, the object  $z \in \cent(\vect_{G})$ is the simple object 
 $z= e_{g}$ with half-braiding determined by a $n$-th root of unity $q$. We assume now that $q$ is primitive. 

By construction, the identity in $\vect_{G}$ provides a morphism  $\unit \cong e_{g}^{\otimes n}$. It is easy to see that this is a morphism in $\cent(\Cat{A})$ and  defines  a   non-zero $q$-relation of length $n$ 
  \begin{equation}
    \label{eq:expl-varphi}
    \varphi: \unit \rightarrow z^{\otimes n}.
  \end{equation}

The corresponding rigid monoidal category is identified with the modules over the Taft algebra:  

  \begin{proposition}
    In the situation above, 
    there is an equivalence of  monoidal categories
    \begin{equation}
      \label{eq:34}
      \mod(Q_{z},\varphi) \cong \mod(H),
    \end{equation}
    where $\mod(H)$ is the category of finite-dimensional modules over the Taft Hopf algebra  $H$ corresponding to $G$ and $q$. 
  \end{proposition}
  \begin{proof} 
    In \cite[Thm. 2.1]{CibQuiv} it  is shown that the algebra $\PQ_{z}/\varphi$ has the structure  of a Hopf algebra, which is isomorphic to the dual of the Taft Hopf algebra. Since the Taft Hopf algebra $H$ is selfdual, we obtain an isomorphism $\PQ_{z}/\varphi \cong H$ of Hopf algebra.
   
    If we denote the vertices of $\PQ$ by $e_{i}$, $i=0,1, \ldots, n-1$ and the arrow from vertex $i$ to $i+1$ ( mod $n$)  by $\lambda_{i}$,  the coproduct is  given by
    \begin{equation}
      \label{eq:def-coprod-cib}
      \Delta(\lambda_{i})= \sum_{j+k=i}\lambda_{k} \otimes e_{j} + \sum_{j+k=i} q^{j} e_{k} \otimes \lambda_{j},
    \end{equation}
    as in the proof of  \cite[Thm. 2.1]{CibQuiv} (after the change of coordinates given by multiplying the arrow $\lambda_{i}$ with $q^{i}$).

    It remains to show, that the monoidal product on $ \mod(Q_{z},\varphi)$ given by $\Delta$ agrees with the  monoidal product $\qoti$.  Let $(V_{s},f_{s})$ for $s=1,2$ with $f_{s}: z \otimes V_{s} \rightarrow V_{s}$ be two objects in $ \mod(Q_{z},\varphi)$.   It suffices to compare the  action of one arrow $\lambda_{i}$ on the respective tensor products. For  $v_{j} \otik w_{k} \in V_{j} \otik W_{k}$ with $j+k=i$, the  action of $\lambda_{i}$ is according to Equation \ref{eq:def-coprod-cib}  given by
    \begin{equation}
      \label{eq:action-Delta-Taft}
      \Delta(\lambda_{i})(v \otik w)= \lambda_{i}(v) \otik w + q^{j} w \otik \lambda_{i}(v).
    \end{equation}
    For the action of $\lambda_{i}$ on $(V_{1},f_{1}) \qoti (V_{2},f_{2})$ we compute \ref{eq:f-on-arrow}, regarding
    $v \otimes w$ as morphism  $e_{g^{i}} \rightarrow V \otimes W$ in  $\vect_{G}$:
    \begin{equation}
      \label{eq:f-on-arrow-Taft}
      (f_{1} \qoti f_{2})(\lambda_{i})(v \otimes w)=
\begin{tikzpicture}[very thick,scale=1.3,color=blue!50!black, baseline]
\draw (0,-1.5) -- (0,1.5); 
\draw[color=green!50!black] (-1,0) .. controls +(0,0.5) and +(-0.5,-0.5) .. (0,1);
\draw[color=green!50!black] (-1,0) .. controls +(0,-0.5) and +(-0.5,0.5) .. (0,-1);
%
%
\fill[color=blue!50!black] (0,1) circle (2.9pt) node[right] (meet1) {{$f_{1} \qoti f_{2} $}};
\fill[color=blue!50!black] (0,0.2) circle (2.9pt) node[right] (v) {{$v \otimes w$}};
\fill[color=blue!50!black] (0,-1) circle (2.9pt) node[right] (meet2) {{$\lambda_{i}$}};
\draw (0,-2.1) node[above] (X) {{$g^{i+1}$}};
\draw[color=green!50!black] (-0.7,-0.56) node[left] (A1) {{$z$}};
\draw[color=blue!50!black] (-0.1,-0.4) node[right] (A1) {{$g^{i}$}};
\draw[color=blue!50!black] (-0.1,0.6) node[right] (A1) {{$ V_{1} \otimes V_{2}$}};
\draw[color=blue!50!black] (-0.1,1.35) node[right] (A1) {{$V_{1} \otimes V_{2}$}};
\end{tikzpicture} 
\, . 
    \end{equation}
Using the definition  of $\qoti$ from Equation \ref{eq:def-qoti} and the  braiding $c_{z,e_{g^{j}}}=q^{j} \id$ from Example \ref{example:first-expl}\ref{item:cyclic-grp}, it follows directly that this action of $\lambda_{i}$ agrees with the action in Equation \ref{eq:action-Delta-Taft}. 
  \end{proof}

\section{Preprojective algebra}
\label{sec:prepr-algebra}
 
We apply the notion of $q$-relations for the fusion quivers $Q_{z}$ to construct rigid monoidal structures
on the modules of certain preprojective algebras. In this section we assume for simplicity
$\chara(\Bbbk)=0$. 

\subsection{Variants of preprojective algebras}

We recall the classical preprojective algebra \cite{GP} of a combinatorial quiver and consider a twisted deformation thereof. Then we  define the category of $\Pi$-modules for a generalized quiver, generalizing the category of modules over the preprojective algebras.

Let $Q$ be a combinatorial quiver with vertex set $I$  and arrow  $\alpha: i \rightarrow j$ in $Q$.  Consider the quiver $\cc{Q}$ with reversed arrows, i.e. $\cc{Q}$ has an arrow   $\cc{\alpha}: j \rightarrow i$ for all  $\alpha: i \rightarrow j$ in $Q$. Finally   the double quiver of $Q$ is the sum   $DQ=Q \oplus \cc{Q}$. 
  Define for all $i \in I$ the following paths in $DQ$
 \begin{equation}
    \label{eq:classical-preproj}
    \omega_{i}= \sum_{\alpha, s(\alpha)=i} \cc{\alpha}\alpha-\sum_{\beta, t(\beta)=i}
    \beta \cc{\beta},
  \end{equation}
  and 
   \begin{equation}
    \label{eq:twisted-prepro}
    \omega_{i}^{\theta}= \sum_{\alpha, s(\alpha)=i} \cc{\alpha}\alpha-\theta\cdot \sum_{\beta, t(\beta)=i}
    \beta \cc{\beta},
  \end{equation}
  \begin{definition}
\label{definition:cl-prepro}
    Let $Q$ be a quiver with basis $\{\alpha: i \rightarrow j\}$ for its paths as above.
    \begin{definitionlist}
    \item   \label{item:cl-prepro} The \emph{classical preprojective algebra $\Pi(Q)$} of  $Q$
is $\Pi(Q)= \PDQ /\Omega$ with 
    the ideal  $\Omega$ generated by the paths $\{\omega_{i}\}_{i \in I}$.  
 \item  \label{item:theta-prepro}For $\theta \in \Bbbk^{\times}$, the \emph{$\theta$-twisted preprojective algebra} is the algebra
  \begin{equation}
    \label{eq:twisted-prepro}
    \Pi^{\theta}(Q)= \PDQ/\Omega^{\theta},
  \end{equation}
  where the ideal $\Omega^{\theta}$ is  generated by the  paths  $\{\omega_{i}^{\theta}\}_{i \in I}$. 
    \end{definitionlist}
      \end{definition}

We proceed to define the corresponding analogue for generalized quivers. Let   $\Cat{M}$ be a linear category and  $Q: \Cat{M} \rightarrow \Cat{M}$ a generalized quiver on $\Cat{M}$. Recall that a \emph{biadjoint} of $Q$  is a functor $Q^{*}: \Cat{M} \rightarrow \Cat{M}$ with has the structure of a left and a right adjoint to $Q$.
We fix the corresponding units and counits $(\eta^{r}, \epsilon^{r})$ and
$(\eta^{l}, \epsilon^{l})$ of the right and left adjunction, respectively.
We consider the endofunctor $D(Q)=Q \oplus Q^{*}$ as quiver and define the following relations on the corresponding modules. 
Note that a  module  $(m,\widetilde{f}: D(Q)(m) \rightarrow m)$ over $D(Q)$  consists of the  sum $\widetilde{f}=f + \cc{f}$ of two  modules $f:Q(m)\rightarrow m$ and $\cc{f}: Q^{*}(m) \rightarrow m$.
\begin{definition}
  \label{definition:double-quiv-prepro}
  Let $Q: \Cat{M} \rightarrow \Cat{M}$ be a generalized quiver with biadjoint $Q^{*}: \Cat{M} \rightarrow \Cat{M}$.
  \begin{definitionlist}
  \item The \emph{double quiver of $Q$} is the generalized quiver $D(Q)= Q \oplus Q^{*}: \Cat{M} \rightarrow \Cat{M}$.
    \item   
      The \emph{category of $\Pi$-modules over $Q$} consists of the full subcategory of those $DQ$-modules $(m,f + \cc{f}: DQ(m) \rightarrow m)$,   with $f:Q(m)\rightarrow m$ and $\cc{f}: Q^{*}(m) \rightarrow m$ such that
      \begin{equation}
        \label{eq:gen-preproj}
        \cc{f} \circ Q^{*}(f) \circ \eta^{r}_{m}- f \circ Q(\cc{f}) \circ \eta^{l}_{m}=0.
      \end{equation}
       \end{definitionlist}
  \end{definition}
  In graphical terms, the relation \eqref{eq:gen-preproj} is
\begin{equation}
\begin{tikzpicture}[very thick,scale=1,color=blue!50!black, baseline]
  \draw (0.7,-1) -- (0.7,1.3); 
  %
  \draw[color=green!50!black] (0.7,0.7) --  (0.7-0.5,0.7-0.5); 
    \draw[color=green!50!black] (0.7,0.15) --  (0.7-0.5,0.15-0.5); 
   \draw[color=green!50!black]  (0.7-0.5,0.7-0.5)  .. controls +(-0.5,-0.5) and +(-0.5,-0.5) .. (0.7-0.5,0.15-0.5);
%
\fill[color=blue!50!black] (0.7,0.15) circle (2.9pt) node[right] (meet2) {{$f$}};
\fill[color=blue!50!black] (0.7,0.7) circle (2.9pt) node[right] (meet2) {{$\cc{f}$}};
%
 \draw (0.65,-0.9) node[right] (X) {{$m$}};
 \draw[color=green!50!black] (0.4,-0.8) node[above] (A1) {{$Q$}};
 \draw[color=green!50!black] (0.2,0.1) node[above] (A1) {{$Q^{*}$}};
\end{tikzpicture} 
-
\begin{tikzpicture}[very thick,scale=1,color=blue!50!black, baseline]
  \draw (0.7,-1) -- (0.7,1.3); 
  %
  \draw[color=green!50!black] (0.7,0.7) --  (0.7-0.5,0.7-0.5); 
    \draw[color=green!50!black] (0.7,0.15) --  (0.7-0.5,0.15-0.5); 
   \draw[color=green!50!black]  (0.7-0.5,0.7-0.5)  .. controls +(-0.5,-0.5) and +(-0.5,-0.5) .. (0.7-0.5,0.15-0.5);
%
\fill[color=blue!50!black] (0.7,0.15) circle (2.9pt) node[right] (meet2) {{$\cc{f}$}};
\fill[color=blue!50!black] (0.7,0.7) circle (2.9pt) node[right] (meet2) {{$f$}};
%
 \draw (0.65,-0.9) node[right] (X) {{$m$}};
 \draw[color=green!50!black] (0.4,-0.85) node[above] (A1) {{$Q^{*}$}};
 \draw[color=green!50!black] (0.2,0.1) node[above] (A1) {{$Q$}};
\end{tikzpicture} 
=0 \, .
\end{equation}

   Note that the category of $\Pi$-modules depends on the choice of units and counits of the biadjunction, as can be seen in the following 
  \begin{example}
    \label{example:preproj-classic}
    Let $\Cat{M}$ be finite semisimple with labeling set $I$ for the simples $(x_{i})_{i \in I}$ and  let $Q \in \End(\Cat{M})$ be a linear functor. We can then choose a biadjunction that reproduces the classical preprojective algebra of $Q$:
    The functor $Q$ is up to isomorphism determined by the matrix of vector spaces $(\Hom_{\Cat{M}}(x_{i}, Q(x_{j})))_{i,j \in I}$, indeed we can assume that
    \begin{equation}
      \label{eq:42}
      Q(x_{j}) = \oplus_{i \in I} \Hom_{\Cat{M}}(x_{i}, Q(x_{j})) \otik x_{i}. 
    \end{equation}
    Define a linear endofunctor $Q^{*}$ of $\Cat{M}$ by
    \begin{equation}
      \label{eq:dual-quiver}
      Q^{*}(x_{j})= \oplus_{i \in I} \Hom_{\Cat{M}}(x_{j},Q(x_{i}))^{*} \otik x_{i}.
    \end{equation}
    In case that $\Cat{M}$ is rigid monoidal and $Q=Q_{x}$ is given by $Q_{x}=x \otimes -$ for $x \in \Cat{M}$, it follows that
    $Q^{*}_{x}\cong Q_{x^{*}}$ as quivers, where $x^{*}$ is the right dual of $x$.
       The functor $Q^{*}$ has the structure of a biadjoint of $Q$:  Since
    \begin{equation}
      \label{eq:49}
     Q^{*} (Q(x_{j}))= \oplus_{i,k \in I} \Hom_{\Cat{M}}(x_{i}, Q(x_{j})) \otik \Hom_{\Cat{M}}(x_{i}, Q(x_{k}))^{*} \otik x_{k}, 
   \end{equation}
   we obtain from the usual coevaluation for finite dimensional vector spaces
linear maps  $\eta_{x_{j}}^{r}: x_{j} \rightarrow Q^{*}Q(x_{j})$ for all $j \in I$ and thus a linear natural transformation $\eta^{r}: \id \rightarrow Q^{*} Q$ by Lemma \ref{lemm:ses-cat-fun}. 
Analogously we obtain from the evaluation maps in $\vect$ a natural transformation $\epsilon^{r}: QQ^{*} \rightarrow \id$, such that the triangle identities hold.
The adjunction $(\epsilon^{l},\eta^{l})$ is defined analogously.
This defines  the biadjoint $Q^{*}$ and upon choosing bases  $\lambda_{i,j, \alpha} \in \Hom_{\Cat{M}}(x_{i}, Q(x_{j}))$ with
    dual bases $\lambda_{i,j}^{\alpha} \in \Hom_{\Cat{M}}(x_{i}, Q(x_{j}))^{*}$,  the functor $D(Q)$ corresponds to the  double of the  combinatorial quiver $Q$ and 
a morphism  $f+\cc{f}: DQ(m) \rightarrow m$  is a $\Pi$-module over $Q$ if and only if for all $j \in I$,
    \begin{equation}
      \label{eq:preproj-rel}
      \sum_{i,\alpha}\cc{f}(\lambda_{i,j}^{\alpha}) \circ f(\lambda_{i,j,\alpha})-f(\lambda_{i,j,\alpha})\circ \cc{f}(\lambda_{i,j}^{\alpha})=0,
    \end{equation}
This is precisely the  relation for the classical preprojective algebra of $Q$ as in  Definition  \ref{definition:cl-prepro} \refitem{item:cl-prepro}.

However, if we for choose some $\theta \in \Bbbk^{\times}$ the unit $\theta \cdot \eta^{l}$ for the left adjoint, with left evaluation $\frac{1}{\theta}\epsilon^{l}$, while keeping
$(\epsilon^{r},\eta^{r})$, we still get a biadjunction, generalizing the $\theta$-twisted preprojective algebra.
The resulting categories differ in general:
In case that $Q$ is the Jordan quiver with one loop, the
    usual $\Pi$-modules consists of a pair of commuting  endomorphism $f, \cc{f}: V \rightarrow V$ of a vector space $V$, while a $\theta \neq 1$ gives the  category of endomorphism $f,\cc{f}$ with $\cc{f}f=\theta f \cc{f}$ (for instance if $V=\Bbbk$, there is no such pair of endomorphisms).  
  \end{example}

In the sequel we consider two classes of preprojective algebras, where the theory of Drinfeld quivers equips the $\Pi$-modules with a rigid monoidal structure.

\subsection{Preprojective algebra on a self-dual quiver}

In the first case let $\Cat{A}$ be a braided monoidal category and $z \in \cent(\Cat{A})$ be a self dual object. 
In this case the quiver $Q_{z}$ is already a double quiver on a suitable subquiver of $Q_{z}$ with half the arrows.
 Under an additional assumptions on  $z$ the preprojective relation is a  $q$-relation for $Q_{z}$.

 \begin{theorem}
   \label{theorem:prepro-selfdual}
  Let $\Cat{A}$ be a  pivotal braided fusion category and $z \in \Cat{A}$ selfdual, with $\FS(z)=-1$ and  $\theta_{z}=1$.
  Let $\varphi=(\id \otimes \phi) \circ \coev{z}: \unit \rightarrow z \otimes z$ for any isomorphism $\varphi:z^{*} \rightarrow z$
  with corresponding ideal $\Omega_{\varphi} \subset \PQ$ for $Q=Q_{z}$.  
  Then  the algebra
  \begin{equation}
    \label{eq:33}
   \Pi_{\varphi}=\PQ/\Omega_{\varphi} 
  \end{equation}
    is a classical preprojective algebra on a subquiver of $Q_{z}$, and the  category of  modules $\mod(\Pi_{\varphi})$ has a rigid monoidal structure. 
 \end{theorem}
\begin{proof}
  By Lemma \ref{lemma:prepro-relation}, $\varphi$ is a $q$-relation with $q=-1$ and thus it follows from Theorem \ref{theorem:main-rel}
  that $\mod(\Pi_{\varphi})$ has a  rigid monoidal structure induced from the monoidal structure on the modules over the Drinfeld quiver $Q_{z}$.
  It remains to show that $\Pi_{\varphi}$ is a preprojective algebra. We use a symplectic pairing to write the quiver $Q$ as a double quiver:
  The space of arrows of $Q$ is  $Q^{1}= \oplus_{i,j \in I}\Hom_{\Cat{A}}(x_{j}, z \otimes x_{i})$, using  the isomorphism in Equation \eqref{eq:arrow-space-iso}. 
It  has a linear symplectic structure
  $<-,->: Q^{1} \otik Q^{1} \rightarrow \Bbbk$, given by
 \begin{equation}
  \label{eq:4}  
\left<
\begin{tikzpicture}[very thick,scale=1,color=blue!50!black, baseline]
\draw (0,-1) -- (0,1); 
\draw[color=green!50!black] (-1,1) .. controls +(0,-0.5) and +(-0.5,0.5) .. (0,0);
\fill[color=blue!50!black] (0,0) circle (2.9pt) node[right] (meet2) {{$\lambda$}};
\draw[color=blue!50!black] (-0.1,-0.5) node[right] (A1) {{$j$}};
  \draw[color=blue!50!black] (-0.1,0.7) node[right] (A1) {{$i$}};
    \draw[color=green!50!black] (-0.8,0.6) node[right] (A1) {{$z$}};
\end{tikzpicture} 
\, ,
\begin{tikzpicture}[very thick,scale=1,color=blue!50!black, baseline]
\draw (0,-1) -- (0,1); 
\draw[color=green!50!black] (-1,1) .. controls +(0,-0.5) and +(-0.5,0.5) .. (0,0);
\fill[color=blue!50!black] (0,0) circle (2.9pt) node[right] (meet2) {{$\mu$}};
\draw[color=blue!50!black] (-0.1,-0.5) node[right] (A1) {{$i$}};
  \draw[color=blue!50!black] (-0.1,0.7) node[right] (A1) {{$j$}};
    \draw[color=green!50!black] (-0.8,0.6) node[right] (A1) {{$z$}};
\end{tikzpicture} 
\right> =
\frac{1}{d_{i}d_{j}}
\tr \left(
\begin{tikzpicture}[very thick,scale=1,color=blue!50!black, baseline,
dot/.style = {circle, fill, minimum size=#1,
              inner sep=0pt, outer sep=0pt},
dot/.default = 5.8pt  
                      ]
\draw (0,-1.5) -- (0,1.5);  
\draw[color=green!50!black] (-0.7,1.5) .. controls +(0,-0.2) and +(-0.5,0.5) .. (0,0.5)  node[pos=0.5,dot,label=below:{$\phi$}] {};   
\draw[color=green!50!black] (-1,-1.5) .. controls +(0,0.5) and +(-0.5,0.7) .. (0,-0.5);    
\fill[color=blue!50!black] (0,-0.5) circle (2.9pt) node[right] (meet2) {{$\mu$}};
\fill[color=blue!50!black] (0,0.5) circle (2.9pt) node[right] (meet2) {{$\lambda$}};
\draw (0,2) node[below] (X) {{$i$}};
\draw (0,-2) node[above] (X) {{$i$}};
\draw[color=green!50!black] (-0.7,2) node[below] (A1) {{$z^{*}$}};
\draw[color=green!50!black] (-1,-2) node[above] (A1) {{$z^{*}$}};
\draw[color=blue!50!black] (0,0) node[right] (A1) {{$j$}};
\end{tikzpicture} 
\right)=
  \frac{1}{d_{i}d_{j}} \cdot \;
  \begin{tikzpicture}[very thick,scale=1,color=blue!50!black, baseline,
    dot/.style = {circle, fill, minimum size=#1,
              inner sep=0pt, outer sep=0pt},
dot/.default = 5.8pt  
                     ]
\draw (0,-1) -- (0,1); 
\draw[color=green!50!black] (-1.4,0.7) .. controls +(0,0.5) and +(-0.5,0.5) .. (0,0.7)
node[pos=0.7,dot,label=below:{$\phi$}] {}; ; 
\draw[color=green!50!black] (-1.4,0.7) .. controls +(0,-0.5) and +(-0.5,0.5) .. (0,-0.7);
\draw[color=blue!50!black] (0,1) .. controls +(0,1) and +(0,1) .. (-2,1); 
\draw[color=blue!50!black] (0,-1) .. controls +(0,-1) and +(0,-1) .. (-2,-1); 
\draw[color=blue!50!black] (-2,1) -- (-2,-1); 
%
\fill[color=blue!50!black] (0,0.7) circle (2.9pt) node[right] (meet1) {{$\lambda$}};
\fill[color=blue!50!black] (0,-0.7) circle (2.9pt) node[right] (meet2) {{$\mu$}};
\draw (-0.1,-1.2) node[right] (X) {{$i$}};
\draw (-0.1,-0) node[right] (X) {{$j$}};
\draw (-0.1,1.2) node[right] (X) {{$i$}};
\draw[color=green!50!black] (-0.7,-0.1) node[below] (A1) {{$z$}};
\end{tikzpicture} 
\end{equation}
  for a choice of isomorphism $\phi: z \rightarrow z^{*}$ and with dimensions $\dim(x_{i})=d_{i}\in \Bbbk^{\times}$.
  Indeed, the non-degeneracy follows from the non-degeneracy of the trace and by the cyclicity of the left trace, we obtain (recall, that we suppress the pivotal structure of $\Cat{A}$)
 
\begin{align*} 
  \label{eq:10}
  <\lambda,\mu>=&
   \frac{1}{d_{i}d_{j}} \cdot \;
  \begin{tikzpicture}[very thick,scale=1,color=blue!50!black, baseline,
    dot/.style = {circle, fill, minimum size=#1,
              inner sep=0pt, outer sep=0pt},
dot/.default = 5.8pt  
                     ]
\draw (0,-1) -- (0,1); 
\draw[color=green!50!black] (-1.4,0.7) .. controls +(0,0.5) and +(-0.5,0.5) .. (0,0.7)
node[pos=0.7,dot,label=below:{$\phi$}] {}; ; 
\draw[color=green!50!black] (-1.4,0.7) .. controls +(0,-0.5) and +(-0.5,0.5) .. (0,-0.7);
\draw[color=blue!50!black] (0,1) .. controls +(0,1) and +(0,1) .. (-2,1); 
\draw[color=blue!50!black] (0,-1) .. controls +(0,-1) and +(0,-1) .. (-2,-1); 
\draw[color=blue!50!black] (-2,1) -- (-2,-1); 
%
\fill[color=blue!50!black] (0,0.7) circle (2.9pt) node[right] (meet1) {{$\lambda$}};
\fill[color=blue!50!black] (0,-0.7) circle (2.9pt) node[right] (meet2) {{$\mu$}};
\draw (-0.1,-1.2) node[right] (X) {{$i$}};
\draw (-0.1,-0) node[right] (X) {{$j$}};
\draw (-0.1,1.2) node[right] (X) {{$i$}};
\draw[color=green!50!black] (-0.7,-0.1) node[below] (A1) {{$z$}};
\end{tikzpicture} 
=
 \frac{1}{d_{i}d_{j}} \cdot \;
  \begin{tikzpicture}[very thick,scale=1,color=blue!50!black, baseline,
    dot/.style = {circle, fill, minimum size=#1,
              inner sep=0pt, outer sep=0pt},
dot/.default = 5.8pt  
                     ]
\draw (0,-1) -- (0,1); 
\draw[color=green!50!black] (-1.4,0.7) .. controls +(0,0.5) and +(-0.5,0.5) .. (0,0.7) ; 
\draw[color=green!50!black] (-1.4,0.7) .. controls +(0,-0.5) and +(-0.5,0.5) .. (0,-0.7)
node[pos=0.7,dot,label=below:{$\phi$}] {};
\draw[color=blue!50!black] (0,1) .. controls +(0,1) and +(0,1) .. (-2,1); 
\draw[color=blue!50!black] (0,-1) .. controls +(0,-1) and +(0,-1) .. (-2,-1); 
\draw[color=blue!50!black] (-2,1) -- (-2,-1); 
%
\fill[color=blue!50!black] (0,0.7) circle (2.9pt) node[right] (meet1) {{$\mu$}};
\fill[color=blue!50!black] (0,-0.7) circle (2.9pt) node[right] (meet2) {{$\lambda$}};
\draw (-0.1,-1.2) node[right] (X) {{$j$}};
\draw (-0.1,-0) node[right] (X) {{$i$}};
\draw (-0.1,1.2) node[right] (X) {{$j$}};
\draw[color=green!50!black] (-0.9,0) node[below] (A1) {{$z$}};
\end{tikzpicture}              
                  \\
  =&
     -
      \frac{1}{d_{i}d_{j}} \cdot \;
  \begin{tikzpicture}[very thick,scale=1,color=blue!50!black, baseline,
    dot/.style = {circle, fill, minimum size=#1,
              inner sep=0pt, outer sep=0pt},
dot/.default = 5.8pt  
                     ]
\draw (0,-1) -- (0,1); 
\draw[color=green!50!black] (-1.4,0.7) .. controls +(0,0.5) and +(-0.5,0.5) .. (0,0.7)
node[pos=0.7,dot,label=below:{$\phi$}] {}; 
\draw[color=green!50!black] (-1.4,0.7) .. controls +(0,-0.5) and +(-0.5,0.5) .. (0,-0.7);
\draw[color=blue!50!black] (0,1) .. controls +(0,1) and +(0,1) .. (-2,1); 
\draw[color=blue!50!black] (0,-1) .. controls +(0,-1) and +(0,-1) .. (-2,-1); 
\draw[color=blue!50!black] (-2,1) -- (-2,-1); 
%
\fill[color=blue!50!black] (0,0.7) circle (2.9pt) node[right] (meet1) {{$\mu$}};
\fill[color=blue!50!black] (0,-0.7) circle (2.9pt) node[right] (meet2) {{$\lambda$}};
\draw (-0.1,-1.2) node[right] (X) {{$j$}};
\draw (-0.1,-0) node[right] (X) {{$i$}};
\draw (-0.1,1.2) node[right] (X) {{$j$}};
\draw[color=green!50!black] (-0.9,0) node[below] (A1) {{$z$}};
\end{tikzpicture}
     =-<\mu,\lambda> \, ,
\end{align*}
  using $\FS(x)=-1$ in the third step.
  We thus find a Darboux-basis of $Q^{1}$ of the form $\{\lambda_{i,j,\alpha}, \lambda_{i,j}^{\beta}\}_{i,j \in I, \alpha,\beta}$ with $(\lambda_{i,j,\alpha})_{\alpha}\in \Hom(x_{j},z \otimes x_{i})$ and $(\lambda_{i,j}^{\beta})_{\beta} \in
  \Hom(x_{i}, z \otimes x_{j})$, such that $<\lambda_{i,j,\alpha}, \lambda_{i,j}^{\beta}>=\delta_{\alpha}^{\beta}$ for all indices. The arrows $\{\lambda_{i,j,\alpha}\}_{i,j,\alpha}$ span a subquiver $Q'$ of $Q$ such that $Q=DQ'$ is the double quiver of $Q'$. 
  If the indices $i,j$ of the basis vectors are clear from the context, they are omitted in the sequel. 

  It follows that for all indices, 

\begin{equation}
  \label{eq:11}
  \begin{tikzpicture}[very thick,scale=1,color=blue!50!black, baseline,
    dot/.style = {circle, fill, minimum size=#1,
              inner sep=0pt, outer sep=0pt},
dot/.default = 5.8pt  
                     ]
\draw (0,-1.5) -- (0,1.5); 
\draw[color=green!50!black] (-1.4,0.7) .. controls +(0,0.5) and +(-0.5,0.5) .. (0,0.5)
node[pos=0.7,dot,label=below:{$\phi$}] {}; 
\draw[color=green!50!black] (-1.4,0.7) .. controls +(0,-0.5) and +(-0.5,0.5) .. (0,-0.5);
%
%
\fill[color=blue!50!black] (0,0.5) circle (2.9pt) node[right] (meet1) {{$\lambda_{\alpha}$}};
\fill[color=blue!50!black] (0,-0.5) circle (2.9pt) node[right] (meet2) {{$\lambda^{\beta}$}};
\draw (-0.1,-1.2) node[right] (X) {{$k$}};
\draw (-0.1,-0) node[right] (X) {{$i$}};
\draw (-0.1,1.2) node[right] (X) {{$k$}};
\draw[color=green!50!black] (-0.9,0.1) node[below] (A1) {{$z$}};
\end{tikzpicture}
=\delta^{\alpha}_{\beta} \, d_{i} \cdot
\begin{tikzpicture}[very thick,scale=1,color=blue!50!black, baseline,
    dot/.style = {circle, fill, minimum size=#1,
              inner sep=0pt, outer sep=0pt},
dot/.default = 5.8pt  
                     ]
\draw (0,-1.5) -- (0,1.5); 
\draw (-0.1,0) node[right] (X) {{$k$}};
\end{tikzpicture}
\, , \quad \text{and} \quad
  \begin{tikzpicture}[very thick,scale=1,color=blue!50!black, baseline,
    dot/.style = {circle, fill, minimum size=#1,
              inner sep=0pt, outer sep=0pt},
dot/.default = 5.8pt  
                     ]
\draw (0,-1.5) -- (0,1.5); 
\draw[color=green!50!black] (-1.4,0.7) .. controls +(0,0.5) and +(-0.5,0.5) .. (0,0.5)
node[pos=0.7,dot,label=below:{$\phi$}] {}; 
\draw[color=green!50!black] (-1.4,0.7) .. controls +(0,-0.5) and +(-0.5,0.5) .. (0,-0.5);
%
%
\fill[color=blue!50!black] (0,0.5) circle (2.9pt) node[right] (meet1) {{$\lambda^{\beta}$}};
\fill[color=blue!50!black] (0,-0.5) circle (2.9pt) node[right] (meet2) {{$\lambda_{\alpha}$}};
\draw (-0.1,-1.2) node[right] (X) {{$k$}};
\draw (-0.1,-0) node[right] (X) {{$i$}};
\draw (-0.1,1.2) node[right] (X) {{$k$}};
\draw[color=green!50!black] (-0.9,0.1) node[below] (A1) {{$z$}};
\end{tikzpicture}
=-\delta^{\alpha}_{\beta} \, d_{i} \cdot
\begin{tikzpicture}[very thick,scale=1,color=blue!50!black, baseline,
    dot/.style = {circle, fill, minimum size=#1,
              inner sep=0pt, outer sep=0pt},
dot/.default = 5.8pt  
                     ]
\draw (0,-1.5) -- (0,1.5); 
\draw (-0.1,0) node[right] (X) {{$k$}};
\end{tikzpicture}
\, .
\end{equation} 
This implies the identity
\begin{equation}
  \label{eq:repl-straight}
  \begin{tikzpicture}[very thick,scale=1,color=blue!50!black, baseline,
    dot/.style = {circle, fill, minimum size=#1,
              inner sep=0pt, outer sep=0pt},
dot/.default = 5.8pt  
               ]
\draw[color=green!50!black] (0,-1.5) -- (0,1.5); 
\draw[color=green!50!black] (-0.1,0) node[right] (X) {{$z$}};
\draw (0.8,-1.5) -- (0.8,1.5); 
\draw (0.8-0.1,0) node[right] (X) {{$i$}};

\end{tikzpicture}
  = \frac{1}{d_{i}} \left( \sum_{k,\alpha}
\begin{tikzpicture}[very thick,scale=1,color=blue!50!black, baseline,
dot/.style = {circle, fill, minimum size=#1,
              inner sep=0pt, outer sep=0pt},
dot/.default = 5.8pt  
                      ]
\draw (0,-1.5) -- (0,1.5);  
\draw[color=green!50!black] (-0.7,1.5) .. controls +(0,-0.2) and +(-0.5,0.5) .. (0,0.5)  ;   
\draw[color=green!50!black] (-1,-1.5) .. controls +(0,0.5) and +(-1,1) .. (0,-0.5)
node[pos=0.9,dot,label=below:{$\phi$}] {};    
\fill[color=blue!50!black] (0,-0.5) circle (2.9pt) node[right] (meet2) {{$\lambda_{\alpha}$}};
\fill[color=blue!50!black] (0,0.5) circle (2.9pt) node[right] (meet2) {{$\lambda^{\alpha}$}};
\draw (0,2) node[below] (X) {{$i$}};
\draw (0,-2) node[above] (X) {{$i$}};
\draw[color=green!50!black] (-0.7,2) node[below] (A1) {{$z$}};
\draw[color=green!50!black] (-1,-2) node[above] (A1) {{$z$}};
\draw[color=blue!50!black] (0,0) node[right] (A1) {{$k$}};
\end{tikzpicture} 
-
\begin{tikzpicture}[very thick,scale=1,color=blue!50!black, baseline,
dot/.style = {circle, fill, minimum size=#1,
              inner sep=0pt, outer sep=0pt},
dot/.default = 5.8pt  
                      ]
\draw (0,-1.5) -- (0,1.5);  
\draw[color=green!50!black] (-0.7,1.5) .. controls +(0,-0.2) and +(-0.5,0.5) .. (0,0.5)  ;   
\draw[color=green!50!black] (-1,-1.5) .. controls +(0,0.5) and +(-1,1) .. (0,-0.5)
node[pos=0.9,dot,label=below:{$\phi$}] {};    
\fill[color=blue!50!black] (0,-0.5) circle (2.9pt) node[right] (meet2) {{$\lambda^{\alpha}$}};
\fill[color=blue!50!black] (0,0.5) circle (2.9pt) node[right] (meet2) {{$\lambda_{\alpha}$}};
\draw (0,2) node[below] (X) {{$i$}};
\draw (0,-2) node[above] (X) {{$i$}};
\draw[color=green!50!black] (-0.7,2) node[below] (A1) {{$z$}};
\draw[color=green!50!black] (-1,-2) node[above] (A1) {{$z$}};
\draw[color=blue!50!black] (0,0) node[right] (A1) {{$k$}};
\end{tikzpicture} 
\right) \; ,
\end{equation}
as can be seen by post-composing both sides  with an arbitrary basis vector.

We now compute the paths $\omega_{i}$ that generate the relation $\Omega_{\varphi}$ as
\begin{align*}
  \label{eq:13}
  \begin{tikzpicture}[very thick,scale=1,color=blue!50!black, baseline,
dot/.style = {circle, fill, minimum size=#1,
              inner sep=0pt, outer sep=0pt},
dot/.default = 5.8pt  
                      ]
    \draw[color=green!50!black] (-0.5,-0.5) .. controls +(0,-0.5) and +(0,-0.5) .. (-1.2,-0.5);
    \draw[color=green!50!black] (-0.5,-0.5)  -- (-0.5,1.5);
    \draw[color=green!50!black] (-1.2,-0.5)  -- (-1.2,1.5)
    node[pos=0.5,dot,label=left:{$\phi^{-1}$}] {}; 
        \draw[color=blue!50!black] (0.3,-1.5)  -- (0.3,1.5);
    %
    %
    \draw[color=green!50!black] (-0.24,-0.6) node[above] (A1) {{$z$}};
    \draw[color=green!50!black] (-1.4,-0.6) node[above] (A1) {{$z^{*}$}};
    \draw (0.3,0) node[left] (X) {{$i$}};
  \end{tikzpicture}
  \stackrel{~\eqref{eq:repl-straight}}{=}&
  \frac{1}{d_{i}} \left( \sum_{k, \alpha}
\begin{tikzpicture}[very thick,scale=1,color=blue!50!black, baseline,
dot/.style = {circle, fill, minimum size=#1,
              inner sep=0pt, outer sep=0pt},
dot/.default = 5.8pt  
                      ]
                      \draw (0,-1.5) -- (0,1.5);
\draw[color=green!50!black] (-0.7,1.5) .. controls +(0,-0.2) and +(-0.5,0.5) .. (0,0.5)  ;   
\draw[color=green!50!black] (-1,-0.7) .. controls +(0,0.5) and +(-1,1) .. (0,-0.5)
node[pos=0.9,dot,label=below:{$\phi$}] {};    
\draw[color=green!50!black] (-1,-0.7) .. controls +(0,-0.5) and +(0,-0.5) .. (-1.5,-0.5);
\draw[color=green!50!black]   (-1.5,-0.5) -- (-1.5,1.5)
node[pos=0.5,dot,label=right:{$\phi^{-1}$}] {};                   
%
\fill[color=blue!50!black] (0,-0.5) circle (2.9pt) node[right] (meet2) {{$\lambda^{\alpha}$}};
\fill[color=blue!50!black] (0,0.5) circle (2.9pt) node[right] (meet2) {{$\lambda_{\alpha}$}};
\draw (0,2) node[below] (X) {{$i$}};
\draw (0,-2) node[above] (X) {{$i$}};
\draw[color=green!50!black] (-0.7,1.5) node[above] (A1) {{$z$}};
\draw[color=green!50!black] (-1.5,1.5) node[above] (A1) {{$z$}};
\draw[color=blue!50!black] (0,0) node[right] (A1) {{$k$}};
\end{tikzpicture} 
-
\begin{tikzpicture}[very thick,scale=1,color=blue!50!black, baseline,
dot/.style = {circle, fill, minimum size=#1,
              inner sep=0pt, outer sep=0pt},
dot/.default = 5.8pt  
                      ]
                      \draw (0,-1.5) -- (0,1.5);
\draw[color=green!50!black] (-0.7,1.5) .. controls +(0,-0.2) and +(-0.5,0.5) .. (0,0.5)  ;   
\draw[color=green!50!black] (-1,-0.7) .. controls +(0,0.5) and +(-1,1) .. (0,-0.5)
node[pos=0.9,dot,label=below:{$\phi$}] {};    
\draw[color=green!50!black] (-1,-0.7) .. controls +(0,-0.5) and +(0,-0.5) .. (-1.5,-0.5);
\draw[color=green!50!black]   (-1.5,-0.5) -- (-1.5,1.5)
node[pos=0.5,dot,label=right:{$\phi^{-1}$}] {}; ;                     
%
\fill[color=blue!50!black] (0,-0.5) circle (2.9pt) node[right] (meet2) {{$\lambda_{\alpha}$}};
\fill[color=blue!50!black] (0,0.5) circle (2.9pt) node[right] (meet2) {{$\lambda^{\alpha}$}};
\draw (0,2) node[below] (X) {{$i$}};
\draw (0,-2) node[above] (X) {{$i$}};
\draw[color=green!50!black] (-0.7,1.5) node[above] (A1) {{$z$}};
\draw[color=green!50!black] (-1.5,1.5) node[above] (A1) {{$z$}};
\draw[color=blue!50!black] (0,0) node[right] (A1) {{$k$}};
\end{tikzpicture} 
     \right)\\
  &=
\frac{1}{d_{i}} \left( \sum_{k, \alpha}
\begin{tikzpicture}[very thick,scale=1,color=blue!50!black, baseline,
dot/.style = {circle, fill, minimum size=#1,
              inner sep=0pt, outer sep=0pt},
            dot/.default = 5.8pt  ]
                      \draw (0,-1.5) -- (0,1.5);
\draw[color=green!50!black] (-0.7,1.5) .. controls +(0,-0.2) and +(-0.5,0.5) .. (0,0.5)  ;  
\draw[color=green!50!black] (-1.5,1.5) .. controls +(0,-0.5) and +(-1,1) .. (0,-0.5);
%
\fill[color=blue!50!black] (0,-0.5) circle (2.9pt) node[right] (meet2) {{$\lambda^{\alpha}$}};
\fill[color=blue!50!black] (0,0.5) circle (2.9pt) node[right] (meet2) {{$\lambda_{\alpha}$}};
\draw (0,2) node[below] (X) {{$i$}};
\draw (0,-2) node[above] (X) {{$i$}};
\draw[color=green!50!black] (-0.7,1.5) node[above] (A1) {{$z$}};
\draw[color=green!50!black] (-1.5,1.5) node[above] (A1) {{$z$}};
\draw[color=blue!50!black] (0,0) node[right] (A1) {{$k$}};
\end{tikzpicture} 
-
\begin{tikzpicture}[very thick,scale=1,color=blue!50!black, baseline,
dot/.style = {circle, fill, minimum size=#1,
              inner sep=0pt, outer sep=0pt},
dot/.default = 5.8pt  
                      ]
                      \draw (0,-1.5) -- (0,1.5);
\draw[color=green!50!black] (-0.7,1.5) .. controls +(0,-0.2) and +(-0.5,0.5) .. (0,0.5)  ;  
\draw[color=green!50!black] (-1.5,1.5) .. controls +(0,-0.5) and +(-1,1) .. (0,-0.5);
%
\fill[color=blue!50!black] (0,-0.5) circle (2.9pt) node[right] (meet2) {{$\lambda_{\alpha}$}};
\fill[color=blue!50!black] (0,0.5) circle (2.9pt) node[right] (meet2) {{$\lambda^{\alpha}$}};
\draw (0,2) node[below] (X) {{$i$}};
\draw (0,-2) node[above] (X) {{$i$}};
\draw[color=green!50!black] (-0.7,1.5) node[above] (A1) {{$z$}};
\draw[color=green!50!black] (-1.5,1.5) node[above] (A1) {{$z$}};
\draw[color=blue!50!black] (0,0) node[right] (A1) {{$k$}};
\end{tikzpicture} 
\right)
\end{align*}
  which agree up to the irrelevant factor $\frac{1}{d_{i}}$ with the generators of the classical preprojective algebra \eqref{eq:classical-preproj} and thus the relation $\Omega_{\varphi}$ coincides with the relation of the classical preprojective algebra. 
\end{proof}
\begin{remark}
  The construction of the theorem can also be applied analogously to the case of a module category $\AM$ over $\Cat{A}$ and the quiver $Q_{z}=z \act -: \Cat{M} \rightarrow \Cat{M}$. In case of the module categories in Example \ref{example:A-actions}\refitem{item:Sl2-module}, this construction of the corresponding preprojective algebra appears already in
   \cite{MOVQuiv}. 
\end{remark}   

We finally discuss  an important class of examples.
A symmetric fusion category category $\Cat{C}$ has a canonical pivotal structure where all twists are the identities. Thus in this case, for a self dual object $x \in \Cat{C}$ with $\FS(x)=-1$, the theorem applies.

We apply this to an extended Dynkin quiver $Q$ with any choice of orientation on the edges. By Example \ref{example:first-expl} \refitem{item:ext-Dyn}, the double $D(Q)$ of $Q$ is a Drinfeld quiver $Q_{V}$ with $V \in \Rep(G)$ for a finite group $G$. Furthermore, $V$ is self-dual with $\FS(V)=-1$. We thus obtain 
 
\begin{corollary}
  Let $Q$ be an extended Dynkin quiver  with any choice of orientation on the edges. Let  $\Pi_{Q}$ be its classical preprojective algebra.
  The category $\mod(\Pi_{Q})$  has a rigid monoidal structure. 
\end{corollary}

\subsection{Preprojective algebra on the restricted double of $Q_{z}$} 

We consider a second situation where the modules over a  certain preprojective algebra carry a rigid  monoidal structure.
Therefore
we consider for a general braided pivotal fusion category $\Cat{C}$ and $x \in \cent(\Cat{C})$ the preprojective  algebra on the double of a certain subquiver of $Q_{x}$.

In a braided monoidal category $\Cat{C}$, each object $x \in \Cat{C}$ canonically corresponds to  two objects in the center $\cent(\Cat{C})$ with underlying object $x$: The object $x$ itself and
the object $\cc{x} \in \cent(\Cat{C})$  with the inverse braiding: $c_{\cc{x}}(y)=c^{-1}_{y,x}: x \otimes y \rightarrow y \otimes x$.
\begin{definition}
  \label{definition:double-z}
  Let $x \in \Cat{C}$ be an object in a braided monoidal category. The \emph{double} of $x$ is the object $Dx=x \oplus \cc{x^{*}} \in \cent(\Cat{C})$. 
\end{definition}

As in Example \ref{example:preproj-classic}, we see that for $x \in \Cat{C}$ the  double quiver $DQ_{x}$ of Definition \ref{definition:double-quiv-prepro} is isomorphic  to the quiver of the double:  $DQ_{x}\cong Q_{Dx}$. 

We use the  canonical isomorphism
\begin{equation*}
 \Hom_{\Cat{C}}(\unit, Dx \otimes Dx) \cong \Hom_{\Cat{C}}(\unit, x \otimes x) \oplus \Hom_{\Cat{C}}(\unit, x \otimes x^{*}) \oplus \Hom_{\Cat{C}}(\unit, x^{*} \otimes x) \oplus \Hom_{\Cat{C}}(\unit, x^{*} \otimes x^{*}) 
\end{equation*}
to define the following morphism in $\Cat{C}$.
  \begin{equation}
    \label{eq:can-state}
    \begin{split}
      \varphi &\in \Hom_{\Cat{C}}(\unit, Dx \otimes Dx),  \\
      \varphi &=
  \begin{tikzpicture}[very thick,scale=1,color=blue!50!black, baseline]
    \draw[color=blue!50!black] (-0.5,-0.5) .. controls +(0,-0.5) and +(0,-0.5) .. (-1.2,-0.5);
    \draw[color=blue!50!black] (-0.5,-0.5)  -- (-0.5,0.6);
    \draw[color=blue!50!black] (-1.2,-0.5)  -- (-1.2,0.6);
    %
    %
    %
    \draw[color=blue!50!black] (-0.24,-0.5) node[above] (A1) {{$x^{*}$}};
\draw[color=blue!50!black] (-1.4,-0.5) node[above] (A1) {{$x$}};
  \end{tikzpicture}
  -
  \begin{tikzpicture}[very thick,scale=1,color=blue!50!black, baseline]
    \draw[color=blue!50!black] (-0.5,-0.5) .. controls +(0,-0.5) and +(0,-0.5) .. (-1.2,-0.5);
    \draw[color=blue!50!black] (-0.5,-0.5)  -- (-0.5,-0.3);
    \draw[color=blue!50!black] (-1.2,-0.5)  -- (-1.2,-0.3);
    %
    \draw[color=blue!50!black] (-0.5,-0.3) .. controls +(0,0.4) and +(0,-0.5) .. (-1.2,0.6);
     \draw[color=white, line width=4pt]  (-1.2,-0.3) .. controls +(0,0.4) and +(0,-0.5) .. (-0.5,0.6); 
    \draw[color=blue!50!black] (-1.2,-0.3) .. controls +(0,0.4) and +(0,-0.5) .. (-0.5,0.6);
    %
    %
    \draw[color=blue!50!black] (-0.24,-0.8) node[above] (A1) {{$x^{*}$}};
    \draw[color=blue!50!black] (-0.35,0) node[above] (A1) {{$x$}};
            \draw[color=blue!50!black] (-1.4,0) node[above] (A1) {{$x^{*}$}};
\draw[color=blue!50!black] (-1.4,-0.8) node[above] (A1) {{$x$}};
  \end{tikzpicture} \; .
    \end{split}
      \end{equation}
      It is important to note, that in general, $\varphi$ is \emph{not} a morphism in $\cent(\Cat{C})$: $\coev{x}:\unit\rightarrow x \otimes x^{*}$ is a morphism
      in $\cent(\Cat{C})$, if we equip $x$ and $x^{*}$ with its original braiding, but it is not a morphism to $x \otimes \cc{x^{*}}$. However, we have
      \begin{lemma}
        \label{lemma:almost-oti-rel}
  The morphism $\varphi$ satisfies
  \begin{equation}
    \label{eq:11}
    c_{Dx}(Dx) \circ \varphi=- \varphi.
  \end{equation}
\end{lemma}
\begin{proof}
 In graphical notation, we compute in $\Cat{C}$:
  \begin{align*}
  \label{eq:4}
  \begin{tikzpicture}[very thick,scale=1,color=blue!50!black, baseline]
%
%
%
\draw[color=blue!50!black](0.5,-0.5) -- (-0.25,0.5);
\draw[color=white, line width=4pt]  (-0.5,-0.5)  -- (0.25,0.5) ;
\draw[color=blue!50!black]  (-0.5,-0.5)  -- (0.25,0.5) ;
\draw[color=blue!50!black]  (-0.25,0.5)  .. controls +(-0.1,0.1) and +(0,0) .. (-0.35,0.8);
\draw[color=blue!50!black]  (0.25,0.5)  .. controls +(0.1,0.1) and +(0,0) .. (0.35,0.8);
%
  \draw (0,-0.7) node[minimum height=0.5cm,minimum width=1cm,draw,fill=white] {{$\varphi$}};
  \draw[color=blue!50!black] (-0.35,0.8) node[above] (X) {{$Dx$}};
\draw[color=blue!50!black] (0.35,0.8) node[above] (A) {{$Dx$}};
\end{tikzpicture} 
&=
\begin{tikzpicture}[very thick,scale=1,color=blue!50!black, baseline]
\draw[color=blue!50!black] (-1.2,0.8) node[above] (X){{$\cc{x^{*}}$}}; 
\draw[color=blue!50!black] (-0.5,0.8) node[above] (A)  {{$x$}};
%
\draw[color=blue!50!black] (-1.2,0.8) -- (-1.2,-0.45); 
\draw[color=blue!50!black] (-0.5,0.8) -- (-0.5,-0.45); 
\draw[color=blue!50!black] (-0.5,-0.45) .. controls +(0,-0.5) and +(0,-0.5) .. (-1.2,-0.45);
  \draw (-0.85,0.2) node[minimum height=0.5cm,minimum width=1cm,draw,fill=white] {{$c_{x}(\cc{x^{*}})$}};
   \draw[color=blue!50!black] (-0.24,-0.8) node[above] (A1) {{$\cc{x^{*}}$}};
    \draw[color=blue!50!black] (-1.35,-0.8) node[above] (A1) {{$x$}};
%
\end{tikzpicture} 
-
\begin{tikzpicture}[very thick,scale=1,color=blue!50!black, baseline]
\draw[color=blue!50!black] (-1.2,0.8) node[above] (X) {{$x$}};
\draw[color=blue!50!black] (-0.5,0.8) node[above] (A) {{$\cc{x^{*}}$}};
%
%
%
\draw[color=blue!50!black] (-1.2,0.8) -- (-1.2,0.45); 
\draw[color=blue!50!black] (-0.5,0.8) -- (-0.5,0.45); 
    \draw[color=blue!50!black] (-0.5,-0.85) .. controls +(0,0.4) and +(0,-0.5) .. (-1.2,0.05);
    \draw[color=white, line width=4pt]  (-1.2,-0.85) .. controls +(0,0.4) and +(0,-0.5) .. (-0.5,0.05); 
    \draw[color=blue!50!black] (-1.2,-0.85) .. controls +(0,0.4) and +(0,-0.5) .. (-0.5,0.05);
    %
\draw[color=blue!50!black] (-0.5,-0.85) .. controls +(0,-0.5) and +(0,-0.5) .. (-1.2,-0.85);
\draw (-0.85,0.2) node[minimum height=0.5cm,minimum width=1cm,draw,fill=white] {{$c_{\cc{x^{*}}}(x)$}};
  \draw[color=blue!50!black] (-0.24,-1.2) node[above] (A1) {{$\cc{x^{*}}$}};
    \draw[color=blue!50!black] (-1.35,-1.2) node[above] (A1) {{$x$}};
%
\end{tikzpicture} 
  =
\begin{tikzpicture}[very thick,scale=1,color=blue!50!black, baseline]
\draw[color=blue!50!black] (-1.2,0.8) node[above] (X){{$\cc{x^{*}}$}}; 
\draw[color=blue!50!black] (-0.5,0.8) node[above] (A)  {{$x$}};
%
\draw[color=blue!50!black] (-1.2,0.8) -- (-1.2,0.75); 
\draw[color=blue!50!black] (-0.5,0.8) -- (-0.5,0.75);
    \draw[color=blue!50!black] (-0.5,-0.15) .. controls +(0,0.4) and +(0,-0.5) .. (-1.2,0.75);
    \draw[color=white, line width=4pt]  (-1.2,-0.15) .. controls +(0,0.4) and +(0,-0.5) .. (-0.5,0.75); 
    \draw[color=blue!50!black] (-1.2,-0.15) .. controls +(0,0.4) and +(0,-0.5) .. (-0.5,0.75);
    %
    \draw[color=blue!50!black] (-1.2,-0.15) -- (-1.2,-0.45); 
\draw[color=blue!50!black] (-0.5,-0.15) -- (-0.5,-0.45);
\draw[color=blue!50!black] (-0.5,-0.45) .. controls +(0,-0.5) and +(0,-0.5) .. (-1.2,-0.45);
   \draw[color=blue!50!black] (-0.24,-0.8) node[above] (A1) {{$\cc{x^{*}}$}};
    \draw[color=blue!50!black] (-1.35,-0.8) node[above] (A1) {{$x$}};
%
\end{tikzpicture} 
-
\begin{tikzpicture}[very thick,scale=1,color=blue!50!black, baseline]
\draw[color=blue!50!black] (-1.2,0.8) node[above] (X) {{$x$}};
\draw[color=blue!50!black] (-0.5,0.8) node[above] (A) {{${x^{*}}$}};
%
%
%
    \draw[color=blue!50!black] (-0.5,-0.85) .. controls +(0,0.4) and +(0,-0.4) .. (-1.2,0.05);
    \draw[color=white, line width=4pt]  (-1.2,-0.85) .. controls +(0,0.4) and +(0,-0.4) .. (-0.5,0.05); 
    \draw[color=blue!50!black] (-1.2,-0.85) .. controls +(0,0.4) and +(0,-0.4) .. (-0.5,0.05);
    %
    \draw[color=blue!50!black] (-1.2,0.05) .. controls +(0,0.4) and +(0,-0.4) .. (-0.5,0.85);   
    \draw[color=white, line width=4pt] (-0.5,0.05) .. controls +(0,0.4) and +(0,-0.4) .. (-1.2,0.85);
     \draw[color=blue!50!black] (-0.5,0.05) .. controls +(0,0.4) and +(0,-0.4) .. (-1.2,0.85);
%
\draw[color=blue!50!black] (-0.5,-0.85) .. controls +(0,-0.5) and +(0,-0.5) .. (-1.2,-0.85);
  \draw[color=blue!50!black] (-0.24,-1.2) node[above] (A1) {{$x^{*}$}};
    \draw[color=blue!50!black] (-1.35,-1.2) node[above] (A1) {{$x$}};
%
  \end{tikzpicture}
\\
&  =
\begin{tikzpicture}[very thick,scale=1,color=blue!50!black, baseline]
\draw[color=blue!50!black] (-1.2,0.8) node[above] (X){{$\cc{x^{*}}$}}; 
\draw[color=blue!50!black] (-0.5,0.8) node[above] (A)  {{$x$}};
%
\draw[color=blue!50!black] (-1.2,0.8) -- (-1.2,0.75); 
\draw[color=blue!50!black] (-0.5,0.8) -- (-0.5,0.75);
    \draw[color=blue!50!black] (-0.5,-0.15) .. controls +(0,0.4) and +(0,-0.5) .. (-1.2,0.75);
    \draw[color=white, line width=4pt]  (-1.2,-0.15) .. controls +(0,0.4) and +(0,-0.5) .. (-0.5,0.75); 
    \draw[color=blue!50!black] (-1.2,-0.15) .. controls +(0,0.4) and +(0,-0.5) .. (-0.5,0.75);
    %
    \draw[color=blue!50!black] (-1.2,-0.15) -- (-1.2,-0.45); 
\draw[color=blue!50!black] (-0.5,-0.15) -- (-0.5,-0.45);
\draw[color=blue!50!black] (-0.5,-0.45) .. controls +(0,-0.5) and +(0,-0.5) .. (-1.2,-0.45);
   \draw[color=blue!50!black] (-0.24,-0.8) node[above] (A1) {{$\cc{x^{*}}$}};
    \draw[color=blue!50!black] (-1.35,-0.8) node[above] (A1) {{$x$}};
\end{tikzpicture} 
-
\begin{tikzpicture}[very thick,scale=1,color=blue!50!black, baseline]
\draw[color=blue!50!black] (-1.2,0.8) node[above] (X) {{$x$}};
\draw[color=blue!50!black] (-0.5,0.8) node[above] (A) {{${x^{*}}$}};
 \draw[color=blue!50!black] (-1.2,0.8) -- (-1.2,-0.85); 
 \draw[color=blue!50!black] (-0.5,0.8) -- (-0.5,-0.85); 
%
\draw[color=blue!50!black] (-0.5,-0.85) .. controls +(0,-0.5) and +(0,-0.5) .. (-1.2,-0.85);
  \draw[color=blue!50!black] (-0.24,-0.9) node[above] (A1) {{$x^{*}$}};
    \draw[color=blue!50!black] (-1.35,-0.9) node[above] (A1) {{$x$}};
  \end{tikzpicture}
     =-\varphi \, ,
\end{align*}
where we used the definition of the half-braiding of $\cc{x^{*}}$ in the second step. 
\end{proof}
With the twist \ref{eq:twist-def}, 
we get

\begin{equation}
  \label{eq:8}
  \varphi=
  \begin{tikzpicture}[very thick,scale=1,color=blue!50!black, baseline]
    \draw[color=blue!50!black] (-0.5,-0.5) .. controls +(0,-0.5) and +(0,-0.5) .. (-1.2,-0.5);
    \draw[color=blue!50!black] (-0.5,-0.5)  -- (-0.5,0.6);
    \draw[color=blue!50!black] (-1.2,-0.5)  -- (-1.2,0.6);
    %
    %
    %
    \draw[color=blue!50!black] (-0.24,-0.5) node[above] (A1) {{$x^{*}$}};
\draw[color=blue!50!black] (-1.4,-0.5) node[above] (A1) {{$x$}};
  \end{tikzpicture}
  -
  \begin{tikzpicture}[very thick,scale=1,color=blue!50!black, baseline]
    \draw[color=blue!50!black] (-0.5,-0.5) .. controls +(0,-0.5) and +(0,-0.5) .. (-1.2,-0.5);
    \draw[color=blue!50!black] (-0.5,-0.5)  -- (-0.5,0.6);
    \draw[color=blue!50!black] (-1.2,-0.5)  -- (-1.2,0.6);
    %
    %
    %
    \draw[color=blue!50!black] (-0.3,-0.7) node[above] (A1) {{$x$}};
        \fill[color=blue!50!black]   (-0.5,0) circle (2.9pt) node[right] (meet2) {{$\theta_x$}};
\draw[color=blue!50!black] (-1.4,-0.7) node[above] (A1) {{$x^{*}$}};
  \end{tikzpicture}
\end{equation}
which for a simple $x \in \Cat{C}$ reduces to
\begin{equation}
  \label{eq:8}
  \varphi=
  \begin{tikzpicture}[very thick,scale=1,color=blue!50!black, baseline=-0.5cm]
    \draw[color=blue!50!black] (-0.5,-0.5) .. controls +(0,-0.5) and +(0,-0.5) .. (-1.2,-0.5);
    \draw[color=blue!50!black] (-0.5,-0.5)  -- (-0.5,0.2);
    \draw[color=blue!50!black] (-1.2,-0.5)  -- (-1.2,0.2);
    %
    %
    %
    \draw[color=blue!50!black] (-0.24,-0.5) node[above] (A1) {{$x^{*}$}};
\draw[color=blue!50!black] (-1.4,-0.5) node[above] (A1) {{$x$}};
  \end{tikzpicture}
  - \theta
  \begin{tikzpicture}[very thick,scale=1,color=blue!50!black, baseline=-0.5cm]
    \draw[color=blue!50!black] (-0.5,-0.5) .. controls +(0,-0.5) and +(0,-0.5) .. (-1.2,-0.5);
    \draw[color=blue!50!black] (-0.5,-0.5)  -- (-0.5,0.2);
    \draw[color=blue!50!black] (-1.2,-0.5)  -- (-1.2,0.2);
    %
    %
    %
    \draw[color=blue!50!black] (-0.3,-0.7) node[above] (A1) {{$x$}};
     \draw[color=blue!50!black] (-1.4,-0.7) node[above] (A1) {{$x^{*}$}};
  \end{tikzpicture}\; ,
\end{equation}
with the twist $\theta=\theta_{x} \in \Bbbk^{\times}$. 

Recall from Definition \ref{definition:cl-prepro} the twisted preprojective algebra. 
\begin{lemma}
  \label{lemma:modules-twisted-prepro}
  Let $x \in \Cat{C}$ be a simple object and $\varphi$ be as above. A module $(m,f) \in \mod(Q_{Dx})$  satisfies the relation $\varphi$, if and only if it is a module over the $\theta_{x}$-twisted preprojective algebra of $Q_{x}$.
\end{lemma}
\begin{proof}
  We define for all $i,j \in I$ the non-degenerate pairing
  \begin{align*}
    \label{eq:paring-prepro}
    <-.-> : \; &\Hom_{\Cat{C}}(x_{j}, x \otimes x_{i}) \otik \Hom_{\Cat{C}}(x_{i}, x^{*} \otimes x_{j} ) \rightarrow \Bbbk, \\
    &(\lambda,\mu) \mapsto
     \frac{1}{d_{i}d_{j}} \cdot \;
\begin{tikzpicture}[very thick,scale=1,color=blue!50!black, baseline]
\draw (0,-1) -- (0,1); 
\draw[color=green!50!black] (-1.4,0.7) .. controls +(0,0.5) and +(-0.5,0.5) .. (0,0.7); 
\draw[color=green!50!black] (-1.4,0.7) .. controls +(0,-0.5) and +(-0.5,0.5) .. (0,-0.7);
\draw[color=blue!50!black] (0,1) .. controls +(0,1) and +(0,1) .. (-2,1); 
\draw[color=blue!50!black] (0,-1) .. controls +(0,-1) and +(0,-1) .. (-2,-1); 
\draw[color=blue!50!black] (-2,1) -- (-2,-1); 
%
\fill[color=blue!50!black] (0,0.7) circle (2.9pt) node[right] (meet1) {{$\lambda$}};
\fill[color=blue!50!black] (0,-0.7) circle (2.9pt) node[right] (meet2) {{$\mu$}};
\draw (-0.1,-1.2) node[right] (X) {{$i$}};
\draw (-0.1,-0) node[right] (X) {{$j$}};
\draw (-0.1,1.2) node[right] (X) {{$i$}};
\draw[color=green!50!black] (-0.4,-0.4) node[left] (A1) {{$x^{*}$}};
\draw[color=green!50!black] (-0.4,0.9) node[left] (A1) {{$x$}};;
\end{tikzpicture} \, .
  \end{align*}
  By cyclicity of the trace, $<\lambda,\mu>=<\mu,\lambda>$
  and if we pick an orthonormal basis $\{\lambda_{i,j, \alpha}\}_{\alpha} \subset \Hom(x_{j}, x \otimes x_{j})$ and $\{ \lambda_{i,j}^{\beta}\}_{\beta} \subset \Hom_{\Cat{C}}(x_{i}, x^{*} \otimes x_{j})$ with
  $<\lambda_{i,j,\alpha}, \lambda_{i,j}^{\beta}>=\delta_{\alpha}^{\beta}$ for all $i,j \in I$ and all $\alpha,\beta$, we obtain
  as in the Proof of Theorem \ref{theorem:main-rel},
  \begin{equation}
  \label{eq:11}
  \begin{tikzpicture}[very thick,scale=1,color=blue!50!black, baseline,
    dot/.style = {circle, fill, minimum size=#1,
              inner sep=0pt, outer sep=0pt},
dot/.default = 5.8pt  
                     ]
\draw (0,-1.5) -- (0,1.5); 
\draw[color=green!50!black] (-1.4,0.7) .. controls +(0,0.5) and +(-0.5,0.5) .. (0,0.5);
\draw[color=green!50!black] (-1.4,0.7) .. controls +(0,-0.5) and +(-0.5,0.5) .. (0,-0.5);
%
%
\fill[color=blue!50!black] (0,0.5) circle (2.9pt) node[right] (meet1) {{$\lambda_{\alpha}$}};
\fill[color=blue!50!black] (0,-0.5) circle (2.9pt) node[right] (meet2) {{$\lambda^{\beta}$}};
\draw (-0.1,-1.2) node[right] (X) {{$i$}};
\draw (-0.1,-0) node[right] (X) {{$j$}};
\draw (-0.1,1.2) node[right] (X) {{$i$}};
\draw[color=green!50!black] (-0.9,0.1) node[below] (A1) {{$x^{*}$}};
\draw[color=green!50!black] (-0.7,0.9) node[below] (A1) {{$x$}};
\end{tikzpicture}
=\delta^{\alpha}_{\beta} \, d_{j} \cdot
\begin{tikzpicture}[very thick,scale=1,color=blue!50!black, baseline,
    dot/.style = {circle, fill, minimum size=#1,
              inner sep=0pt, outer sep=0pt},
dot/.default = 5.8pt  
                     ]
\draw (0,-1.5) -- (0,1.5); 
\draw (-0.1,0) node[right] (X) {{$i$}};
\end{tikzpicture} 
\quad, \text{and} \quad \quad
%
  \begin{tikzpicture}[very thick,scale=1,color=blue!50!black, baseline,
    dot/.style = {circle, fill, minimum size=#1,
              inner sep=0pt, outer sep=0pt},
dot/.default = 5.8pt  
                     ]
\draw (0,-1.5) -- (0,1.5); 
\draw[color=green!50!black] (-1.4,0.7) .. controls +(0,0.5) and +(-0.5,0.5) .. (0,0.5);
\draw[color=green!50!black] (-1.4,0.7) .. controls +(0,-0.5) and +(-0.5,0.5) .. (0,-0.5);
%
%
\fill[color=blue!50!black] (0,0.5) circle (2.9pt) node[right] (meet1) {{$\lambda^{\beta}$}};
\fill[color=blue!50!black] (0,-0.5) circle (2.9pt) node[right] (meet2) {{$\lambda_{\alpha}$}};
\draw (-0.1,-1.2) node[right] (X) {{$i$}};
\draw (-0.1,-0) node[right] (X) {{$j$}};
\draw (-0.1,1.2) node[right] (X) {{$i$}};
\draw[color=green!50!black] (-0.9,0.1) node[below] (A1) {{$x$}};
\draw[color=green!50!black] (-0.7,0.9) node[below] (A1) {{$x^{*}$}};
\end{tikzpicture}
=\delta^{\alpha}_{\beta} \, d_{j} \cdot
\begin{tikzpicture}[very thick,scale=1,color=blue!50!black, baseline,
    dot/.style = {circle, fill, minimum size=#1,
              inner sep=0pt, outer sep=0pt},
dot/.default = 5.8pt  
                     ]
\draw (0,-1.5) -- (0,1.5); 
\draw (-0.1,0) node[right] (X) {{$i$}};
\end{tikzpicture}
\end{equation}
and thus also
\begin{equation}
  \label{eq:12}
\begin{tikzpicture}[very thick,scale=1,color=blue!50!black, baseline,
    dot/.style = {circle, fill, minimum size=#1,
              inner sep=0pt, outer sep=0pt},
dot/.default = 5.8pt  
                     ]
\draw[color=green!50!black] (0,-1.5) -- (0,1.5); 
\draw[color=green!50!black] (-0.1,0) node[right] (X) {{$x$}};
\draw (0.8,-1.5) -- (0.8,1.5); 
\draw (0.8-0.1,0) node[right] (X) {{$i$}};

\end{tikzpicture}
  = \frac{1}{d_{i}}  \sum_{j,\alpha}
\begin{tikzpicture}[very thick,scale=1,color=blue!50!black, baseline,
dot/.style = {circle, fill, minimum size=#1,
              inner sep=0pt, outer sep=0pt},
dot/.default = 5.8pt  
                      ]
\draw (0,-1.5) -- (0,1.5);  
\draw[color=green!50!black] (-0.7,1.5) .. controls +(0,-0.2) and +(-0.5,0.5) .. (0,0.5)  ;   
\draw[color=green!50!black] (-1,-1.5) .. controls +(0,0.5) and +(-1,1) .. (0,-0.5);
\fill[color=blue!50!black] (0,-0.5) circle (2.9pt) node[right] (meet2) {{$\lambda_{\alpha}$}};
\fill[color=blue!50!black] (0,0.5) circle (2.9pt) node[right] (meet2) {{$\lambda^{\alpha}$}};
\draw (0,2) node[below] (X) {{$i$}};
\draw (0,-2) node[above] (X) {{$i$}};
\draw[color=green!50!black] (-0.7,2) node[below] (A1) {{$x$}};
\draw[color=green!50!black] (-1,-2) node[above] (A1) {{$x$}};
\draw[color=blue!50!black] (0,0) node[right] (A1) {{$j$}};
\end{tikzpicture}
\quad , \text{and} \quad \quad
\begin{tikzpicture}[very thick,scale=1,color=blue!50!black, baseline,
    dot/.style = {circle, fill, minimum size=#1,
              inner sep=0pt, outer sep=0pt},
dot/.default = 5.8pt  
                     ]
\draw[color=green!50!black] (0,-1.5) -- (0,1.5); 
\draw[color=green!50!black] (-0.1,0) node[right] (X) {{$x^{*}$}};
\draw (0.8,-1.5) -- (0.8,1.5); 
\draw (0.8-0.1,0) node[right] (X) {{$i$}};

\end{tikzpicture}
  = \frac{1}{d_{i}}  \sum_{j,\alpha}
\begin{tikzpicture}[very thick,scale=1,color=blue!50!black, baseline,
dot/.style = {circle, fill, minimum size=#1,
              inner sep=0pt, outer sep=0pt},
dot/.default = 5.8pt  
                      ]
\draw (0,-1.5) -- (0,1.5);  
\draw[color=green!50!black] (-0.7,1.5) .. controls +(0,-0.2) and +(-0.5,0.5) .. (0,0.5)  ;   
\draw[color=green!50!black] (-1,-1.5) .. controls +(0,0.5) and +(-1,1) .. (0,-0.5);
\fill[color=blue!50!black] (0,-0.5) circle (2.9pt) node[right] (meet2) {{$\lambda^{\alpha}$}};
\fill[color=blue!50!black] (0,0.5) circle (2.9pt) node[right] (meet2) {{$\lambda_{\alpha}$}};
\draw (0,2) node[below] (X) {{$i$}};
\draw (0,-2) node[above] (X) {{$i$}};
\draw[color=green!50!black] (-0.7,2) node[below] (A1) {{$x^{*}$}};
\draw[color=green!50!black] (-1,-2) node[above] (A1) {{$x^{*}$}};
\draw[color=blue!50!black] (0,0) node[right] (A1) {{$j$}};
\end{tikzpicture}\; . 
\end{equation}
We conclude, that for all $i \in I$,

\begin{align*}
  \label{eq:8}
  \begin{tikzpicture}[very thick,scale=1,color=blue!50!black, baseline]
    \draw[color=green!50!black](0.4,-0.5) -- (0.4,1.5);
    \draw[color=green!50!black]  (-0.4,-0.5)  -- (-0.4,1.5) ;
    %
    \draw[color=blue!50!black](1.1,-2) -- (1.1,1.5);
    \draw[color=green!50!black]  (0,-0.7) node[minimum height=0.5cm,minimum width=1cm,draw,fill=white] {{$\varphi$}};
    %
    \draw[color=green!50!black] (0.4,1) node[left] (A1) {{$z$}};
\draw[color=green!50!black] (-0.4,1) node[left] (A1) {{$z$}};
    \draw[color=blue!50!black] (1.3,0) node[below] (A1) {{$i$}};
\end{tikzpicture}
=&
  \begin{tikzpicture}[very thick,scale=1,color=blue!50!black, baseline]
    \draw[color=blue!50!black] (-0.5,-0.5) .. controls +(0,-0.5) and +(0,-0.5) .. (-1.2,-0.5);
    \draw[color=blue!50!black] (-0.5,-0.5)  -- (-0.5,1.5);
    \draw[color=blue!50!black] (-1.2,-0.5)  -- (-1.2,1.5);
      \draw[color=blue!50!black](0.2,-2) -- (0.2,1.5);
    %
    %
    %
    \draw[color=blue!50!black] (-0.24,-0.5) node[above] (A1) {{$x^{*}$}};
    \draw[color=blue!50!black] (-1.4,-0.5) node[above] (A1) {{$x$}};
        \draw[color=blue!50!black] (0.4,0) node[below] (A1) {{$i$}};
  \end{tikzpicture}
  - \theta
  \begin{tikzpicture}[very thick,scale=1,color=blue!50!black, baseline]
    \draw[color=blue!50!black] (-0.5,-0.5) .. controls +(0,-0.5) and +(0,-0.5) .. (-1.2,-0.5);
    \draw[color=blue!50!black] (-0.5,-0.5)  -- (-0.5,1.5);
    \draw[color=blue!50!black] (-1.2,-0.5)  -- (-1.2,1.5);
          \draw[color=blue!50!black](0.2,-2) -- (0.2,1.5);
    %
    %
    %
    \draw[color=blue!50!black] (-0.3,-0.7) node[above] (A1) {{$x$}};
    \draw[color=blue!50!black] (-1.4,-0.7) node[above] (A1) {{$x^{*}$}};
            \draw[color=blue!50!black] (0.4,0) node[below] (A1) {{$i$}};
  \end{tikzpicture}
  =
\frac{1}{d_{i}} \left( \sum_{k, \alpha}
\begin{tikzpicture}[very thick,scale=1,color=blue!50!black, baseline,
dot/.style = {circle, fill, minimum size=#1,
              inner sep=0pt, outer sep=0pt},
            dot/.default = 5.8pt  ]
                      \draw (0,-1.5) -- (0,1.5);
\draw[color=green!50!black] (-0.7,1.5) .. controls +(0,-0.2) and +(-0.5,0.5) .. (0,0.5)  ;  
\draw[color=green!50!black] (-1.5,1.5) .. controls +(0,-0.5) and +(-1,1) .. (0,-0.5);
%
\fill[color=blue!50!black] (0,-0.5) circle (2.9pt) node[right] (meet2) {{$\lambda^{\alpha}$}};
\fill[color=blue!50!black] (0,0.5) circle (2.9pt) node[right] (meet2) {{$\lambda_{\alpha}$}};
\draw (0,2) node[below] (X) {{$i$}};
\draw (0,-2) node[above] (X) {{$i$}};
\draw[color=green!50!black] (-0.7,1.5) node[above] (A1) {{$z$}};
\draw[color=green!50!black] (-1.5,1.5) node[above] (A1) {{$z$}};
\draw[color=blue!50!black] (0,0) node[right] (A1) {{$k$}};
\end{tikzpicture} 
- \theta
\begin{tikzpicture}[very thick,scale=1,color=blue!50!black, baseline,
dot/.style = {circle, fill, minimum size=#1,
              inner sep=0pt, outer sep=0pt},
dot/.default = 5.8pt  
                      ]
                      \draw (0,-1.5) -- (0,1.5);
\draw[color=green!50!black] (-0.7,1.5) .. controls +(0,-0.2) and +(-0.5,0.5) .. (0,0.5)  ;  
\draw[color=green!50!black] (-1.5,1.5) .. controls +(0,-0.5) and +(-1,1) .. (0,-0.5);
%
\fill[color=blue!50!black] (0,-0.5) circle (2.9pt) node[right] (meet2) {{$\lambda_{\alpha}$}};
\fill[color=blue!50!black] (0,0.5) circle (2.9pt) node[right] (meet2) {{$\lambda^{\alpha}$}};
\draw (0,2) node[below] (X) {{$i$}};
\draw (0,-2) node[above] (X) {{$i$}};
\draw[color=green!50!black] (-0.7,1.5) node[above] (A1) {{$z$}};
\draw[color=green!50!black] (-1.5,1.5) node[above] (A1) {{$z$}};
\draw[color=blue!50!black] (0,0) node[right] (A1) {{$k$}};
\end{tikzpicture} 
   \right)
  \\
  =&\frac{1}{d_{i}} \omega_{i}^{\theta},
\end{align*}
  using the path $\omega_{i}^{\theta}$ from Equation \eqref{eq:twisted-prepro}. Thus, the ideal corresponding to $\varphi$ coincides with the ideal of the $\theta_{x}$-twisted preprojective algebra from Definition \ref{definition:cl-prepro} \refitem{item:theta-prepro}. 
 \end{proof}
    
Since $\varphi$ is in general not a morphism in $\cent(\Cat{C})$, in general, we do not get a rigid monoidal structure on all modules that satisfy the relation $\varphi$ by our general construction.      
We remedy this issue by restricting the quiver $Q_{x}$ to the full subquiver on the vertices that projectively centralize $x$ \cite[Def 8.22.4]{EGNObook}. 
\begin{definition}
  Let $\Cat{C}$ be a braided monoidal and $x \in \Cat{C}$.
  An object $y \in \Cat{C} $   \emph{projective centralizes $x$}, if there is a $\lambda \in \Bbbk^{\times}$, such
  that $c_{y,x} \circ c_{x,y}=\lambda \cdot \id_{x \otimes y}$.  
  The \emph{projective centralizer} $\PCent_{x}$ of $x$ is the full subcategory of $\Cat{C}$ on the objects $y \in \Cat{C}$ that projective centralize $x$. 
\end{definition}
  It follows from the definition that $\PCent_{x}$ is a full braided monoidal subcategory of $\Cat{C}$, which is rigid (semisimple) if $\Cat{C}$ is rigid (semisimple).
Furthermore, an object $y$ projective centralizes $x$ if and only if it projective  centralizes $x^{*}$.

  Now suppose that $\Cat{C}$ is again a braided pivotal fusion category. Denote by $I_{\Cent_{x}} \subset I$ the representatives for the simple objects of $\PCent_{x}$, i.e. $I_{\Cent_{x}}$ represents those simple objects of $\Cat{C}$ that projective centralize $x$. 
  \begin{definition}
    Let $x \in \Cat{C}$ with double $z = Dx \in \cent(\Cat{C})$. 
    The \emph{double centralizer quiver of $x$} is the restriction of $Q_{z}: \Cat{C} \rightarrow \Cat{C}$ to the projective centralizer of $x$,
    \begin{equation}
      \label{eq:68}
      \CentQ_{z}: \PCent_{x} \rightarrow \PCent_{x}. 
    \end{equation}
    \end{definition}
  That is, $\CentQ_{z}$ has vertices $i,j \in I_{\Cent_{x}}$ and arrow spaces $\CentQ_{z}(i,j)= \Hom_{\Cat{C}}(x_{i}, z \otimes x_{j})$ for $i,j \in I_{\Cent_{x}}$.

 For the morphism $\varphi$ from Equation \eqref{eq:can-state}   we obtain categories $\mod(Q_{z};\varphi)$ and $\mod(\CentQ_{z};\varphi)$ of modules over $Q_{z}$ and
  $\CentQ_{z}$ that satisfy the canonical $\varphi$-relation, respectively.
  \begin{theorem}
    \label{theorem:double-cent-rigid-mon}
    Let $\Cat{C}$ be a braided pivotal fusion category and $x \in \cent(\Cat{C})$ with double $z=Dx$ and double centralizer $ \CentQ_{z}: \PCent_{x} \rightarrow \PCent_{x}$. 
    \begin{theoremlist}
    \item The category of modules  $\mod(\CentQ_{z};\varphi)$ is a rigid monoidal category with the $\qoti$ monoidal product.
    \item The category of modules $\mod(Q_{z};\varphi)$ is a left module category over $\mod(\CentQ_{z};\varphi)$: The $\qoti$-product provides an action functor
      \begin{equation}
        \label{eq:10}
        \qoti: \mod(\CentQ_{z};\varphi) \times \mod(Q_{z};\varphi) \rightarrow \mod(Q_{z};\varphi).
      \end{equation}
    \end{theoremlist}
  \end{theorem}
  \begin{proof}
    For the first part we need to show, that for modules $(a,f)$ and $(b,g)$ of $\CentQ_{z}$ which satisfy the relation $\varphi$, also $f \qoti g$ satisfies the relation $\varphi$.
    The morphism $f: Dx \otimes a \rightarrow a$ is the direct sum of two morphism $f_{1}:x \otimes a \rightarrow a$ and $f_{2}: \cc{x} \otimes a \rightarrow a$. Analogously we have $g=g_{1}\oplus g_{2}$  for
    $g_{1}: x \otimes b \rightarrow b$
    and $g_{2}: \cc{x} \otimes b \rightarrow b$.    
 We need to show that
\begin{equation}
  \label{eq:varphi-q-rel}
\begin{tikzpicture}[very thick,scale=1,color=blue!50!black, baseline]
  \draw (0,-1.8) -- (0,1.3); 
  \draw (0.7,-1.8) -- (0.7,1.3); 
\draw[color=green!50!black] (-0.7,-0.75) .. controls +(0,0.5) and +(-0.5,-0.5) .. (0,0.15); 
\draw[color=green!50!black] (-1.2,-0.75) .. controls +(0,0.5) and +(-0.5,-0.5) .. (0,0.7);
%
\fill[color=blue!50!black] (0,0.15) circle (2.9pt) node[right] (meet2) {{$f$}};
\fill[color=blue!50!black] (0,0.7) circle (2.9pt) node[right] (meet2) {{$f$}};
 \draw (-0.95,-1) node[minimum height=0.4cm,minimum width=0.7cm,draw,fill=white] {{$\varphi$}};
%
 \draw (0,-2.3) node[above] (X) {{$a$}};
 \draw (0.7,-2.3) node[above] (X) {{$b$}};
\draw[color=green!50!black] (-0.25,-0.55) node[above] (A1) {{$z$}};
\draw[color=green!50!black] (-0.83,-0.1) node[above] (A1) {{$z$}};
\end{tikzpicture} 
\, +
\begin{tikzpicture}[very thick,scale=1,color=blue!50!black, baseline]
  \draw (0.7,-1.8) -- (0.7,1.3); 
    \draw (0,-1.8) -- (0,1.3); 
    \draw[color=white, line width=4pt] (-0.7,-0.75) .. controls +(0,0.5) and +(-0.5,-0.5) .. (0.7,0.15);  
        \draw[color=white, line width=4pt]  (-1.2,-0.75) .. controls +(0,0.5) and +(-0.5,-0.5) .. (0.7,0.7);  
  \draw[color=green!50!black] (-0.7,-0.75) .. controls +(0,0.5) and +(-0.5,-0.5) .. (0.7,0.15); 
\draw[color=green!50!black] (-1.2,-0.75) .. controls +(0,0.5) and +(-0.5,-0.5) .. (0.7,0.7);
%

\fill[color=blue!50!black] (0.7,0.15) circle (2.9pt) node[right] (meet2) {{$g$}};
\fill[color=blue!50!black] (0.7,0.7) circle (2.9pt) node[right] (meet2) {{$g$}};
 \draw (-0.95,-1) node[minimum height=0.4cm,minimum width=0.7cm,draw,fill=white] {{$\varphi$}};
%
 \draw (0,-2.3) node[above] (X) {{$a$}};
 \draw (0.7,-2.3) node[above] (X) {{$b$}};
\draw[color=green!50!black] (-0.25,-0.77) node[above] (A1) {{$z$}};
\draw[color=green!50!black] (-0.66,-0.2) node[above] (A1) {{$z$}};
\end{tikzpicture} 
+
\begin{tikzpicture}[very thick,scale=1,color=blue!50!black, baseline]
  \draw (0.7,-1.8) -- (0.7,1.3); 
  \draw (0,-1.8) -- (0,1.3); 
  \draw[color=green!50!black] (-0.7,-0.75) .. controls +(0,0.5) and +(-0.5,-0.5) .. (0,0.15); 
  \draw[color=white, line width=4pt]  (-1.2,-0.75) .. controls +(0,1) and +(-1,-0.5) .. (0.7,0.7);   
\draw[color=green!50!black] (-1.2,-0.75) .. controls +(0,1) and +(-1,-0.5) .. (0.7,0.7);
%
\fill[color=blue!50!black] (0,0.15) circle (2.9pt) node[right] (meet2) {{$f$}};
\fill[color=blue!50!black] (0.7,0.7) circle (2.9pt) node[right] (meet2) {{$g$}};
 \draw (-0.95,-1) node[minimum height=0.4cm,minimum width=0.7cm,draw,fill=white] {{$\varphi$}};
%
 \draw (0,-2.3) node[above] (X) {{$a$}};
 \draw (0.7,-2.3) node[above] (X) {{$b$}};
\draw[color=green!50!black] (-0.25,-0.5) node[above] (A1) {{$z$}};
\draw[color=green!50!black] (-0.66,0.1) node[above] (A1) {{$z$}};
\end{tikzpicture} 
\, +
\begin{tikzpicture}[very thick,scale=1,color=blue!50!black, baseline]
  \draw (0.7,-1.8) -- (0.7,1.3); 
 \draw (0,-1.8) -- (0,1.3); 
   \draw[color=white, line width=4pt]   (-0.7,-0.75) .. controls +(0,0.8) and +(-0.5,-0.5) .. (0.7,0.15);      
  \draw[color=green!50!black] (-0.7,-0.75) .. controls +(0,0.8) and +(-0.5,-0.5) .. (0.7,0.15); 
\draw[color=green!50!black] (-1.2,-0.75) .. controls +(0,0.5) and +(-0.5,-0.5) .. (0,0.7);
%
\fill[color=blue!50!black] (0.7,0.15) circle (2.9pt) node[right] (meet2) {{$g$}};
\fill[color=blue!50!black] (0,0.7) circle (2.9pt) node[right] (meet2) {{$f$}};
 \draw (-0.95,-1) node[minimum height=0.4cm,minimum width=0.7cm,draw,fill=white] {{$\varphi$}};
%
 \draw (0,-2.3) node[above] (X) {{$a$}};
 \draw (0.7,-2.3) node[above] (X) {{$b$}};
\draw[color=green!50!black] (-0.25,-0.6) node[above] (A1) {{$z$}};
\draw[color=green!50!black] (-0.66,0.1) node[above] (A1) {{$z$}};
\end{tikzpicture}
=0
\end{equation}
    By Lemma \ref{lemma:almost-oti-rel}, the sum of the last two terms is zero as in the proof of Theorem \ref{theorem:main-rel}. Since $f$ satisfy relation $\varphi$, the first term in Equation \eqref{eq:varphi-q-rel} vanishes.
    The second term computes as
      \begin{align*}
\begin{tikzpicture}[very thick,scale=1,color=blue!50!black, baseline]
  \draw (0.7,-1.8) -- (0.7,1.3); 
  \draw (0,-1.8) -- (0,1.3); 
  \draw[color=white, line width=4pt] (-0.7,-0.75) .. controls +(0,0.5) and +(-0.5,-0.5) .. (0.7,0.15); 
  \draw[color=green!50!black] (-0.7,-0.75) .. controls +(0,0.5) and +(-0.5,-0.5) .. (0.7,0.15); 
  \draw[color=white, line width=4pt]  (-1.2,-0.75) .. controls +(0,0.5) and +(-0.5,-0.5) .. (0.7,0.7);
\draw[color=green!50!black] (-1.2,-0.75) .. controls +(0,0.5) and +(-0.5,-0.5) .. (0.7,0.7);
%
\fill[color=blue!50!black] (0.7,0.15) circle (2.9pt) node[right] (meet2) {{$g$}};
\fill[color=blue!50!black] (0.7,0.7) circle (2.9pt) node[right] (meet2) {{$g$}};
 \draw (-0.95,-1) node[minimum height=0.4cm,minimum width=0.7cm,draw,fill=white] {{$\varphi$}};
%
 \draw (0,-2.3) node[above] (X) {{$a$}};
 \draw (0.7,-2.3) node[above] (X) {{$b$}};
\draw[color=green!50!black] (-0.25,-0.77) node[above] (A1) {{$z$}};
\draw[color=green!50!black] (-0.66,-0.2) node[above] (A1) {{$z$}}; 
\end{tikzpicture} 
&=
\begin{tikzpicture}[very thick,scale=1,color=blue!50!black, baseline]
  \draw (0.7,-1.8) -- (0.7,1.3); 
    \draw (0,-1.8) -- (0,1.3); 
  \draw[color=white, line width=4pt] (0.7,0.7) --  (0.7-1,0.7-1); 
  \draw[color=green!50!black] (0.7,0.7) --  (0.7-1,0.7-1); 
  \draw[color=white, line width=4pt] (0.7,0.15) --  (0.7-1,0.15-1); 
    \draw[color=green!50!black] (0.7,0.15) --  (0.7-1,0.15-1); 
   \draw[color=green!50!black]  (0.7-1,0.7-1)  .. controls +(-0.5,-0.5) and +(-0.5,-0.5) .. (0.7-1,0.15-1);
\fill[color=blue!50!black] (0.7,0.15) circle (2.9pt) node[right] (meet2) {{$g_1$}};
\fill[color=blue!50!black] (0.7,0.7) circle (2.9pt) node[right] (meet2) {{$g_2$}};
%
 \draw (0,-2.5) node[above] (X) {{$a$}};
 \draw (0.7,-2.5) node[above] (X) {{$b$}};
 \draw[color=green!50!black] (0.3,0.25) node[above] (A1) {{$x$}};
 \draw[color=green!50!black] (0.4,-0.96) node[above] (A1) {{$\cc{x^*}$}};
\end{tikzpicture} 
- \theta
\begin{tikzpicture}[very thick,scale=1,color=blue!50!black, baseline]
  \draw (0.7,-1.8) -- (0.7,1.3); 
    \draw (0,-1.8) -- (0,1.3); 
  \draw[color=white, line width=4pt]  (0.7,0.7) --  (0.7-1,0.7-1); 
  \draw[color=green!50!black] (0.7,0.7) --  (0.7-1,0.7-1); 
  \draw[color=white, line width=4pt]  (0.7,0.15) --  (0.7-1,0.15-1); 
    \draw[color=green!50!black] (0.7,0.15) --  (0.7-1,0.15-1); 
   \draw[color=green!50!black]  (0.7-1,0.7-1)  .. controls +(-0.5,-0.5) and +(-0.5,-0.5) .. (0.7-1,0.15-1);
%
\fill[color=blue!50!black] (0.7,0.15) circle (2.9pt) node[right] (meet2) {{$g_2$}};
\fill[color=blue!50!black] (0.7,0.7) circle (2.9pt) node[right] (meet2) {{$g_1$}};
%
 \draw (0,-2.5) node[above] (X) {{$a$}};
 \draw (0.7,-2.5) node[above] (X) {{$b$}};
 \draw[color=green!50!black] (0.3,0.25) node[above] (A1)  {{$\cc{x^*}$}};
 \draw[color=green!50!black] (0.4,-0.8) node[above] (A1) {{$x$}};
\end{tikzpicture} 
=
\begin{tikzpicture}[very thick,scale=1,color=blue!50!black, baseline]
  \draw (0.7,-1.8) -- (0.7,1.3); 
  \draw[color=green!50!black] (0.7,0.15) --  (0.7-1,0.15-1); 
   \draw[color=white, line width=4pt] (0,-1.8) -- (0,1.3);   
    \draw (0,-1.8) -- (0,1.3); 
     \draw[color=white, line width=4pt] (0.7,0.7) --  (0.7-1,0.7-1); 
  \draw[color=green!50!black] (0.7,0.7) --  (0.7-1,0.7-1); 
         \draw[color=green!50!black]  (0.7-1,0.7-1)  .. controls +(-0.5,-0.5) and +(-0.5,-0.5) .. (0.7-1,0.15-1);
\fill[color=blue!50!black] (0.7,0.15) circle (2.9pt) node[right] (meet2) {{$g_1$}};
\fill[color=blue!50!black] (0.7,0.7) circle (2.9pt) node[right] (meet2) {{$g_2$}};
%
 \draw (0,-2.5) node[above] (X) {{$a$}};
 \draw (0.7,-2.5) node[above] (X) {{$b$}};
 \draw[color=green!50!black] (0.3,0.25) node[above] (A1) {{$x$}};
 \draw[color=green!50!black] (0.4,-0.85) node[above] (A1) {{${x^*}$}};
\end{tikzpicture} 
- \theta
\begin{tikzpicture}[very thick,scale=1,color=blue!50!black, baseline]
  \draw (0.7,-1.8) -- (0.7,1.3); 
    \draw[color=green!50!black] (0.7,0.7) --  (0.7-1,0.7-1); 
\draw[color=white, line width=4pt] (0,-1.8) -- (0,1.3);   
  \draw (0,-1.8) -- (0,1.3); 
  \draw[color=white, line width=4pt] (0.7,0.15) --  (0.7-1,0.15-1); 
    \draw[color=green!50!black] (0.7,0.15) --  (0.7-1,0.15-1); 
   \draw[color=green!50!black]  (0.7-1,0.7-1)  .. controls +(-0.5,-0.5) and +(-0.5,-0.5) .. (0.7-1,0.15-1);
\fill[color=blue!50!black] (0.7,0.15) circle (2.9pt) node[right] (meet2) {{$g_2$}};
\fill[color=blue!50!black] (0.7,0.7) circle (2.9pt) node[right] (meet2) {{$g_1$}};
%
 \draw (0,-2.5) node[above] (X) {{$a$}};
 \draw (0.7,-2.5) node[above] (X) {{$b$}};
 \draw[color=green!50!black] (0.3,0.25) node[above] (A1)  {{${x^*}$}};
 \draw[color=green!50!black] (0.4,-0.8) node[above] (A1) {{$x$}};
\end{tikzpicture} 
\\
    &=\frac{1}{\lambda} \left(
\begin{tikzpicture}[very thick,scale=1,color=blue!50!black, baseline]
  \draw (0.7,-1.8) -- (0.7,1.3); 
   \draw (0,-1.8) -- (0,1.3); 
  \draw[color=white, line width=4pt]  (0.7,0.7) --  (0.7-1,0.7-1); 
  \draw[color=green!50!black] (0.7,0.7) --  (0.7-1,0.7-1); 
  \draw[color=white, line width=4pt]  (0.7,0.15) --  (0.7-1,0.15-1); 
    \draw[color=green!50!black] (0.7,0.15) --  (0.7-1,0.15-1); 
   \draw[color=green!50!black]  (0.7-1,0.7-1)  .. controls +(-0.5,-0.5) and +(-0.5,-0.5) .. (0.7-1,0.15-1);
%
\fill[color=blue!50!black] (0.7,0.15) circle (2.9pt) node[right] (meet2) {{$g_1$}};
\fill[color=blue!50!black] (0.7,0.7) circle (2.9pt) node[right] (meet2) {{$g_2$}};
%
 \draw (0,-2.5) node[above] (X) {{$a$}};
 \draw (0.7,-2.5) node[above] (X) {{$b$}};
 \draw[color=green!50!black] (0.3,0.25) node[above] (A1) {{$x$}};
 \draw[color=green!50!black] (0.4,-0.96) node[above] (A1) {{${x^*}$}};
\end{tikzpicture} 
- \theta
\begin{tikzpicture}[very thick,scale=1,color=blue!50!black, baseline]
  \draw (0.7,-1.8) -- (0.7,1.3); 
   \draw (0,-1.8) -- (0,1.3); 
  (0,-1.8) -- (0,1.3);   
   \draw[color=white, line width=4pt]  (0.7,0.7) --  (0.7-1,0.7-1); 
   \draw[color=green!50!black] (0.7,0.7) --  (0.7-1,0.7-1); 
    \draw[color=white, line width=4pt] (0.7,0.15) --  (0.7-1,0.15-1); 
    \draw[color=green!50!black] (0.7,0.15) --  (0.7-1,0.15-1); 
   \draw[color=green!50!black]  (0.7-1,0.7-1)  .. controls +(-0.5,-0.5) and +(-0.5,-0.5) .. (0.7-1,0.15-1);
%
\fill[color=blue!50!black] (0.7,0.15) circle (2.9pt) node[right] (meet2) {{$g_2$}};
\fill[color=blue!50!black] (0.7,0.7) circle (2.9pt) node[right] (meet2) {{$g_1$}};
%
 \draw (0,-2.5) node[above] (X) {{$a$}};
 \draw (0.7,-2.5) node[above] (X) {{$b$}};
 \draw[color=green!50!black] (0.3,0.25) node[above] (A1)  {{${x^*}$}};
 \draw[color=green!50!black] (0.4,-0.8) node[above] (A1) {{$x$}};
\end{tikzpicture} 
                               \right)=0,
  \end{align*}
    by definition of the half-braiding of $\cc{x^{*}}$ in the second step and using in the third step that  $a \in \PCent_{x}$, which
    implies
    \begin{equation}
  \label{eq:braiding-proj-cent} 
\begin{tikzpicture}[very thick,scale=1,color=blue!50!black, baseline]
\draw[color=blue!50!black] (-0.5,-0.5) node[below] (A1) {{$x^*$}};
\draw[color=blue!50!black] (0.5,-0.5) node[below] (A2) {{$a$}};
\draw[color=blue!50!black] (-0.5,0.5) node[above] (B1) {}; 
\draw[color=blue!50!black] (0.5,0.5) node[above] (B2) {};
\draw[color=blue!50!black] (A2) -- (B1);
\draw[color=white, line width=4pt] (A1) -- (B2);
\draw[color=blue!50!black] (A1) -- (B2);
\end{tikzpicture} 
\, , \quad
=   \frac{1}{\lambda}
\begin{tikzpicture}[very thick,scale=1,color=blue!50!black, baseline]
\draw[color=blue!50!black] (-0.5,-0.5) node[below] (A1) {{$x^*$}};
\draw[color=blue!50!black] (0.5,-0.5) node[below] (A2) {{$a$}};
\draw[color=blue!50!black] (-0.5,0.5) node[above] (B1) {}; 
\draw[color=blue!50!black] (0.5,0.5) node[above] (B2) {};

\draw[color=blue!50!black] (A1) -- (B2);
\draw[color=white, line width=4pt] (A2) -- (B1);
\draw[color=blue!50!black] (A2) -- (B1);
\end{tikzpicture} 
\, .   
\end{equation}
    for a $\lambda \in \Bbbk^{\times}$. This concludes the first part.
  
    For the second part, it follows analogously that for $(a,f) \in \mod(\CentQ_{z};\varphi)$ and $(c,h) \in \mod(Q_{z}, \varphi)$ also $f \qoti h \in \mod(Q_{z}, \varphi)$, using $a \in \PCent_{x}$. 
  \end{proof}

  For simple $x \in \Cat{C}$, the category  $\mod(\CentQ_{Dx};\varphi)$ is by Lemma \ref{lemma:modules-twisted-prepro} a category of modules over a $\theta$-twisted preprojective algebra by Definition \ref{definition:cl-prepro} \refitem{item:theta-prepro}:
  \begin{corollary}
    \label{corollary:q-preproj}
    Let $x \in \cent(\Cat{C})$ be a simple object with twist $\theta \in \Bbbk^{\times}$ and let $\CentQ_{Dx}$ be its double centralizer quiver. Let $Q=x \otimes -: \PCent_{x} \rightarrow \PCent_{x}$
    be the quiver given by left multiplication with $x$. 
    The rigid monoidal category     $\mod(\CentQ_{Dx};\varphi)$ is equivalent to the category $\mod(\Pi^{\theta}(Q))$ of  modules over the $\theta$-twisted preprojective algebra of $Q$, which is thereby rigid monoidal.  
   \end{corollary}
   \begin{example}
  For a symmetric fusion category $\Cat{C}$,  the centralizer of any object $x \in \cent(\Cat{C})$ is all of $\Cat{C}$.  If we equip $\Cat{C}$ with the canonical pivotal structure, the  twists are the identities and thus, the $q$-preprojective algebras are the usual preprojective algebras whose categories of modules are by Theorem \ref{theorem:double-cent-rigid-mon} rigid monoidal.
   \end{example}

\appendix
  
\section{Categorical notions}
\label{sec:categorical-notions}

For  the  convenience of  the  reader we collect some pertinent notions and  results.
\subsection{Linear categories}

All categories and functors in this article are  $\Bbbk$-linear and abelian, $\Hom$-spaces are finite dimensional $\Bbbk$-vector spaces and composition is bilinear, without further mentioning.
Recall that a category $\Cat{M}$ is called \emph{semisimple}, if each object is isomorphic to a
direct sum of simple objects, and it is called \emph{finite}, if there are only finitely many isomorphism classes of simple objects in $\Cat{M}$. 
The cardinality of these isomorphism classes is called the \emph{rank of the category}. A finite semisimple category $\Cat{M}$ of rank $n$
is equivalent to the category of $I$-graded finite-dimensional vector spaces
$\vect^I$ upon choosing representatives   $\{x_i\}_{i \in I}$ for the isomorphism classes of simple objects. We call $I$  a \emph{labeling set for $\Cat{M}$} and assume in the sequel that we have chosen a labeling set for $\Cat{M}$. 
For an object $m \in \Cat{M}$ and a vector space $V \in \vect$, there exists an object $V \otik m \in \Cat{M}$, determined up to unique isomorphism as an object representing the functor
\begin{equation}
  \label{eq:repres-vect}
  \op{M} \rightarrow \vect, \quad n \mapsto V \otik \Hom_{\Cat{M}}(n,m), 
\end{equation}
i.e. there is a natural isomorphism $V \otik \Hom_{\Cat{M}}(n,m) \cong \Hom_{\Cat{M}}(n, V \otik m)$.
For each object $m \in \Cat{M}$ there exist vector spaces $(V_{i})_{i \in I}$ and an isomorphism $m \cong  \oplus_{i \in I} V_i \otimes x_i$. 

We collect some more useful facts on (semisimple) categories.
\begin{lemma}
  \label{lemma:facts-ses}
  Let $\Cat{M}$ be a linear finite semisimple category with labeling set $I$ for representatives $(x_{i})_{i \in I}$
  of the simple objects.
  \begin{lemmalist}
  \item For all objects $m,n \in \Cat{M}$, there is a canonical natural isomorphism
    \begin{equation}
      \label{eq:ses-otik}
\oplus_{i \in I}      \Hom_{\Cat{M}}(m,x_{i}) \otik \Hom_{\Cat{M}}(x_{i},n) \cong \Hom_{\Cat{M}}(m,n)
\end{equation}
\item \label{item:coend-ses} For all objects $m \in \Cat{M}$, there are canonical natural isomorphisms
  \begin{equation}
    \label{eq:60}
    \oplus_{i \in I} \Hom_{\Cat{M}}(x_{i},m)\otik x_{i} \cong m, \quad \text{and} \quad
      \oplus_{i \in I} \Hom_{\Cat{M}}(m,x_{i})^{*}\otik x_{i} \cong m.
  \end{equation} 
  \end{lemmalist}
\end{lemma}

In addition, using that $\Bbbk$ is algebraically closed,  for all $i \in I$ there is a  non-degenerate pairing
\begin{equation}
  \label{eq:pairing-ses}
  \Hom_{\Cat{M}}(m,x_{i}) \times \Hom_{\Cat{M}}(x_{i},m) \rightarrow \Bbbk, \quad (f,g) \mapsto <f,g>,
\end{equation}
with $<f,g>$ the unique element in $\Bbbk$ such that $f \circ g = <f,g> \cdot \id_{x_{i}}$. 
This defines  the natural isomorphism
\begin{equation}
  \label{eq:nat-iso-dual}
  \Hom_{\Cat{M}}(x_{i},m)^{*} \cong \Hom_{\Cat{M}}(m,x_{i}). 
\end{equation}

For linear categories $\Cat{A}$ and $\Cat{B}$, the category of linear functors and linear natural transformation is denoted $\Fun(\Cat{A}, \Cat{B})$,
with  $\End(\Cat{A})=\Fun(\Cat{A}, \Cat{A})$ the category of endofunctors.

Functors and natural transformations involving semisimple categories  are determined by their  value on simple objects in the following sense, the proofs are standard:
\begin{lemma}
  \label{lemm:ses-cat-fun}
  Let $\Cat{A}$ be a semisimple linear category with representatives of simple objects $(x_{i})_{i \in I}$ and let $\Cat{B}$ be a linear category.
  \begin{lemmalist}
  \item  A functor $F: \Cat{A} \rightarrow \Cat{B}$ is uniquely determined up to isomorphism by the objects $F(x_{i})_{i \in I}$.
    Conversely, for any choice  of objects  $F(x_{i})_{i \in I}$ there exists a corresponding linear functor $F: \Cat{A} \rightarrow \Cat{B}$.
  \item Let $F,G: \Cat{A} \rightarrow \Cat{B}$ be linear functors. A natural transformation $\eta: F \rightarrow G$ is uniquely determined by the morphisms
    $(\eta_{x_{i}}: F(x_{i}) \rightarrow G(x_{i}) )_{i \in I}$ and for any such collection there exists a corresponding natural transformation.
  \end{lemmalist}
\end{lemma}

Finally, we recall the Deligne product of linear categories $\Cat{M}$ and $\Cat{N}$: The linear category $\Cat{M} \boti \Cat{N}$ is defined by a universal property with respect to
right exact bilinear functor, \cite{DelCat}. However, in case that $M$ and $N$ are finitely semisimple with labeling sets $I_{\Cat{M}}$ and $I_{\Cat{N}}$, the category $\Cat{M} \boti \Cat{N}$
is again finitely semisimple with labeling set $I_{\Cat{M}} \times I_{\Cat{N}}$.

We use the following graphical calculus for categories: An object $m \in \Cat{M}$ is represented as labeled vertical line
$
\begin{tikzpicture}[very thick,scale=1,color=blue!50!black, baseline]
\draw (0,-0.7) -- (0,0.7); 
%
\draw[color=blue!50!black] (-0.1,0) node[right] (A1) {{$m$}};
\end{tikzpicture} 
\, , 
$
a morphism $f: m \rightarrow n$ as a box or a dot
$
\begin{tikzpicture}[very thick,scale=1,color=blue!50!black, baseline]
\draw (0,-0.7) -- (0,0.7); 
%
%
\fill[color=blue!50!black] (0,0) circle (2.9pt) node[right] (v) {{$f$}};
\draw[color=blue!50!black] (-0.1,-0.4) node[right] (A1) {{$m$}};
\draw[color=blue!50!black] (-0.1,0.4) node[right] (A1) {{$n$}};
\end{tikzpicture} 
\, . 
$

A functor $F: \Cat{M} \rightarrow \Cat{N}$ is represented as a green line, such that
$
\begin{tikzpicture}[very thick,scale=1,color=blue!50!black, baseline]
  \draw (0,-0.7) -- (0,0.7);
  \draw[color=green!50!black] (-0.5,-0.7) -- (-0.5,0.7);
%
  \draw[color=blue!50!black] (-0.1,0) node[right] (A1) {{$m$}};
  \draw[color=green!50!black] (-0.6,0) node[right] (A1) {{$F$}};
\end{tikzpicture} 
\,  
$
represents the object $F(m) \in \Cat{N}$ and
$
\begin{tikzpicture}[very thick,scale=1,color=blue!50!black, baseline]
\draw (0,-1) -- (0,1); 
\draw[color=green!50!black] (-1,-1) .. controls +(0,0.5) and +(-0.5,-0.5) .. (0,0.25);
\fill[color=blue!50!black] (0,0.25) circle (2.9pt) node[right] (meet2) {{$f$}};
\draw[color=blue!50!black] (-0.1,-0.5) node[right] (A1) {{$m$}};
  \draw[color=blue!50!black] (-0.1,0.7) node[right] (A1) {{$n$}};
    \draw[color=green!50!black] (-0.8,-0.5) node[right] (A1) {{$F$}};
  \end{tikzpicture}
  $
a morphism $f: F(m) \rightarrow n$ in $\Cat{N}$. 

\subsection{Monoidal categories  and  their representations}
\label{sec:mono-categ-their}

We refer to \cite{EGNObook} for the definition of a monoidal category $\Cat{A}$ with unit $\unit \in \Cat{A}$. 
In the setting  of linear categories we demand that the monoidal structure
$\otimes: \Cat{A} \times \Cat{A} \rightarrow \Cat{A}$ is bilinear.
The monoidal product of two morphisms $f: x \rightarrow y$ and $g: a \rightarrow b$ is graphically represented as 
\begin{equation}
f \otimes g=
  \begin{tikzpicture}[very thick,scale=1,color=blue!50!black, baseline]

    \draw[color=blue!50!black] (-1,1)-- (-1,-1); 
    \draw[color=blue!50!black] (0.3,1)-- (0.3,-1); 
\draw (-1,0) node[minimum height=0.4cm,minimum width=0.7cm,draw,fill=white] {{$f$}};
\draw (0.3,0) node[minimum height=0.4cm,minimum width=0.7cm,draw,fill=white] {{$g$}};
\draw[color=blue!50!black] (-0.9,-0.7) node[left] (A1) {{$x$}};
\draw[color=blue!50!black] (-0.9,0.7) node[left] (A1) {{$y$}};
\draw[color=blue!50!black] (0.2,-0.7) node[right] (A1) {{$a$}};
    \draw[color=blue!50!black] (0.2,0.7) node[right] (A1) {{$b$}};
\end{tikzpicture} \; .
\end{equation}

\begin{definition}
  A monoidal category  $\Cat{A}$ is called \emph{rigid}, if every object $x \in \Cat{A}$ has both a left and right dual object: A \emph{right dual $x^{*}$} of $x$ is an object of $\Cat{A}$ together with evaluation and coevaluation morphisms, $\ev{x}: x^{*} \otimes x \rightarrow \unit$ and $\coev{x}: x \otimes x^{*} \rightarrow \unit$, which satisfy the usual triangle identities \cite[Def. 2.10.1]{EGNObook}.
  A left dual object ${}^{*}x$ is defined analogously. 
\end{definition}
Graphically, we represent these structures as
\begin{equation}
  \label{eq:ev-coev}
  \coev{x}=
  \begin{tikzpicture}[very thick,scale=1,color=blue!50!black, baseline=-0.5cm]
    \draw[color=blue!50!black] (-0.5,-0.5) .. controls +(0,-0.5) and +(0,-0.5) .. (-1.2,-0.5);
    \draw[color=blue!50!black] (-0.5,-0.5)  -- (-0.5,0);
    \draw[color=blue!50!black] (-1.2,-0.5)  -- (-1.2,0);
    %
    %
    %
    \draw[color=blue!50!black] (-0.24,-0.5) node[above] (A1) {{$x^{*}$}};
\draw[color=blue!50!black] (-1.4,-0.5) node[above] (A1) {{$x$}};
\end{tikzpicture}, \quad
\ev{x}=
  \begin{tikzpicture}[very thick,scale=1,color=blue!50!black, baseline=-0.7cm]
    \draw[color=blue!50!black] (-0.5,-0.5) .. controls +(0,0.5) and +(0,0.5) .. (-1.2,-0.5);
   \draw[color=blue!50!black] (-0.5,-0.5)  -- (-0.5,-1);
    \draw[color=blue!50!black] (-1.2,-0.5)  -- (-1.2,-1);
    %
    %
    %
    \draw[color=blue!50!black] (-0.3,-0.8) node[above] (A1) {{$x$}};
\draw[color=blue!50!black] (-1.5,-0.8) node[above] (A1) {{$x^{*}$}};
\end{tikzpicture}\,
\end{equation}
with triangle identities,
\begin{equation}
  \label{eq:triangle}
  \begin{tikzpicture}[very thick,scale=1,color=blue!50!black, baseline=-0.5cm,xscale=-1]
    \draw[color=blue!50!black] (-0.5,-0.5) .. controls +(0,-0.5) and +(0,-0.5) .. (-1,-0.5);
    \draw[color=blue!50!black] (-0.5,-0.5)  -- (-0.5,0);
    \draw[color=blue!50!black] (-1,-0.5)  -- (-1,0);
    \draw[color=blue!50!black] (-1,0) .. controls +(0,0.5) and +(0,0.5) .. (-1.5,0);
     \draw[color=blue!50!black] (-0.5,0)  -- (-0.5,0.7);
    \draw[color=blue!50!black] (-1.5,0)  -- (-1.5,-1);
    %
    %
    %
    \draw[color=blue!50!black] (-0.4,0.5) node[right] (A1) {{$x$}};
        \draw[color=blue!50!black] (-1.1,-0.7) node[right] (A1) {{$x$}};
\end{tikzpicture}
=
  \begin{tikzpicture}[very thick,scale=1,color=blue!50!black, baseline=-0.5cm]
    \draw[color=blue!50!black] (0,-1)  -- (0,0.7);
    \draw[color=blue!50!black] (-0.1,-0.2) node[right] (A1) {{$x$}};
  \end{tikzpicture}\, , \quad
  \begin{tikzpicture}[very thick,scale=1,color=blue!50!black, baseline=-0.5cm]
    \draw[color=blue!50!black] (-0.5,-0.5) .. controls +(0,-0.5) and +(0,-0.5) .. (-1,-0.5);
    \draw[color=blue!50!black] (-0.5,-0.5)  -- (-0.5,0);
    \draw[color=blue!50!black] (-1,-0.5)  -- (-1,0);
    \draw[color=blue!50!black] (-1,0) .. controls +(0,0.5) and +(0,0.5) .. (-1.5,0);
     \draw[color=blue!50!black] (-0.5,0)  -- (-0.5,0.7);
    \draw[color=blue!50!black] (-1.5,0)  -- (-1.5,-1);
    %
    %
    %
    \draw[color=blue!50!black] (-0.24,-0.5) node[above] (A1) {{$x^{*}$}};
\draw[color=blue!50!black] (-1.6,-0.7) node[right] (A1) {{$x^{*}$}};
\end{tikzpicture}
=
  \begin{tikzpicture}[very thick,scale=1,color=blue!50!black, baseline=-0.5cm]
    \draw[color=blue!50!black] (0,-1)  -- (0,0.7);
    \draw[color=blue!50!black] (-0.1,-0.2) node[right] (A1) {{$x^{*}$}};
\end{tikzpicture}\, .
\end{equation}
Dual objects are unique up to unique isomorphism and provide adjunctions
\begin{equation}
  \label{eq:5}
  \Hom_{\Cat{A}}(y \otimes x, z) \cong \Hom_{\Cat{A}}(y, z \otimes x^{*}) \quad \text{and} \quad
  \Hom_{\Cat{A}}( x \otimes y, z) \cong \Hom_{\Cat{A}}(y, {}^{*}x \otimes z).
\end{equation}
It follows, that the functors $x \otimes -: \Cat{A} \rightarrow \Cat{A}$ and $- \otimes x: \Cat{A} \rightarrow \Cat{A}$ have both adjoints and are thus exact functors. 

The fact that the category $\Cat{A}$ possesses a rigid monoidal structure has further homological consequences. The proof of  \cite[Prop.2.3]{FinTen} in case of a finite-dimensional algebra generalizes  directly:
\begin{proposition}
  \label{proposition:self-inj}
  Let $\Cat{A}=\mod(A)$ be the category of finite-dimensional modules over a (not necessarily finite-dimensional)
  $\Bbbk$-algebra $A$. Suppose there exists a rigid monoidal structure on $\Cat{A}$. Then
  all finite-dimensional projective modules are injective and vice versa.
\end{proposition}

By a \emph{multi-fusion category} we mean a finite semisimple rigid monoidal category $\Cat{A}$. It is called a \emph{fusion category}, if the monoidal unit is simple.
For a multi-fusion category $\Cat{A}$, the \emph{Grothendieck ring} $\Gr(\Cat{A})$ is the free $\Z$-module on  the isomorphism classes of objects of $\Cat{A}$. It acquires  the structure of a based ring, see \cite[Sec. 4.9]{EGNObook} with multiplication written $ab$ for $a, b \in \Gr(\Cat{A})$. 

\begin{example}
  \label{example:vectG}
  For a finite group $G$ the finite-dimensional $G$-graded vector spaces $(V_{g})_{g \in G}$ form a fusion category $\vect_{G}$ with  monoidal structure given by $(V \otimes W)_{g}=\oplus_{ab=g}V_{a}\otimes W_{b}$. As representative for the simple objects we consider the objects $e_{g}$ with as only non-zero vector space $\Bbbk$ in degree $g \in G$.
\end{example}

Let  $\Cat{A}$ be a rigid monoidal category. A \emph{pivotal structure} on $\Cat{A}$ is a monoidal natural isomorphism $\gamma_{x}: x \stackrel{\cong}{\longrightarrow} x^{**}$ for $x \in \Cat{A}$. Given a pivotal structure, the left trace of an endomorphism $f \in \End_{\Cat{A}}(x)$ is
\begin{equation}
  \label{eq:51}
  \tr(f)=
  \begin{tikzpicture}[very thick,scale=1,color=blue!50!black, baseline]
    \draw[color=blue!50!black] (-0.5,-0.5) .. controls +(0,-0.5) and +(0,-0.5) .. (-1.2,-0.5);
    \draw[color=blue!50!black] (-0.5,-0.5)  -- (-0.5,0.6);
    \draw[color=blue!50!black] (-1.2,-0.5)  -- (-1.2,0.6);
        \draw[color=blue!50!black] (-0.5,0.6) .. controls +(0,0.5) and +(0,0.5) .. (-1.2,0.6);

    %
    %
    %
    \draw[color=blue!50!black] (-0.3,-0.7) node[above] (A1) {{$x$}};
        \fill[color=blue!50!black]   (-0.5,0) circle (2.9pt) node[right] (meet2) {{$f$}};
\draw[color=blue!50!black] (-1.4,-0.7) node[above] (A1) {{$x^{*}$}};
  \end{tikzpicture}.
\end{equation}
It induces a cyclic non-degenerate pairing of the $\Hom$-spaces, if $\Cat{A}$ is semisimple. In particular,
the dimensions $\dim(x)=\tr(\id_{x})$ of the simple objects $x \in \Cat{A}$ are non-zero.
In the sequel, the pivotal structure is used to identify the left and right duals and is suppressed from the notation and the  diagrams.

An object $x \in \Cat{A}$ is called \emph{self-dual}, if there is an isomorphism $\phi: x^{*}  \rightarrow x$ in $\Cat{A}$. For a simple self-dual object $x$, there is a unique number $\FS(x) \in \{ \pm 1\}$, called the \emph{Frobenius-Schur indicator of $x$}, such that

\begin{equation}
  \label{eq:FS-def}
  \begin{tikzpicture}[very thick,scale=1,color=blue!50!black, baseline]
    \draw[color=blue!50!black] (-0.5,-0.5) .. controls +(0,-0.5) and +(0,-0.5) .. (-1.2,-0.5);
    \draw[color=blue!50!black] (-0.5,-0.5)  -- (-0.5,0.5);
    \draw[color=blue!50!black] (-1.2,-0.5)  -- (-1.2,0.5);
    %
    %
    \fill[color=blue!50!black] (-0.5,0.1) circle (2.9pt) node[right] (meet2) {{$\phi$}};
    \draw[color=blue!50!black] (-0.24,-0.6) node[above] (A1) {{$x^{*}$}};
        \draw[color=blue!50!black] (-0.24,0.15) node[above] (A1) {{$x$}};
\draw[color=blue!50!black] (-1.4,-0.6) node[above] (A1) {{$x$}};
  \end{tikzpicture}
  = \FS(x)
  \begin{tikzpicture}[very thick,scale=1,color=blue!50!black, baseline]
    \draw[color=blue!50!black] (-0.5,-0.5) .. controls +(0,-0.5) and +(0,-0.5) .. (-1.2,-0.5);
    \draw[color=blue!50!black] (-0.5,-0.5)  -- (-0.5,0.5);
    \draw[color=blue!50!black] (-1.2,-0.5)  -- (-1.2,0.5);
    %
    %
    \fill[color=blue!50!black] (-1.2,0.1) circle (2.9pt) node[right] (meet2) {{$\phi$}};
    \draw[color=blue!50!black] (-0.24,-0.65) node[above] (A1) {{$x$}};
    \draw[color=blue!50!black] (-1.45,-0.65) node[above] (A1) {{$x^{*}$}};
    \draw[color=blue!50!black] (-1.4,0.15) node[above] (A1) {{$x$}};
  \end{tikzpicture}
  \; .
\end{equation}
 
 Let $\Cat{A}, \Cat{B}$ be monoidal categories.
    A functor $F: \Cat{A} \rightarrow \Cat{B}$ is
    called \emph{lax monoidal}, if it is equipped with natural morphisms
    \begin{equation}
      \label{eq:lax-mon}
f_{a,b}:      F(a) \otimes F(b) \rightarrow F(a \otimes b),
\end{equation}
which are coherent with respect to the monoidal product of three objects and the unit.
It is called an \emph{oplax monoidal} functor, if it is equipped with coherent natural morphisms
  \begin{equation}
      \label{eq:lax-mon}
f_{a,b}:     F(a \otimes b) \rightarrow   F(a) \otimes F(b). 
\end{equation}
It is called \emph{strong monoidal}, if it is (op)lax monoidal such that the coherence morphisms are invertible.

\paragraph{Braided monoidal categories}

A \emph{braiding} on a monoidal category $\Cat{C}$ is a collection of natural isomorphisms
$c_{z,x}:z \otimes x \rightarrow x \otimes z$ for all objects $x,z$, that satisfy a hexagon diagram.
The coherence ensures, that the graphical expression
\begin{equation}
  \label{eq:braiding-def}
  c_{z,x}=
  \begin{tikzpicture}[very thick,scale=1,color=blue!50!black, baseline]
\draw[color=blue!50!black] (-0.5,-0.5) node[below] (A1) {{$z$}};
\draw[color=blue!50!black] (0.5,-0.5) node[below] (A2) {{$x$}};
\draw[color=blue!50!black] (-0.5,0.5) node[above] (B1) {}; 
\draw[color=blue!50!black] (0.5,0.5) node[above] (B2) {};
\draw[color=blue!50!black] (A2) -- (B1);
\draw[color=white, line width=4pt] (A1) -- (B2);
\draw[color=blue!50!black] (A1) -- (B2);
\end{tikzpicture} 
\end{equation}
for $c_{z,x}$ is meaningful. 

In a braided pivotal category, the \emph{twist} on an object $x$ is defined as
\begin{equation}
  \label{eq:twist-def}
  \begin{tikzpicture}[very thick,scale=1,color=blue!50!black, baseline=1cm]
\draw[color=blue!50!black] (-1.4,1.1) .. controls +(0,-0.2) and +(0,-0.2) .. (-1.7,1.1);
\draw[color=blue!50!black] (-1.4,1.1) .. controls +(0,0.25) and +(0,0.25) .. (-1.9,1.1);
%
  \draw[color=blue!50!black]  (-1.9,1.1) -- (-1.9,0.2); 
    \draw[color=white, line width=4pt]  (-1.7,1.1) -- (-1.7,1.8) ; 
    \draw[color=blue!50!black]    (-1.7,1.1) -- (-1.7,1.8) ; 
\draw[color=blue!50!black] (-1.9,1.3) node[above] (A1) {{$x$}};
\end{tikzpicture}
=
  \begin{tikzpicture}[very thick,scale=1,color=blue!50!black, baseline=1cm]
    \draw[color=blue!50!black]    (-1.7,0.2) -- (-1.7,1.8) ; 
    \draw[color=blue!50!black] (-1.9,1.3) node[above] (A1) {{$x$}};
    \fill[color=blue!50!black]  (-1.7,1) circle (2.9pt) node[right] (meet2) {{$\theta_x$}};
    \end{tikzpicture}
  \end{equation}
  In case that $x$ is simple, the twist is a scalar $\theta_{x} \in \Bbbk^{\times}$ times the identity on $x$. 

For a monoidal category $\Cat{A}$, its \emph{Drinfeld center} $\cent(\Cat{A})$ is the category with objects pairs
$(z, c_{z})$ of an object $z \in \Cat{A}$ together with a \emph{half-braiding} $c_{z}(x): z \otimes x  \rightarrow x \otimes z$: For all $x \in \Cat{A}$ these are natural isomorphisms, which satisfy a hexagon diagram. The morphisms in $\cent(\Cat{A})$ are morphisms in $\Cat{A}$ that are compatible with the half-braidings.
The category $\cent(\Cat{A})$ is naturally a braided monoidal category, which is rigid (pivotal), if $\Cat{A}$ is, 
see \cite[Sec.8.5]{EGNObook} for more details. 

Diagrams for (braided, pivotal) monoidal categories can be used to perform calculations graphically: A diagram can be evaluated to a morphism in the category and the evaluation is invariant under suitable isotopies, see \cite{Sel} for a summary. 
\paragraph{Module categories}
Let $\Cat{A}$ be a monoidal category. 
A  category $\Cat{M}$ is called a \emph{$\Cat{A}$-left module category}, if it is equipped with a bifunctor
  \begin{equation}
    \label{eq:14}
    \act: \Cat{A} \times \Cat{M} \rightarrow \Cat{M}, 
  \end{equation}
  and natural isomorphisms $(a \otimes b)\act m \cong a \act (b \act m)$, which are coherent with respect to the action of three objects and the unit of $\Cat{A}$. 
A  module 

For  $\AM$ and $\AN$  left $\Cat{A}$-module categories, a \emph{lax module functor}
  $F: \Cat{M} \rightarrow \Cat{N}$ is a functor with coherent natural morphisms
  \begin{equation}
    \label{eq:12}
    a \act F(m) \rightarrow F(a \act m),
  \end{equation}
  for all $a \in \Cat{A}$ and $m \in \Cat{M}$.
  It is called a \emph{oplax module functor}, if it is equipped with coherent natural morphisms 
  \begin{equation}
    \label{eq:oplax-module}
  F(a \act m) \rightarrow    a \act F(m),
  \end{equation}
  for all $a \in \Cat{A}$ and $m \in \Cat{M}$. It is called a \emph{strong module functor}, if it is a lax module functor such that the coherence morphisms are invertible.
  
 For $F, G: \Cat{M} \rightarrow \Cat{N}$ two lax module functors, a \emph{module natural transformation}
  $\eta: F \rightarrow G$ is a natural transformation, such that for all $a \in \Cat{A}$ and $m \in \Cat{M}$, the diagram
  \begin{equation}
    \label{eq:16}
    \begin{tikzcd}
      a \act F(m) \ar{r}{a \act \eta_{m}} \ar{d}{} & a \act G(m) \ar{d}{}\\
      F(a \act m) \ar{r}{\eta_{a \act m}} & G(a \act m)
    \end{tikzcd}
  \end{equation}
  commutes.

  By considering for two monoidal categories $\Cat{A}$ and $\Cat{B}$ the monoidal category $\Cat{A} \times \Cat{B}^{\rev}$, where we consider  the reversed  monoidal structure on $\Cat{B}$, we obtain the  corresponding notions of $(\Cat{A},\Cat{B})$-bimodule category and bimodule functor as well as bimodule natural transformations.

  We collect these structures in the following definition, see the next subsection for a review on bicategories. 
  \begin{definition}\label{definition:bicat-mod-cat}
Let $\Cat{A}$   be  a monoidal category.  The  $\Cat{A}$-left module categories, oplax module functors and module natural transformations
  form a bicategory ${}_{\Cat{A}}\Modoplax$.  For two monoidal categories $\Cat{A},\Cat{B}$   bimodule categories and the corresponding morphisms form a bicategory
  ${}_{\Cat{A}}\BiModoplax_{\Cat{B}}$.
  \end{definition}

  In particular, for module categories $\AM$ and $\AN$, the strong module functors and module natural transformations form a category $\Fun_{\Cat{A}}(\AM,\AN)$.
  If we consider $\Cat{A}$ as a bimodule category over itself, we have a canonical equivalence of monoidal categories \cite[Prop.7.13.8]{EGNObook}:
  \begin{equation}
    \label{eq:Drinf-mod-fun}
    \End_{\Cat{A},\Cat{A}}(\AAA) \cong \cent(\Cat{A}).
  \end{equation}
  Indeed, the assignment of $ z \in \cent(\Cat{A})$ to the endofunctor $Q_z(a)= z \otimes a$ provides the claimed equivalence.

\subsection{Bicategories}
\label{sec:bicategories}

In this short subsection we fix our conventions regarding bicategories and the morphisms between them. 
For a more complete review, see  \cite{Benabou}.

A \emph{bicategory} $\Cat{B}$ consists of objects $a,b \in \Cat{B}$ and for any two objects of a category $\Cat{B}(a,b)$ of 1- and 2-morphisms. Moreover, there is a composition functor
\begin{equation}
  \label{eq:comp-bicat}
 \circ: \Cat{B}(b,c) \times \Cat{B}(a,b) \rightarrow \Cat{B}(a,c) 
\end{equation}
 and a unit $1_{a} \in \Cat{B}(a,a,)$ for all objects, which are weakly associative and unital analogous to the  case of monoidal categories. Indeed, a bicategory with one object is the same as a monoidal category. 

A \emph{2-functor} $\Gamma: \Cat{B} \rightarrow \Cat{D}$ between bicategories maps objects $a \in \Cat{B}$ to objects $\Gamma(a) \in \Cat{D}$ and the higher morphisms with functors
$\Gamma: \Cat{B}(a,b) \rightarrow \Cat{D}(\Gamma(a),\Gamma(b))$. Moreover, it has coherent isomorphism $\Gamma_{F,G}: \Gamma(F) \circ \Gamma(G) \rightarrow \Gamma(F \circ G)$ and $1_{\Gamma(a)} \cong \Gamma(1_{a})$.

A \emph{lax natural 2-transformation} $\rho: \Gamma_{1} \rightarrow \Gamma_{2}$ between 2-functors
$\Gamma_{i}: \Cat{B} \rightarrow \Cat{D}$, $i=1,2$ consists of a collection of 1-morphisms
$\rho_{a} : \Gamma_{1} (a) \rightarrow \Gamma_{2}(a)$ for all $a \in \Cat{B}$ and for all 1-morphisms $F: a \rightarrow b$ a 2-morphism $\rho_{a}: \rho_{b} \circ \Gamma_{1}(F) \rightarrow \Gamma_{2}(F) \circ \rho_{a}$, which are natural and coherent.

An \emph{oplax natural 2-transformation} $\sigma: \Gamma_{1} \rightarrow \Gamma_{2}$ between the same 2-functors has again for all objects $a \in \Cat{B}$ a 1-morphism $\sigma_{a}: \Gamma_{1}(a) \rightarrow \Gamma_{2}(a)$, but coherent natural 2-morphisms
$\sigma_{F}: \Gamma_{2}(F) \circ \rho_{a} \rightarrow \rho_{b} \circ \Gamma_{a}(F)$ for all 1-morphism $F: a \rightarrow b$.

Finally, for two lax natural 2-transformations $\rho_{1},\rho_{2}: \Gamma_{1} \rightarrow \Gamma_{2}$, a \emph{modification} $\phi: \rho_{1} \rightarrow \rho_{2}$ is a collection of 2-morphisms $\phi_{a}: (\rho_{1})_{a} \rightarrow (\rho_{2})_{a}$ for all objects $a \in \Cat{B}$, which obeys an obvious condition for all 1-morphisms $F: a \rightarrow b$.
Similarly, we define modifications between two oplax natural 2-transformations.
Between bicategories $\Cat{B}$ and $\Cat{D}$  there is a bicategory $\twoFunol(\Cat{B},\Cat{D})$
of 2-functors, oplax natural 2-transformations and modification and similar for the lax natural 2-transformations.

Examples of bicategories are the bicategory $\Lincat$ with objects linear categories, 1-morphisms linear functors and 2-morphisms linear natural transformations and the sub bicategory $\lincatses$ of finite semisimple categories.
A bicategory with one object is the same as a monoidal category.


\end{document}